\documentclass{amsart}
\usepackage{amssymb,mathrsfs}
\usepackage{color}
\usepackage{amsmath,amsfonts,amssymb}
\usepackage{graphicx,epsfig,subfigure}
\usepackage{enumitem}

\usepackage{lineno}

\usepackage[textheight=25 cm, textwidth=16cm, headheight=10pt, headsep=20pt, footskip=10pt, left=1.5cm, right=1.5cm]{geometry}

\def\bR {\mathbb{R}}

\def\fH {\mathfrak{H}}

\def\\fH {\mathcal{H}}

\def\cP {\mathcal{P}}

\def\\hbar {{\\hbarilon}}

\newcommand{\Lip}{\operatorname{Lip}}

\newcommand{\ba}{\begin{aligned}}
\newcommand{\ea}{\end{aligned}}

\newcommand{\be}{\begin{equation}}
\newcommand{\ee}{\end{equation}}

\newcommand{\bea}{\begin{eqnarray}}
\newcommand{\eea}{\end{eqnarray}}

\newcommand{\MKd}
{W_2}
\newcommand{\MKu}
{W_1}

\newtheorem{Thm}{Theorem}[section]

\renewcommand{\hbar}{{}}

\usepackage{hyperref}
\hypersetup{pdfborder=0 0 0}

\newcommand{\gen}{\gamma}
\newcommand{\vla}{\nu}

\newcommand{\rhoi}{\rho^{in}}
\newcommand{\rhoin}{(\rhoi)^{\otimes N}}

\newcommand{\macphin}{\varphi^{in}}

\newcommand{\rc}{{\tau_c}}

\newcommand{\REVISION}[1]{\textcolor{black}{#1}}
{}
\newcommand{\COMMENT}[1]
{}

\begin{document}

\title[Microscopic, kinetic and hydrodynamic hybrid models of collective motions with
chemotaxis]{\Large Microscopic, kinetic and hydrodynamic hybrid models of collective motions with
chemotaxis: a numerical study}

\author[M. Menci]{Marta Menci}\address[M.M]{Università Campus Bio-Medico di Roma
Via Àlvaro del Portillo 00128 Roma, Italy}
\email{m.menci@unicampus.it}
\author[R. Natalini]{Roberto Natalini}
\address[R.N]{IAC, Via dei Taurini, 19, 00185 Roma RM Italy}
\email{roberto.natalini@cnr.it}
\author[T. Paul]{Thierry Paul}
\address[T.P.]{CNRS  Laboratoire Ypatia des Sciences Mathématiques (LYSM),  Roma, Italia \& Sorbonne Universit\'e, CNRS, Universit\'e \ \ \ \ Paris Cit\'e, Laboratoire Jacques-Louis Lions (LJLL), F-75005 Paris, France}
\email{thierry.paul@upmc.fr}

\begin{abstract}
\Large
A general class of hybrid models has been introduced recently, gathering the advantages of multiscale descriptions. Concerning biological applications, the particular coupled structure fits to collective cell migrations and pattern formation phenomena due to intercellular and chemotactic stimuli.
In this context, cells are modelled as discrete entities and their dynamics are given by ODEs, while the chemical signal influencing the motion is considered as a continuous signal which solves a diffusive equation.
From the analytical point of view, this class of model has been recently proved to have a mean-field limit in the Wasserstein distance towards a system given by the coupling of a Vlasov-type equation with the chemoattractant equation. 
Moreover, a pressureless nonlocal Euler-type system has been derived for these models, rigorously equivalent to the Vlasov one for monokinetic initial data. 
For applications, the monokinetic assumption is quite strong and far from real experimental setting. 
The aim of this paper is to introduce a numerical approach to the hybrid coupled structure at the different scales, investigating the case of general initial data. 
Several scenarios will be presented, aiming at exploring the role of the different terms on the overall dynamics.
Finally, the pressureless nonlocal Euler-type system is generalized by means of an additional pressure term.   
\end{abstract}
\LARGE
\maketitle

\tableofcontents
\section{Introduction and rigorous results previously obtained}\label{intro}

Collective motions modelling is attracting the interest of different research fields, due to the great variety of living and non-living systems exhibiting collective behaviors (see the seminal paper in \cite{vicsek-collective} for an introduction the the field). 
Different approaches have been proposed, depending on the features models aim at reproducing and on the scale of observation. 
At the microscopic scale, the majority of models are based on elementary mechanical interactions among agents. Alignment, repulsion and attraction kind of interactions are usually considered to model the characteristic tendence to move towards the same direction, keeping cohesion in the group and avoiding collisions \cite{vicsek, CS1, CS2, dorsogna, strombom}.
Concerning alignment models, the seminal model in \cite{CS1} originally proposed to describe the dynamics in flocks of birds was extended in different directions \cite{carr4, rt, prt}.

In the biological field, collective behaviors occur in living processes involving cells dynamics, see  \cite{arboleda, belmonte, hatzikirou, vicsek-collective-cell} for seminal review papers of the field. The main feature is that collective cells migration is also driven by a chemical stimuli, and not only by mechanical interactions among agents. 
On the one hand, the microscopic approach allows to model in fine details mechanical interactions among cells. On the other hand, when modelling the evolution in time of a chemical signal influencing the overall dynamics, the microscopic approach is actually not convenient, and a continuum approach, based on reaction-diffusion equation better fullfill the requirement.
The structure of hybrid coupled models here considered gather the advantages of microscopic and macroscopic descriptions.
Focusing on alignment models, in \cite{dicostanzo} a model for the morphogenesis in the zebrafish lateral line primordium was proposed. Based on the experimental data in \cite{gilmour06, gilmour08}, the model couple alignment and attraction-repulsion kind of interactions with chemotactic effect. 
The model is hybrid in the following sense:
particles are considered as circular shaped entities, whereas the chemical signal is modelled as a continuum. During time evolution, the positions and velocities of each cells are described by second-order Ode, whereas the concentration of the chemical signal is described as the solution of a parabolic equation with term of source and degradation.
The original structure has been extended including stochasticity and/or further cell mechanism for different biological phenomena \cite{dicostanzo_cardio,axioms}. 
From the analytical point of view, a simplified version of the model in  \cite{dicostanzo} was proposed in \cite{dicostanzo2} to allow a full analytical investigation of the asymptotic behavior, together with well posedeness results in $\mathbb{R}^2$. Further analytical results on generalized version of hybrid systems can be found in \cite{mp, mp2, mppsJDE}. 
In all of the previous papers, the hybrid structure considered is at particle level. 
We begin our work summarizing the main results of the literature concerning the passage from particle to kinetic and macroscopic scales. 
%
%
Consider on $\bR^{2dN}\ni 
((x_i(t))_{i=1,\dots,N},(v_i(t))_{i=1,\dots,N}):=(X(t),V(t))$ the following vector field
\be\label{eq1}
\left\{
\begin{array}{l}
\dot {x_i}(t)=v_i,\\
\dot {v_i}(t)=F_i(t,X(t),V(t))
\end{array}
 \ 
 \right.
\ee

where 
\be\label{defG}
F_i(t,X,V)=\frac1N\sum_{j=1}^N\gamma(v_i-v_j,x_i-x_j)+\eta
\nabla_{x}\varphi^t(x_i)
+F_{ext}(x_i),
\ee
for any $ i=1,\dots,N$. 

Here  $F_{ext}$ denotes an external force, $x_i,v_i$ are the position and velocity of the $i-th$ cell and $\varphi$ stands for a generic chemical signal produced by the cells themselves.
In particular, $\varphi$ satisfies the equation
\be\label{defeqphi0}
\partial_s\varphi^s(x)=D\Delta \varphi^s-\kappa\varphi^s +f(x,X(s)),
\ s\in[0,t],
\ee 
for some $\kappa, D,\eta\geq0$ and function $f$ of the form
\be\label{deff}
f(x,X)=\frac1N\sum_{j=1}^N\chi(x-x_j),\ \ \chi\in\mathcal C_c^1.
\ee
\REVISION{The choice of function $\chi$ is justified by the fact that every cell produces the chemical signal, i.e. cells are seen as regions from which the signal arises. 
The assumption of a compact support recovers the modelling choice of \cite{axioms, dicostanzo_cardio}, where each cell $j=1,..,N$ is modeled as a disk, centered in its position $x_j$.}
The function $\gamma:\bR^d\times \bR^d\to\bR^d$ models the interactions among agents and it is supposed to be Lipschitz continuous\footnote{
{through this paper we define $\Lip(f)$ for $f:\ \bR^{n}\to\bR^{m},m,n\in\mathbb{N},$ as $\Lip(f):=
\sqrt{\sum\limits_{i=1}^m(\Lip(f_i)^2}$.}}\label{deffone}.
Initial data are given by $(X(0),V(0), \varphi(0))=(X^{in},V^{in}, \varphi^{in})$.

Let us describe the main steps and results of  \cite{rt2}, in order to introduce the corresponding kinetic and hydrodynamic limit derived for the particle hybrid structure.

For any fixed function $\macphin$ and any  $t,N$ the mapping $\Phi^t_N=\Phi^t$ is defined as 
\be\label{defphitN}
\left\{
\begin{array}{rrcl}
\Phi^t_N:&\bR^{2dN}&\longrightarrow&\bR^{2dN}\\
& (X^{in},V^{in})&\longrightarrow& (X(t),V(t))\mbox{ solution of $(\ref{eq1},\ref{defG},\ref{defeqphi0},\ref{deff})$.}
\end{array}
\right.
\ee
Note that $\Phi^t_N$ is not a  flow.

In \cite{rt2} a kinetic model was derived, corresponding to system $(\ref{eq1},\ref{defG},\ref{defeqphi0}, \ref{deff}$).
In particular, the following non local in time Vlasov system is derived:

\be\label{eq1V}
\partial_t\rho^t+v\cdot\nabla_x\rho^t= \REVISION{\nabla_v \cdot}  (
\vla(t,x,v)\rho^t),\ \rho^{0}=\rhoi
\ee
where
\be\label{defGV}
\vla(t,x,v)=\gamma*\rho^t(x,v)+\eta\nabla_x\psi^t(x)+F_{ext}(x)
\ee
and  $\psi^s$ satisfies
\be\label{defeqphiV}
\partial_s\psi^s(x)=D\Delta_x\psi\REVISION{^s}-\kappa\psi^s+g(x,\rho^s),\ \psi^0=\macphin
\ee
with
\be\label{defgvlasov}
g(x,\rho^s)=
\int_{\bR^{2d}}\chi(x-y)\rho^s(y,\xi)dyd\xi.
\ee

The hydrodynamic limit of Cucker-Smale models has provided up to now a large literature, whose exhaustive quotation is beyond the scope of the present paper. 
The approach and results in \cite{rt2} differ from previous results of the literature. In particular authors do not make use of empirical measures formalism. 

Indeed one easily sees that, in the case where $\rho^{in}$ is monokinetic, i.e.
\begin{equation}
\rho^{in}(x,v)=\mu^{in}(x)\delta(v-u^{in}(x)),
\end{equation}
 the monokinetic form is preserved by the Vlasov equation \eqref{eq1V} and  the solution is furnished by the solution of the 
Euler type system:
\be\label{systhydro}
\left\{
\begin{array}{l}
\partial_t\mu^t+\REVISION{\nabla_x \cdot} (u^t \mu^t) = 0\\
\partial_t(\mu^tu^t)+
\REVISION{\nabla_x \cdot} (\mu^t(u^t)^{\otimes 2}) = 
\mu^t\int\gamma(\cdot-y,u^t(\cdot)-u^t(y))\mu^t(y)dy+
\eta\mu^t\nabla_{x}\psi^t
+\mu^tF_{ext}
\end{array}
\right.
\ee
where
\begin{equation}
\partial_s\psi^s = D\Delta_x\psi\REVISION{^s}-\kappa\psi\REVISION{^s} +
\chi*\mu^s,\ s\in[0,t],
\ \Psi^0=\macphin.
\end{equation}

Let the Wasserstein distance of order two between two probability measures $\REVISION{\tilde \mu, \tilde \nu}$ on $\bR^m$ with finite second moments be defined as
\begin{equation}
\MKd(\REVISION{\tilde \mu, \tilde \nu})^2
=
\inf_{\gamma\in\Gamma(\REVISION{\tilde \mu, \tilde \nu})}\int_{\bR^m\times\bR^m}
|x-y|^2\gamma(dx,dy)
\end{equation}
where $\Gamma(\REVISION{\tilde \mu, \tilde \nu})$ is the set of probability measures on $\bR^m\times\bR^m$ whose marginals on the two factors are $\REVISION{\tilde \mu}$ and $\REVISION{\tilde \mu, \tilde \nu}$ (see \cite{VillaniAMS,VillaniTOT}).
For the sake of completeness we here sketch the two main results of \cite{rt2}:
\begin{Thm}\label{mainthm}

Let $\rho^{in}$ be a compactly supported 
probability on $\bR^{2dN}$,  let $\Phi^t_N$ be the mapping generated by the particles system $(\ref{eq1},\ref{defG},\ref{defeqphi0},\ref{deff})$ as defined by \eqref{defphitN},
 and let 
 $\tau_{\rho^{in}}$ be the function defined in \cite[Formula (41), \REVISION{Proposition 5.1}]{rt2}.
 
Then, for any $t\geq 0$,
$$
\MKd\left((\Phi^t_N\#(\rho^{in})^{\otimes N})_{N;1},\rho^t\right)^2\leq 
\tau_{\rho^{in}}(t)
\left\{
\begin{array}{ll}
N^{-\frac12}&d=1\\
N^{-\frac12}{\log{N}}&d=2\\
N^{-\frac1{d}}&d>2
\end{array}
\right.
$$
where $\rho^t$ is the solution of the   Vlasov equation 
$(\ref{eq1V},\ref{defGV},
\ref{defeqphiV},\ref{defgvlasov})$ with initial condition $\rho^{in}$  provided by \cite[Theorem 8.1]{rt2}
\footnote{\REVISION{We recall that the pushforward of a measure $\tilde\mu$ by a measurable function $\Phi$ is $\Phi\# \tilde \mu$ defined by $\int h d(\Phi\# \tilde\mu):=\int (h \circ \Phi)d\tilde\mu$ for every test function $h$.
Moreover, $(\rho^{in})^{\otimes N}$ stands for $\underbrace{\rho^{in} \otimes ...  \otimes \rho^{in}}_\text{N times}$}. }.
Moreover, let us 
denote by $\varphi^t_{Z^{in}}$ the chemical density solution of $(\ref{eq1},\ref{defG},\ref{defeqphi0},\ref{deff})$ with initial data $(Z^{in}=(X^{in},V^{in}),\macphin)$ and by $\psi^t_{\rhoi}$ the one solution of $(\ref{eq1V},\ref{defGV},
\ref{defeqphiV},\ref{defgvlasov})$ with initial data $(\rhoi,\macphin)$.

Then
$$
\int_{\bR^{2dN}}\|\nabla\varphi^t_{Z^{in}}-\nabla\psi^t_{\rhoi}\|^2_\infty
(\rhoi)^{\otimes N}(dZ^{in})
\leq 
\rc(t)
\left\{
\begin{array}{ll}
N^{-\frac12}&d=1\\
N^{-\frac12}{\log{N}}&d=2\\
N^{-\frac1{d}}&d>2
\end{array}
\right.
$$
where $\rc$ is defined below by \cite[Formula (54)]{rt2}.

Finally, the functions $\tau(t),\rc(t)$ depend only on $t$, $\Lip(\gen),\Lip(\chi),\Lip(\nabla\chi)$, and the supports of $\Phi^t_N\#\rhoin$ and $\rho^t$, and satisfies the following estimate for all $t\in\bR$,
$$
\tau_{\rhoi}(t)\leq e^{e^{Ct}},\ \rc(t)\leq e^{e^{C_ct}}
$$
for some constants $C,C_c$,depending  on  $\Lip(\gen),\Lip(\chi),\Lip(\nabla \chi)$ and $|supp(\rhoi)|$.
\end{Thm}
\begin{Thm}\label{hydro}
Let $\mu^t,u^t,\psi^t$ be a solution to the following system
\be\label{systhydro}
\left\{
\begin{array}{l}
\partial_t\mu^t+ \REVISION{\nabla_x \cdot} (u^t \mu^t) = 0\\
\partial_t(\mu^tu^t)+
\REVISION{\nabla_x \cdot} (\mu^t(u^t)^{\otimes 2}) = 
\mu^t\int\gamma(\cdot-y,u^t(\cdot)-u^t(y))\mu^t(y)dy+
\eta\mu^t\nabla_{x}\psi^t
+\mu^tF_{ext}\\
\partial_s\psi^s = D\Delta_x \psi\REVISION{^s}-\kappa \psi\REVISION{^s}+
\chi*\mu^s,\ s\in[0,t],
\\ 
(\mu^0,u^0,\psi^0) = (\mu^{in},u^{in},\psi^{in})\in H^s,\ s>\tfrac d2+1.
\end{array}
\right.\nonumber
\ee
where 
$
\mu^t, u^t\in C([0,t];H^s)\cap C^1([0,T];H^{s-1}),\ 
\psi^t  \in C([0,t];H^{s})\cap C^1([0,T];H^{s-2}) \cap L^2(0,T; H^{s+1})$

Then $\rho^t(x,v):=\mu^t(x)\delta(v-u^t(x))$ solves 
the following system
\be\label{vlasobhydro}
 \left\{
\begin{array}{l}
\partial_t\rho^t+v\cdot\nabla_x\rho^t=
\REVISION{\nabla_v \cdot} (
\vla(t,x,v)\rho^t),\ 
\\ 
\vla(t,x,v)=\gamma(x,v)*\rho^t+\eta\nabla_x\psi^t(x)+F_{ext}(x),\\ 
\partial_s\psi^s(\REVISION{x})=D\Delta_x \psi\REVISION{^s}-\kappa\psi\REVISION{^s}+g(\REVISION{x},\rho^s),\ 
\psi^0=\psi^{in},\\
\rho^{0}(x,v)=\mu^{in}(x)\delta(v-u^{in}(x)).
\end{array}
\right.\nonumber
\ee
\end{Thm}

We remark that the derivation of the models is obtained assuming a monokinetic structure of the initial data at the microscopic level. This infers a very strong constraint on the initial data:  they have to be included in an $Nd$ submanifold of the $2Nd$ dimensional phase space, instead of having the freedom of filling the entire phase space $\bR^{2Nd}$.
Moreover, when dealing with applications, the monokinetic structure is clearly far from reality.
The following question arises naturally: \REVISION{neglecting} monokinetic assumption, is there any relation among the different scales? Are the \REVISION{zero-th} and first moments (in $v$) of the solution to the Vlasov system with general initial data $\rhoi$ still well approximated by the solution of the Euler system with initial data the \REVISION{zero-th} and first moments of the initial condition $\rhoi$?
One of the main objectives of the present article is to gain some insights on the answer to these questions, presenting a numerical study over different scenarios.

The paper is organized as follows: Section \ref{micres} reviews the previous numerical results established for the particle, Vlasov and Euler systems inspired by collective dynamics. In Section \ref{setmotnum} we state the setting and motivations of our computations, whose precise results are contained in Section \ref{mr}. Conclusions and perspectives conclude the paper in Section \ref{compers}. Details of the numerical schemes adopted for the different simulations are exposed in Appendix \ref{schemes}.

Let us conclude this section by recalling the three following dynamics involved in this paper, denoted by (P) for \textit{(Particles)}, (V) for \textit{(Vlasov)} and (E) for \textit{(Euler)}. 

$$
\hspace{-4cm}
\begin{array}{cl}
(P) \left\{ 
\begin{array}{l}
\dot {x_i}=v_i,\  
\dot {v_i}=F_i(t,X(t),V(t)),\ 
(X(0),V(0))\in\bR^{2dN} 
\\ \\
F_i(t,X,V)=\frac1N\sum\limits_{j=1}^N\gamma(v_i-v_j,x_i-x_j)+\eta\nabla_{x}
\varphi^t(x_i )
\\ \\
\partial_s\varphi^s(x)
=D\Delta_x\varphi\REVISION{^s}-\kappa\varphi\REVISION{^s} +f(x,X(s)),\ s\in[0,t],\ \varphi^0=\macphin, \\\ \\
f(x,X)=\sum\limits_{j=1}^N\chi (x-x_j).
\end{array} 
\right.
\nonumber
\end{array} 
$$

$$
\hspace{-5cm}
\begin{array}{cl}
(V) \left\{
\begin{array}{l}
\partial_t\rho^t+v\cdot\nabla_x\rho^t=
\REVISION{\nabla_v \cdot} (\vla(t,x,v)\rho^t),\ \rho^0=\rho^{in}\in\cP(\bR^{2d})\\ \\
\vla(t,x,v)=\gamma*\rho^t(x,v)+\eta\nabla_x\psi^t(x)
-\alpha v\\ \\
\partial_s\psi^s(x)=D\Delta_x \psi\REVISION{^s}-\kappa \psi\REVISION{^s}+g(x,\rho^s),\ \psi^0=\macphin,\\ \\ 
g(x,\rho^s)=
\chi*\rho^s(x).
\end{array}
\right.\nonumber
\end{array}
$$

$$
\begin{array}{cl}
(E)
\left\{
\begin{array}{l}
\partial_t\mu^t+ \REVISION{\nabla_x \cdot}(u^t \mu^t) = 0\\ \\
\partial_t(\mu^tu^t)+
\REVISION{\nabla_x \cdot} (\mu^t(u^t)^{\otimes 2}) = 
\mu^t\int\gamma(\cdot-y,u^t(\cdot)-u^t(y))\mu^t(y)dy+
\eta\mu^t\nabla_{x}\psi^t
-\alpha\mu^t u^t\\ \\
\partial_s\psi^s = D\Delta_x\psi\REVISION{^s}-\kappa\psi\REVISION{^s} +
\chi*\mu^s,\ s\in[0,t].
\end{array}
\right.\nonumber
\end{array}
$$

\section{Previous numerical results}\label{micres}

In this Section we present a short overview of previous seminal works of the literature concerning numerical simulations of the models in $(P)$, $(V)$ and $(E)$, focusing in particular on the ones that will be recalled in the following sections. 
As a preliminary observation, note that numerical simulations of Vlasov-type equations arising in collective dynamics modelling can be found in the literature in absence of chemotaxis. To the best of our knowledge this is the first work showing numerical simulations of the coupled system in $(V)$. The same holds for $(E)$.
A complete lists of reference dealing with numerical approaches to Vlasov and Euler equations goes beyond the scope of the paper. In this Section we refer in particular to some seminal and/or recent papers of the literature, highlighting the difference with the structure of the models considered in our paper and with the numerical scheme adopted in the different simulations. 
Moreover, the numerical tests exhibited in some of the cited works will be used to test the numerical scheme adopted in this paper over well established results.

\textbf{Particle dynamics}\label{partsys}
Concerning the particle level, numerical simulations of $(P)$ can be found in \cite{dicostanzo2}.
In fact, the discrete-continuum hybrid system investigated corresponds to $(P)$ with $d=1$, choosing a Cucker-Smale like interaction function
\begin{equation}\label{gamma_particle}
\gamma\left(v_i-v_j, x_i-x_j\right)= \frac{1}{\left(1+\frac{\left\|x_i-x_j\right\|^2}{R^2}\right)^{\beta}}(v_j-v_i)
\end{equation}
where $\beta>0$ is the alignment parameter and $R>0$ is the radius of the disk modelling each particle.

Initial data are given by initial position and velocity for each particle:
\begin{align}
	X(0)&=X_0,\quad 	 V(0)=V_0,
\end{align}
and by the initial concentration of signal, $\varphi_0=0$.

On a linearised version of the system, authors analytically prove that the particles' aggregate exponentially converges to a state in which all agents share the same position, and the velocity converges to zero. 
The behavior of the full nonlinear system is investigated with a numerical approach, showing a good agreement between numerical simulations and the theoretical results on the linearised version.
For \REVISION{the} generalized structure of $(P)$, introducing different cell mechanisms and stochasticity, numerical simulations oriented to applications on biological phenomena can be found in \cite{dicostanzo2,dicostanzo_cardio, axioms}.

\textbf{Vlasov equation}
Several numerical methods for Vlasov equations have been proposed over the years (see \cite{Filbet1, Elchina, AlbiPareschi, CarrilloVlasov, CarrilloVecil, DiMarco} and the references therein).
Focusing on Vlasov equations arising in collective dynamics phenomena, in \cite{AlbiPareschi} several stochastic Monte Carlo methods, improving classic Monte Carlo approaches, are presented, whereas in \cite{CarrilloVecil} deterministic splitting techniques, together with semi-Lagrangian and flux balance methods, are analyzed.
To the best of our knowledge, the structure of system $(V)$, where a Vlasov-like equation is coupled with a parabolic diffusion equation, has not been investigated from a numerical point of view.

\textbf{Euler system}
Two are the main features of $(E)$: \REVISION{the presence of the non-local integral term} in the Euler equations, and its coupling with the parabolic chemotactic equation.
Neglecting chemotaxis, a numerical approach to pressureless hydrodynamic alignment models has been presented in \cite{CarrilloEuler}. Here authors adopt a Lagrangian numerical scheme to approximate the solution of one dimensional hydrodynamic system. In particular, numerical simulations with different initial data are performed, in order to investigate regularity of the solutions, comparing global solutions and solutions exhibiting finite time blow up.
Adding a control term, a finite volume scheme to simulate the proposed non-local hydrodynamic system with and without control can be found in \cite{AlbiControllo}.
Concerning the coupling aspect, the seminal paper in \cite{GambaPercolation} and reference therein can be compared to the structure of $(E)$, but replacing the non-local integral term with the gradient of a phenomenological pressure.
Far from the numerical approach, the recent paper \cite{tadmor2022} presents a systematic study of the long-time swarming behavior of hydrodynamic Euler system accounting both the alignment interaction and pressure term.
A detailed description of the mixed pressure-alignment structure will be the core of Section \ref{VEpressure}.

\section{Setting and motivation of the present numerical study}\label{setmotnum}
In the present paper, we study numerically the solutions of the Vlasov and Euler systems starting with different initial data, both with and without chemotaxis.

Precisely, we study the system at particle level \eqref{sys-cuck-1D} for different values of the parameters $\beta$ and $\eta$. As we have seen in Section \ref{partsys}, parameter $\beta$ leads to conditional or unconditional flocking for the Cucker-Smale model, which can be obtained from \eqref{sys-cuck-1D} considering $\eta=0$. Values of $\eta\neq 0$ introduces the coupling of the microscopic system to the chemical activity. Here we are interested in positive chemotactic phenomena ($\eta >0$), since inspired by biological phenomena, cells move towards higher concentrations of chemicals.
The path
$$
microscopic\longrightarrow Vlasov \longrightarrow Euler
$$
corresponds to successive lacks of knowledge of the system: the meanfield limit, leading to the Vlasov equation, is obtained out of the microscopic system by averaging on all particles but one, and the passage from Vlasov to Euler retains only the moments of order zero and one in the momentum of the solution.

On the other side, the benefit of these succesive transformation is clear: the kinetic Vlasov equation deals with probabilities on finite dimensional phase spaces instead of the infinite dimensional ones of the original system, and the Euler setting describes the system by mean of densities and velocity fields on physical (configuration) spaces, that is quantities directly observable.

Finally, as explained in Section \ref{intro}, the passage form Vlasov to Euler was rigorously obtained in \cite{rt2} by a monokinetic hypothesis on the solution to Vlasov of the form $\rho^t(x,\xi)=\mu^t(x)\delta(\xi-u^t(x))$, shown to be rigorously equivalent for the couple $(\mu^t,u^t)$ satisfying the Euler system. 

\REVISION{Inspired by the results obtained for $\rho^t$ with monokinetic form,
one might wonder if, for more general (non monokinetic) solution to Vlasov, the quantities defined by \eqref{mom} still satisfy, at least approximatively, the Euler system of equations.}

From a numerical point of view, this would be even more appealing for computations based on finite difference schemes, since solving a system of two PDEs, each of them in a $m-$ dimensional space, is much more economical than solving one, having dimension $2m$.

The motivation of our work is therefore fourfold:
\begin{enumerate}[label=(\roman*)]
\item\label{i} showing that the different features of the microscopic dynamics are still present in the numerics of Vlasov and Euler situations.
\item\label{ii} comparing, for $\rho^{in}$ non-monokinetic, the moments of the solution to Vlasov as defined by \eqref{mom}  with the solution of the Euler system with initial data
\be\label{mom}
\left\{
\begin{array}{l}
\mu^{in}(x)=\int\rho^{in}(x,\xi)d\xi,\\\\
\mu^{in}(x) u^{in}(x)=\int\xi\rho^{in}(x,\xi)d\xi.
\end{array}
\right.
\ee
\item\label{iii} evaluating on the settings \ref{i} and \ref{ii} the influence of the chemical coupling by varying the value of the parameter $\eta$ from $0$ to positive values.
\item\label{iv} investigating the role of an additional pressure term in the non local Euler system and tuning its effect in order to recover agreement with the kinetic scale.
\end{enumerate}

\section{Main results}\label{mr}

In the following section we present our numerical study across the different scales, simulating different scenarios for $(P)$, $(V)$ and $(E)$. 
The aim of this paper is to analyze numerically
some of the open problems related to generalizations of theoretical results obtained in \cite{rt2}. 
We present some numerical insights concerning general initial data, far from the monokinetic ones, and we investigate the role of the different terms, especially at kinetic and hydrodynamic level.
We restrict ourselves to the 1-dimensional case ($d=1$), leaving higher dimensional cases for future researches. 
\REVISION{In particular}, for partial version of the models investigated (i.e. without chemotaxis and/or damping), we \REVISION{compare our results with the ones of the literature, when available, running the simulations starting with the same initial data}. 
Details of the numerical implementation can be found in the Appendix \ref{schemes}. 
We run the numerical codes on a laptop equipped with an Intel Core i7-1060NG7 processor and 16 GB RAM.

\vskip 1cm

\subsection{Numerical Simulations of (P) and (V)}\label{sec:PeV}

We start our analysis comparing particle and kinetic levels. To the best of our knowledge, this is the first time that a Vlasov-like equation is coupled to a chemotaxis parabolic equation. 
Numerical simulations of $(V)$, neglecting chemotaxis and damping, can be found in \cite{ AlbiPareschi, CarrilloVlasov}.  Starting with the same initial conditions of these papers, we recover the expected behavior using a different numerical scheme. Then we add the coupling with the chemotaxis and we show the effect of the gradient term. In conclusion, we simulate the scenario of a pure chemotactic interaction.

\subsubsection{alignment without chemotaxis}
\vskip 0.5cm
 \REVISION{Neglecting} chemotaxis and considering only an alignment kind of interaction, systems $(P)$ and $(V)$ reduce to a Cucker-Smale model in a particle and kinetic regime, respectively: 

\begin{align}\label{sys-cuck-nochemo}
\left\{
	\begin{array}{l}
	\dot{x}_i=v_i,\\
	\dot{v}_i=\displaystyle\frac{1}{N}\displaystyle\sum_{j=1}^N\frac{1}{\left(1+|x_i-x_j|^2\right)^{\beta}}(v_j-v_i)	\ \ \ \ i=1,..., N
	\end{array}
	\right.
\end{align}

\begin{eqnarray}\label{Vlasov_nochemo}
\left\{
\begin{array}{l}
\partial_t\rho^t+v\partial_x\rho^t=\partial_v(
\vla(t,x,v)\rho^t),\\ \\ 
\vla(t,x,v)= \displaystyle \int_{\bR^{2}} \frac{v-w}{(1+\left| x-y\right|^2)^\beta} \rho^t(y,w) dy dw. 
\end{array}
\right.
\end{eqnarray}
In the literature, numerical simulations of \eqref{Vlasov_nochemo} are often based on Montecarlo techniques. Here, dealing with a 1-dimensional case, we implemented a first-order scheme (see Appendix \ref{schemes} for details). The drawback of the constrain given by a CFL condition is overcome by the easier implementation of the scheme.
As preliminary tests, we run the the numerical code starting with the same initial distribution of \cite{AlbiPareschi}, i.e.
\begin{equation}\label{rho0}
\rho^0(x,v)= \frac{1}{2 \pi \sigma_x  \sigma_v} e^{\frac{-x^2}{2\sigma_x^2}}
\left(e^{\frac{-(v+v_0)^2}{2\sigma_v^2}} + e^{\frac{-(v-v_0)^2}{2\sigma_v^2}} \right),
\end{equation}
with $v_0=3.5$, $\sigma_x=\sqrt{0.1}$, $\sigma_v=\sqrt{0.5}$.

We consider two different scenarios, depending on the value of $\beta$, 
\REVISION{which is the key parameter of the alignment kind of interaction, as proved in \cite{haliu}. In particular, the convergence to a uniform velocity, regardless of the initial state of the system, is ensured only for value of $\beta \le 0.5$.}
\REVISION{In the following Test 1 - Test 4} the phase space representation is obtained \REVISION{considering} the space-velocity domain $\left[ -20, 20\right] \times \left[ -5, 5\right] $, with $\Delta x = \Delta v = 0.01$. The simulation runs over the time interval $\left[ 0, T\right]=[0, 5]$ with $\Delta t = 0.001$.

\begin{figure}[h!]
	\centering
		\subfigure{\includegraphics[width=6cm, height=2.3cm]{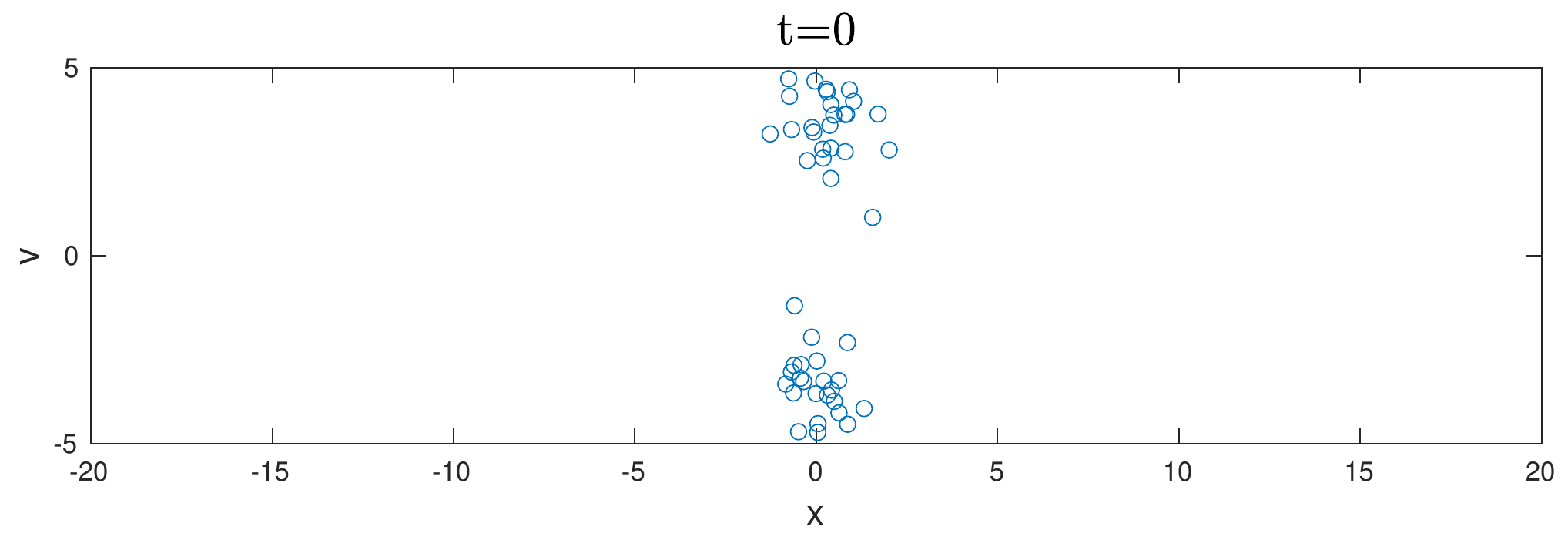}}
	\hspace{0 mm}
	\subfigure{\includegraphics[width=6cm, height=2.3cm]{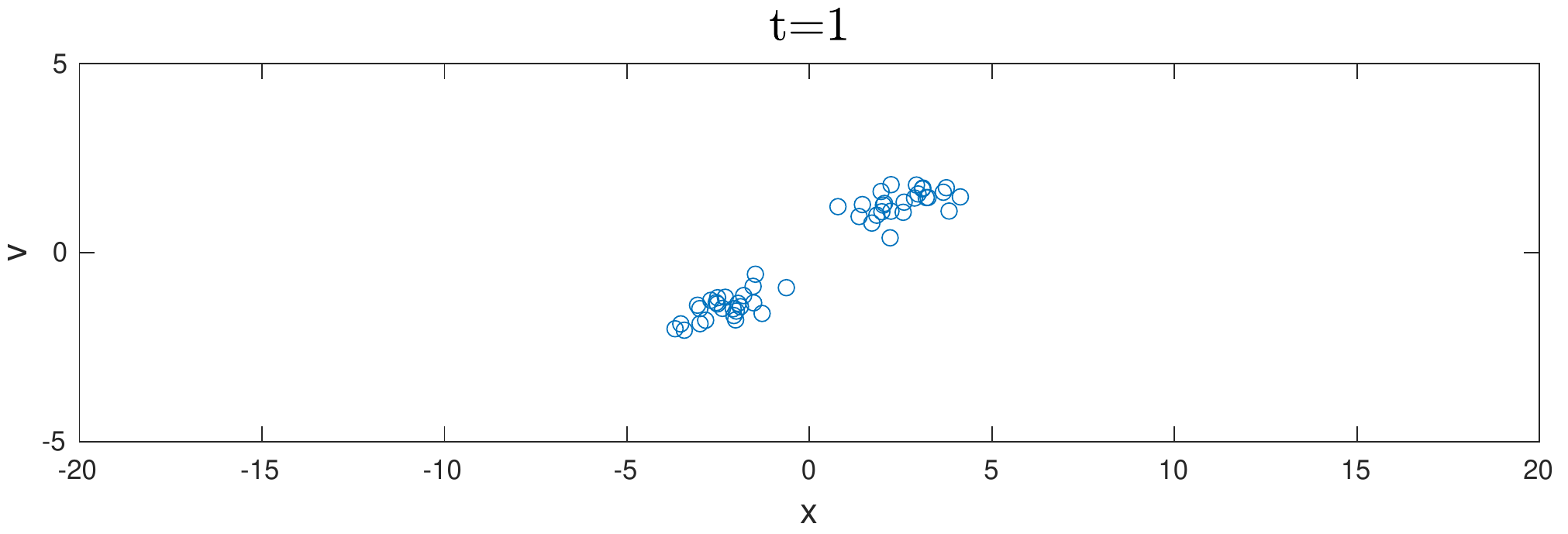}}
	\hspace{0 mm}
	\subfigure{\includegraphics[width=6cm, height=2.3cm]{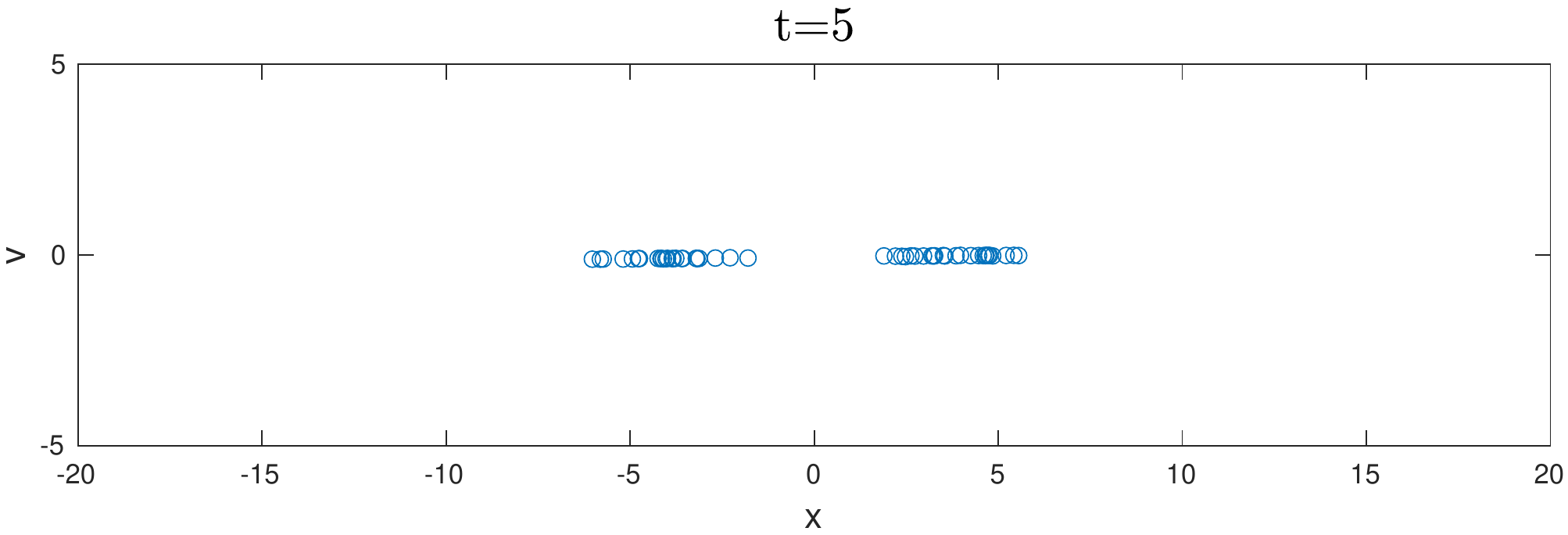}}\\
	\subfigure{\includegraphics[width=6cm, height=2.3cm]{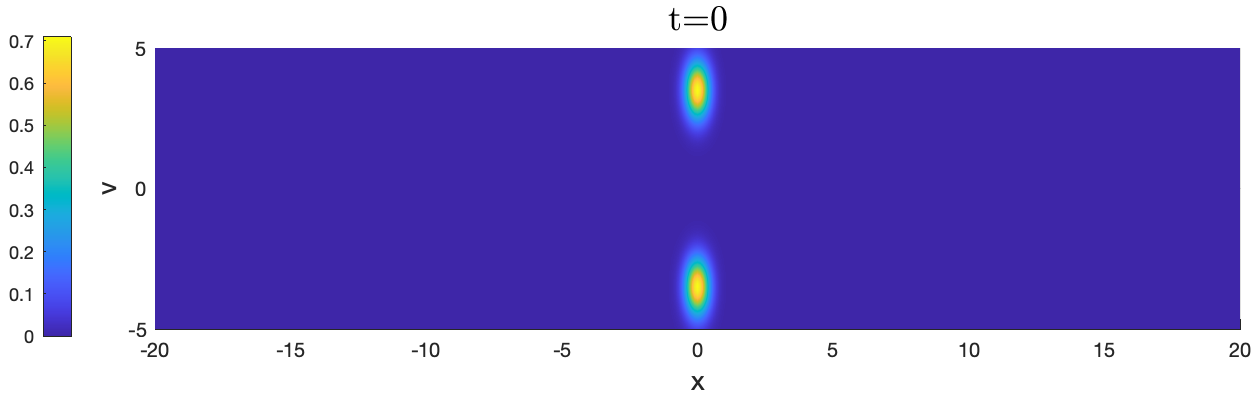}}
	\hspace{0 mm}
	\subfigure{\includegraphics[width=6cm, height=2.3cm]{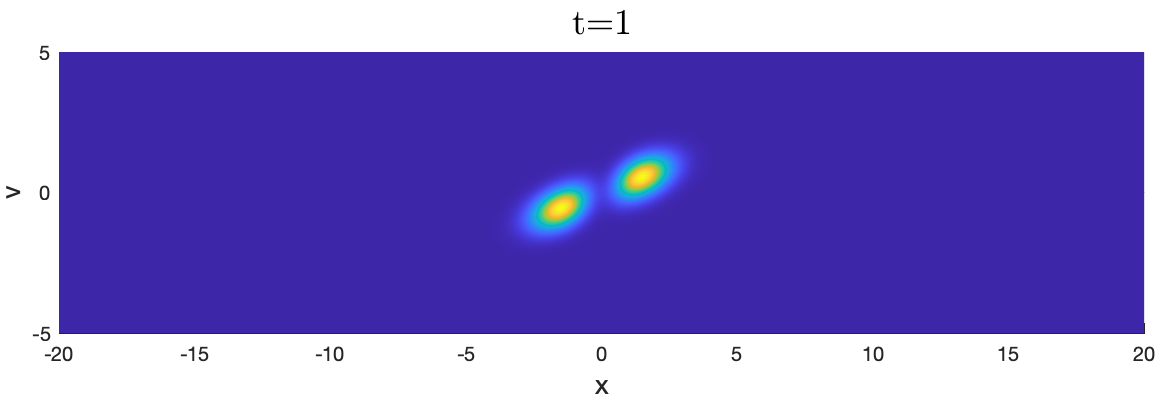}}
	\hspace{0 mm}
	\subfigure{\includegraphics[width=6cm, height=2.3cm]{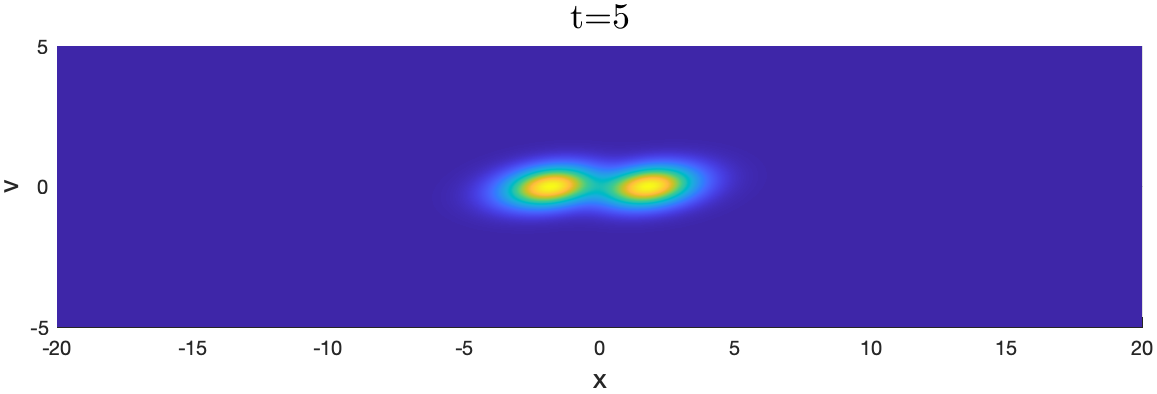}}
	\caption{Test 1: Numerical simulation of Cucker-Smale model with $\beta=0.05$, at particle level (first line) and kinetic (second line) level.}
	\label{fig:Vlasov_b005}
\end{figure}

\begin{figure}[h!]
	\centering
		\subfigure{\includegraphics[width=6cm, height=2.3cm]{figs/CS_b095_t0_revised.pdf}}
	\hspace{0 mm}
	\subfigure{\includegraphics[width=6cm, height=2.3cm]{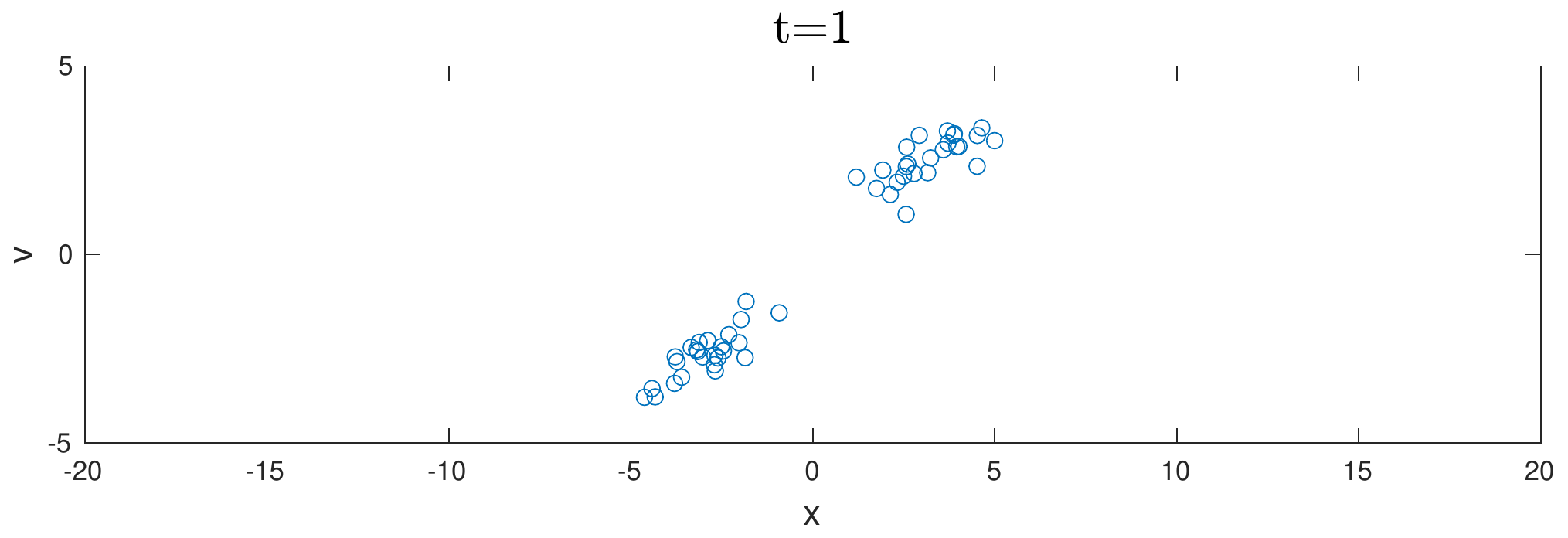}}
	\hspace{0 mm}
	\subfigure{\includegraphics[width=6cm, height=2.3cm]{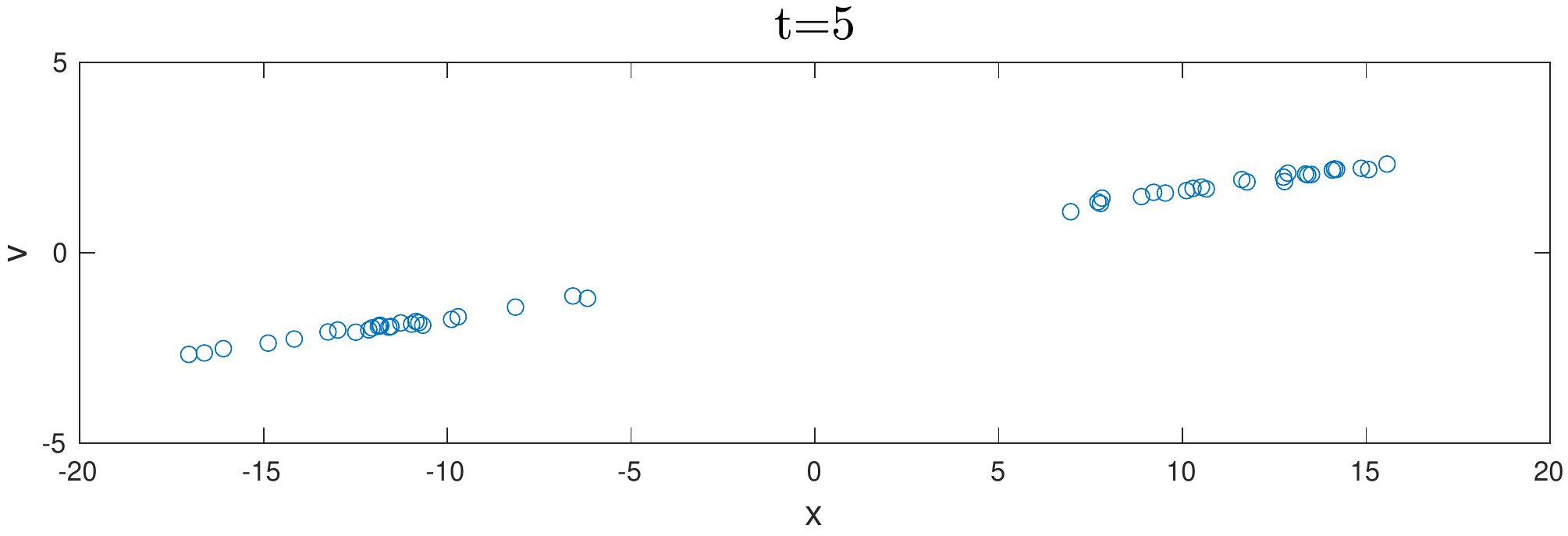}}\\
	\subfigure{\includegraphics[width=6cm, height=2.3cm]{figs/Vlas_b005e095_t0_revised.png}}
	\hspace{0 mm}
	\subfigure{\includegraphics[width=6cm, height=2.3cm]{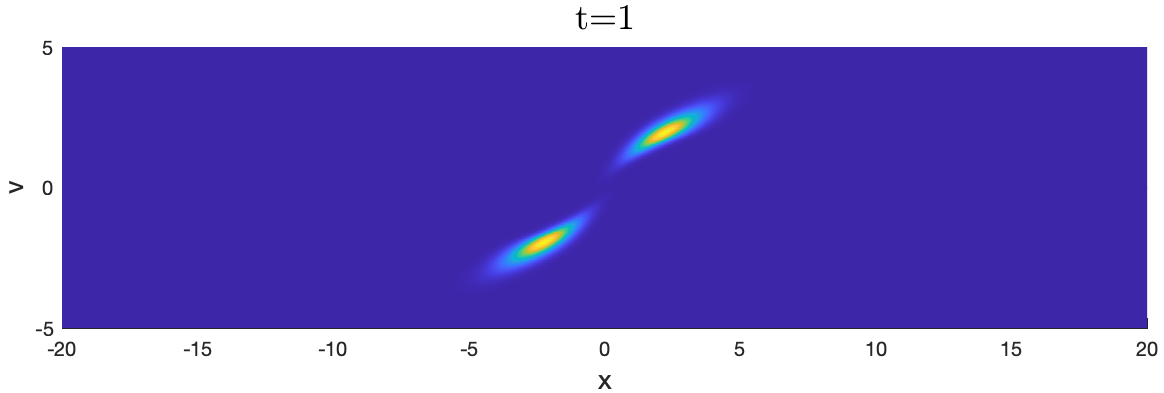}}
	\hspace{0 mm}
	\subfigure{\includegraphics[width=6cm, height=2.3cm]{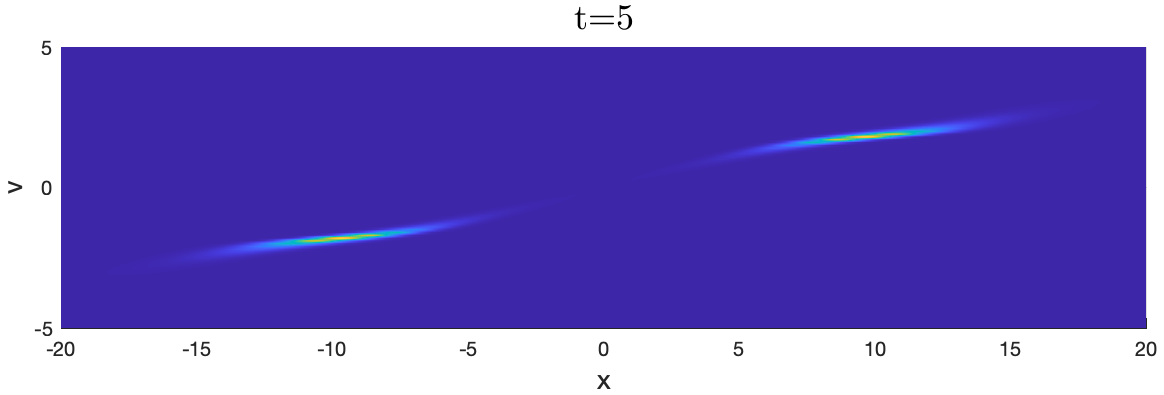}}
	\caption{Test 2: Numerical simulation of Cucker-Smale model with $\beta=0.95$, at particle level (first line) and kinetic (second line) level.}
	\label{fig:Vlasov_b095}
\end{figure}

Figure \ref{fig:Vlasov_b005} shows three different screenshots of the numerical simulation  \REVISION{ of \eqref{sys-cuck-nochemo} (first line) and \eqref{Vlasov_nochemo} (second line) performed with $\beta=0.05$.}
\REVISION{Since $\beta \le 0.5$, at particle level we observe that a flocking state is reached, and particles' velocities tend to a common (null) one.}
In a similar way, we observe that the distribution function tends to be distributed only along the spatial dimension, concentrating around the null value in the velocity space. 
 Figure \ref{fig:Vlasov_b095} shows three screenshots of the numerical simulation performed with $\beta=0.95$. 
\REVISION{In this case, since $\beta > 0.5$, the convergence to a uniform velocity is no longer reproduced.
At particle level, the two groups start migrating with positive or negative velocities, without reaching an asymptotic flocking state. The microscopic behavior is coherently reproduced at kinetic level.}

\subsubsection{alignment with chemotaxis}

On the microscopic scale, the coupling between Cucker-Smale model and chemotaxis reads:

\begin{align}\label{sys-cuck-1D}
\left\{
	\begin{array}{l}
	\dot{x}_i=v_i,\\
	\dot{v}_i=\displaystyle\frac{1}{N}\displaystyle\sum_{j=1}^N\frac{1}{\left(1+|x_i-x_j|^2\right)^{\beta}}(v_j-v_i)+\eta \partial_x \varphi \REVISION{^t} (x_i), \ \ \ \ i=1,..., N \\
	\partial_s \varphi \REVISION{^s}=D\partial_{x}^{2} \varphi \REVISION{^s} - \kappa \varphi  \REVISION{^s}+ \sum\limits_{j=1}^N\chi_{R} (x-x_j),  \ \ \ \ \REVISION{ s \in [0,t],}
	\end{array}
	\right.
\end{align}
where $R, \eta, \kappa, D >0$. In particular $\eta>0$ tunes the effect of the chemotactic gradient, and $R$ denotes the radius of the characteristic function centered in $x_j$, \REVISION{modelling the fact that the chemical signal is produced by the particles themselves}.  
Let now focus on the novel coupled Vlasov-chemotaxis system: 

\begin{eqnarray}\label{VlasovChemo}
\begin{array}{cl}
\left\{
\begin{array}{l}
\partial_t\rho^t+v\partial_x\rho^t=\partial_v(
\vla(t,x,v)\rho^t),\\ \\
\vla(t,x,v)=\displaystyle \int_{\bR^{2}} \frac{v-w}{(1+\left| x-y\right|^2)^\beta} \rho^t(y,w) dy dw+\eta\partial_x\psi^t(x),  \\ \\
\partial_s\psi^s(x)=D\partial_{x}^{2}\psi \REVISION{^s}-\kappa\psi \REVISION{^s}+ \displaystyle \int_{\mathbb{R}} \int_{x-R}^{x+R}\rho^s(y,w)dy dw, \ \ \ \ \REVISION{ s \in [0,t].} \\ 
\end{array}
\right.
\end{array}
\end{eqnarray}

To the best of our knowledge, this is the first paper presenting numerical simulations of this coupled system. The obtained results, and the comparison with the corresponding particle simulations, give some interesting insights.

\begin{figure}[h!]
	\centering
\hspace{-13 mm}
	\subfigure{\includegraphics[width=6cm, height=2.3cm]{figs/CS_b095_t0_revised.pdf}}
	\hspace{0 mm}
	\subfigure{\includegraphics[width=6cm, height=2.3cm]{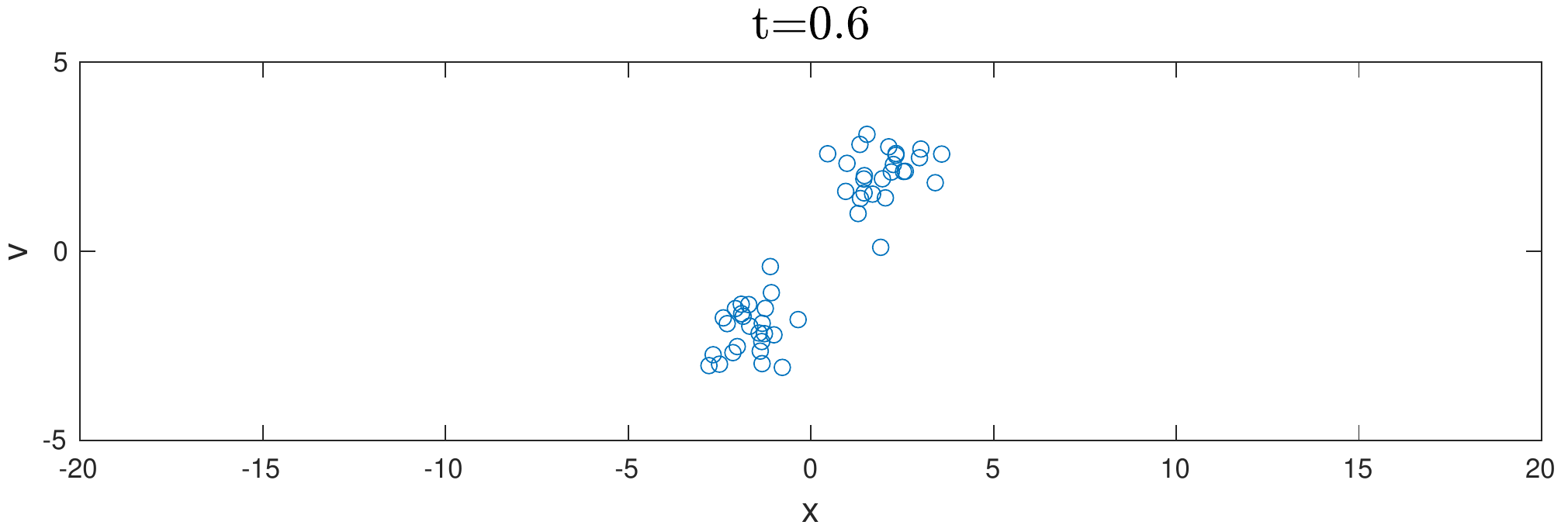}}
	\hspace{0 mm}
	\subfigure{\includegraphics[width=6cm, height=2.3cm]{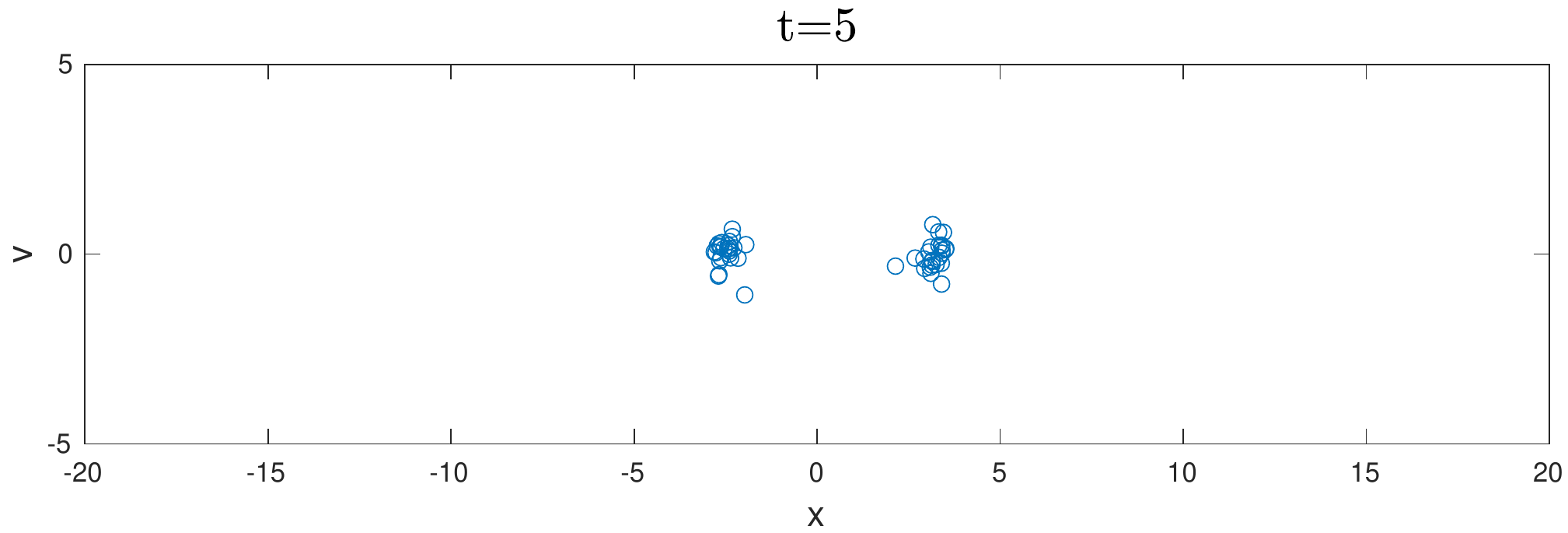}}\\
	\hspace{-13 mm}
	\subfigure{\includegraphics[width=6cm, height=2.3cm]{figs/Vlas_b005e095_t0_revised.png}}
	\hspace{0 mm}
	\subfigure{\includegraphics[width=6cm, height=2.3cm]{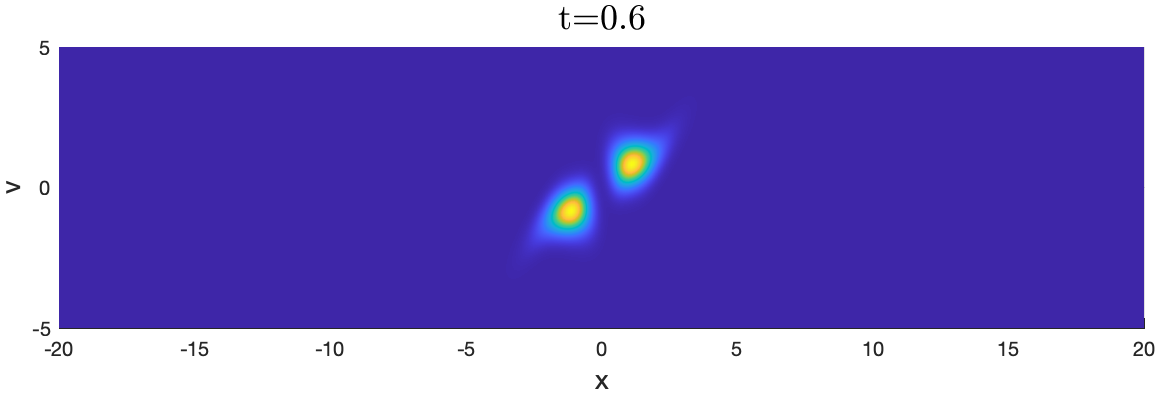}}
	\hspace{0 mm}
	\subfigure{\includegraphics[width=6cm, height=2.3cm]{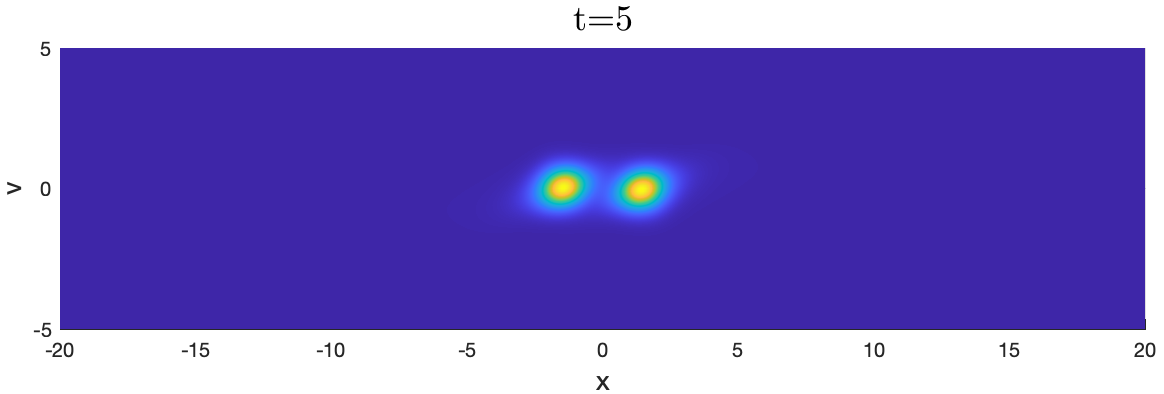}}
	\caption{Test 3: Numerical simulation of Cucker-Smale model with chemotaxis at particle level (first line) and kinetic (second line) level, with parameters as in Test 2 and $\eta=1.4$.}
	\label{fig:Vlasov_b095_chemo}
\end{figure}

In Test 3 we consider the same initial scenario and parameters of Test 2, adding chemotaxis. In particular we chose $R=10\Delta x$, $D=1$, $\kappa=0.01$, $\eta=1.4$, and initial null concentration $\psi^0=0$.
Figure \ref{fig:Vlasov_b095_chemo} shows three screenshots of the performed numerical simulation.
At particle level, in \cite{dicostanzo2} a rigorous proof of the asymptotic behavior of the solution to \eqref{sys-cuck-1D} is given. In particular, it has been proved that the introduction of the chemotactic effect ensures convergence to a flocking state, even in the cases of conditional flocking ($\beta>0.5$). 
Numerical evidence show that the behavior analytically proved on the microscopic scale, is recovered also in the kinetic regime.

\subsubsection{chemotaxis without alignment}
To conclude our analysis, we present a numerical simulation in the case of a pure chemotactic kind of interaction, hence neglecting the alignment term. The system reads:

\begin{eqnarray}\label{VlasovOnlyChemo}
\begin{array}{cl}
\left\{
\begin{array}{l}
\partial_t\rho^t+v\partial_x\rho^t=\partial_v(
\vla(t,x,v)\rho^t),\\ \\
\vla(t,x,v)=\displaystyle \eta\partial_x\psi^t(x),\\ \\
\partial_s\psi^s(x)=D\partial_{x}^{2}\psi\REVISION{^s}-\kappa\psi\REVISION{^s}+ \displaystyle  \int_{\mathbb{R}} \int_{x-R}^{x+R}\rho^s(y,w)dy dw, \ \ \ \ \REVISION{ s \in [0,t].}\\ 
\end{array}
\right.
\end{array}
\end{eqnarray}

\begin{figure}[h!]
	\centering
	\hspace{-13 mm}
	\subfigure{\includegraphics[width=6cm, height=2.3cm]{figs/CS_b095_t0_revised.pdf}}
	\hspace{0 mm}
	\subfigure{\includegraphics[width=6cm, height=2.3cm]{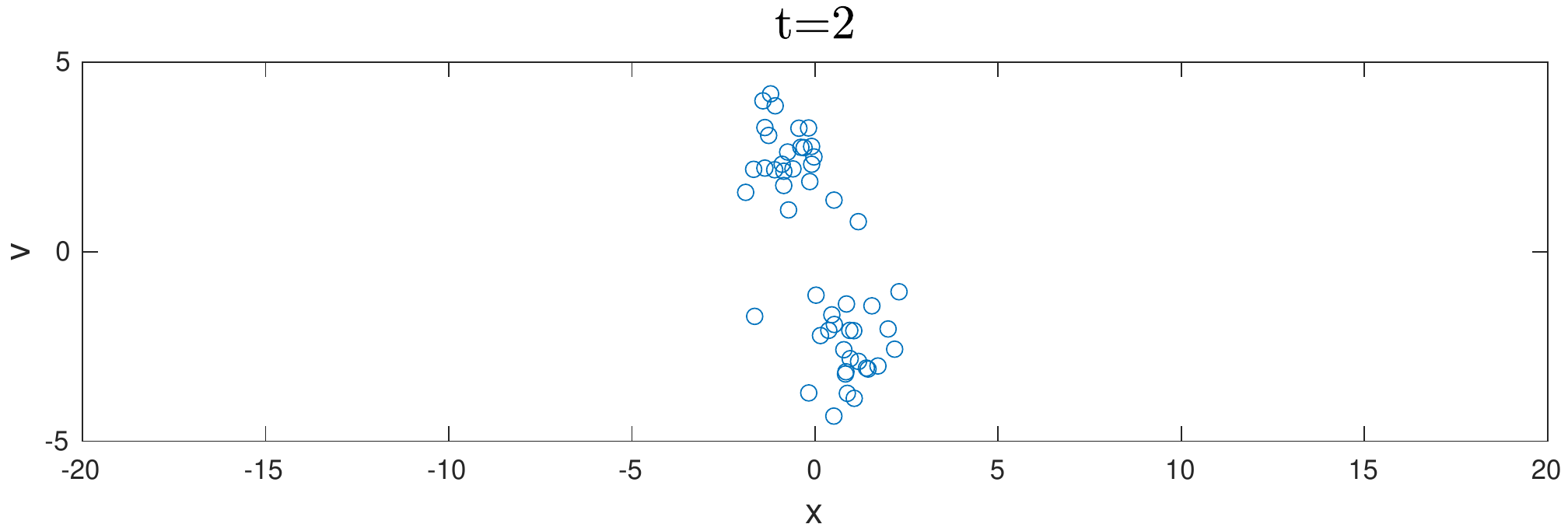}}
	\hspace{0 mm}
	\subfigure{\includegraphics[width=6cm, height=2.3cm]{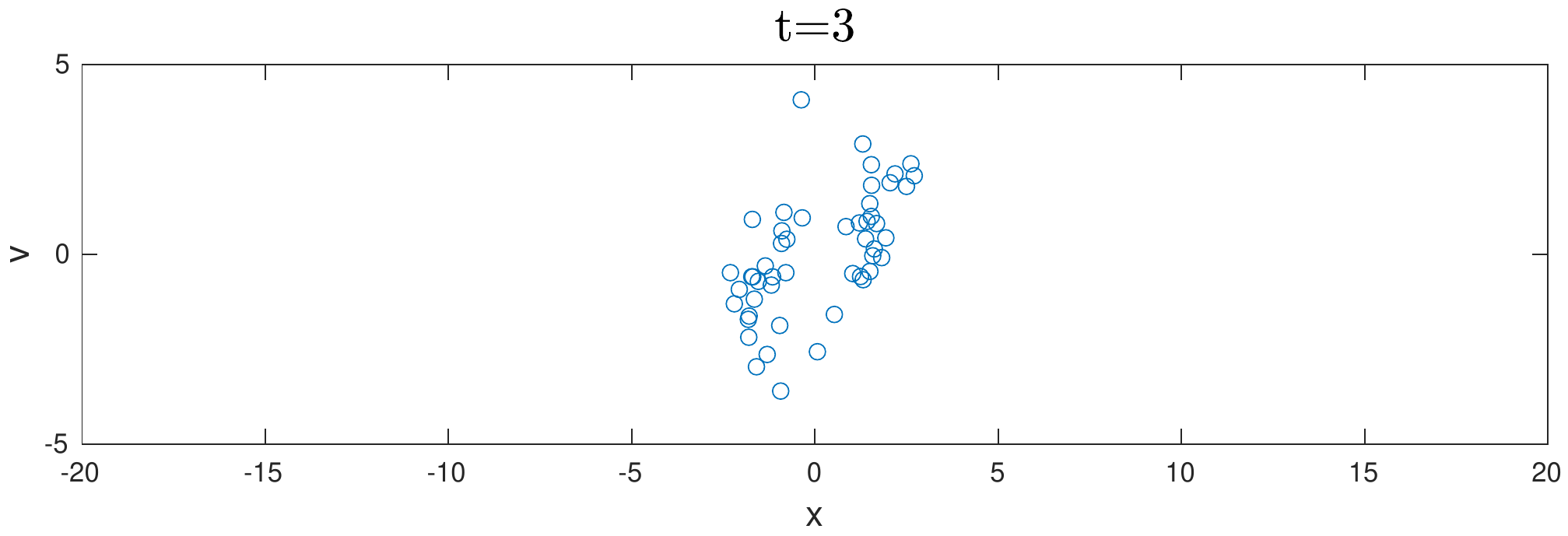}}\\
	\hspace{-13 mm}
	\subfigure{\includegraphics[width=6cm, height=2.3cm]{figs/Vlas_b005e095_t0_revised.png}}
	\hspace{0 mm}
	\subfigure{\includegraphics[width=6cm, height=2.3cm]{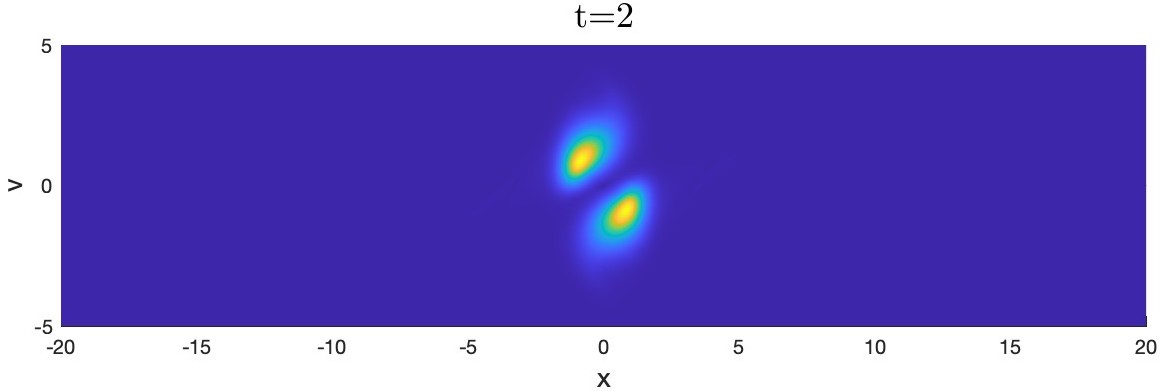}}
	\hspace{0 mm}
	\subfigure{\includegraphics[width=6cm, height=2.3cm]{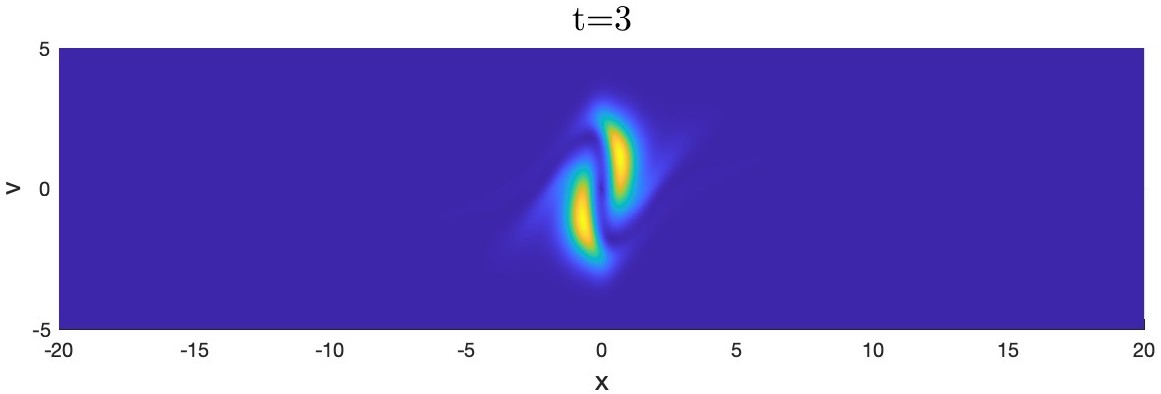}}
	\caption{Test 4:
	Numerical simulation of \eqref{VlasovOnlyChemo} at particle level (first line) and kinetic (second line) level, with parameters as in Test 3.}
	\label{fig:Vlasov_onlychemo}
\end{figure}

Figure \ref{fig:Vlasov_onlychemo} shows three screenshots of the performed numerical simulation. Both at particle and kinetic level,  the absence of the alignment kind of interaction results in an oscillating behavior around the center of mass. 
As in Test 3, due to the chemotactic gradient, particles tend to concentrate, but a consensus state is not reached.

We remark that the initial distribution of Tests 1-4 is clearly far from the monokinetic case. The comparison shown in this section shows a good level of correspondence between particle and kinetic scale in the different scenarios, and it represents the preliminary step of our analysis.
In the next Sections we focus on numerical simulations of the hydrodynamic system $(E)$, recovering previous results of the literature and exploring the novel coupling with the chemotaxis equation.

\subsection{Numerical Simulations of (E)}\label{simu_carrillo}

As seen in Section \ref{micres}, pressureless Euler systems have been investigated in the literature, both from a theoretical and numerical perspective. Also at this level, the novel aspect lies in the coupling with the chemotaxis equation. 
The model here considered writes as a hyperbolic-parabolic non local system for macroscopic quantities: density, momentum and concentration of a chemoattractant:

\begin{eqnarray}\label{Euler1D}
\begin{array}{cl}
\left\{
\begin{array}{l}
\partial_t\mu^t+\partial_x(u^t \mu^t) = 0\\ \\
\partial_t(\mu^tu^t)+
\partial_x(\mu^t(u^t)^{2}) = 
\mu^t\int\gamma(\cdot-y,u^t(\cdot)-u^t(y))\mu^t(y)dy\\ \\
\ \ \ \ \ \ \ \ \ \ \ \ \ \ \ \ \ \ \ \ \ \ \ \ \ \ \ \ \ \ \ \ +\eta\mu^t\partial_{x}\psi^t\\ \\
\partial_s\psi^s = D\partial^2_{x}\psi\REVISION{^s}-\kappa\psi\REVISION{^s} +
\displaystyle \int_{x-R}^{x+R}\mu^s(y)dy,\ s\in[0,t].
\end{array}
\right.
\end{array}
\end{eqnarray}

We investigate the role of chemotaxis, comparing the case with $\eta=0$ and $\eta \neq0$.
The numerical scheme implemented is based on relaxation techniques originally proposed in \cite{AregbaNatalini} (see Appendix \ref{schemes} for the details).
As a preliminary test, we simulate the system in absence of chemotaxis, starting with the initial condition chosen as in \cite{CarrilloEuler}. Then we explore the role of chemotaxis and damping in the same scenario.

\subsubsection{without chemotaxis: finite time blow up of solution}

Neglecting chemotaxis, system (\ref{Euler1D}) reduces to a 1-dimensional pressureless Euler-alignment system:

\COMMENT{Note: Carrillo et al. do not consider damping}
\begin{eqnarray}\label{Euler1D_nochemo}
\begin{array}{cl}
\left\{
\begin{array}{l}
\partial_t\mu^t+\partial_x(u^t \mu^t) = 0\\ \\
\partial_t(\mu^tu^t)+
\partial_x(\mu^t(u^t)^{2}) = 
\mu^t\int\gamma(\cdot-y,u^t(y)-u^t(\cdot))\mu^t(y)dy.\\ 
\end{array}
\right.
\end{array}
\end{eqnarray}

Analytical results and numerical insights on the stability of flock solutions for this system can be found in \cite{CarrilloEuler2, CarrilloEuler}.
In \cite{CarrilloEuler2}, authors distinguish between \textit{supercritical } and \textit{subcritical region}: initial data lying in the \textit{supercritical region} lead to blow up of the solution, whereas starting with data in the \textit{subrcritical region} global existence of the solution is ensured.
In \cite{CarrilloEuler} different numerical simulations \REVISION{in the spatial domain $[-0.75,0.75 ]$} with data lying in both regions are performed.
\REVISION{In Test 5 we approximate the solution of \eqref{Euler1D_nochemo} starting with the same initial data defined in \cite{CarrilloEuler}: }

\begin{equation}\label{dens_0}
\mu^0(x)= c_1\cos\left (\frac{\pi x}{1.5} \right),
\end{equation}

\begin{equation}\label{vel_0}
u^0(x)= -c_2 \sin \left(\frac{\pi x}{1.5} \right)
\end{equation}
where $c_1$ is the normalization factor, and $c_2>0$ is varied in the different simulations.

\begin{figure}[h!]
	\centering
\subfigure{\includegraphics[width=5cm, height=5cm]{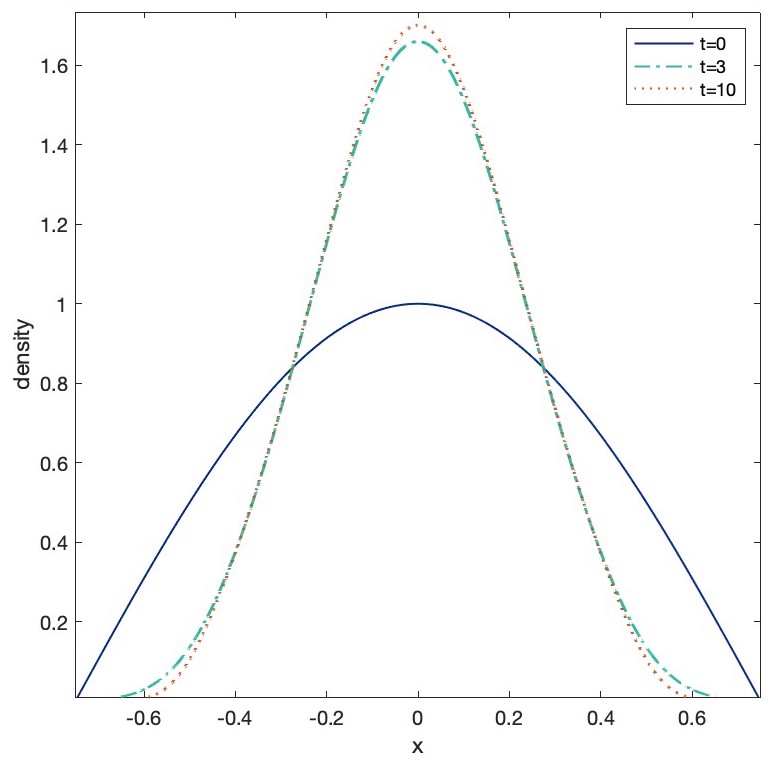}}
\subfigure{\includegraphics[width=5cm, height=5cm]{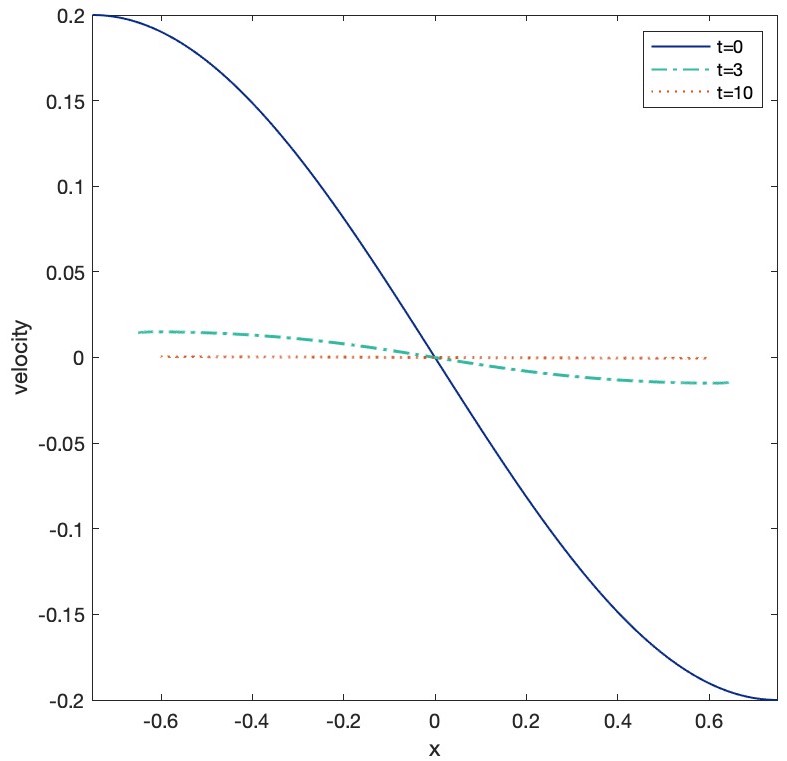}}\\	
\subfigure{\includegraphics[width=5cm, height=5cm]{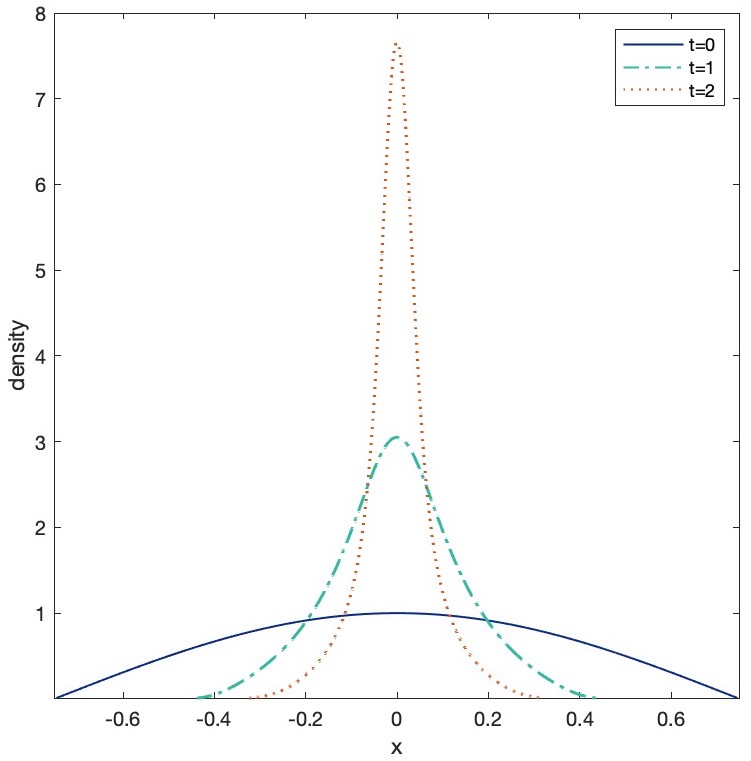}}
\subfigure{\includegraphics[width=5cm, height=5cm]{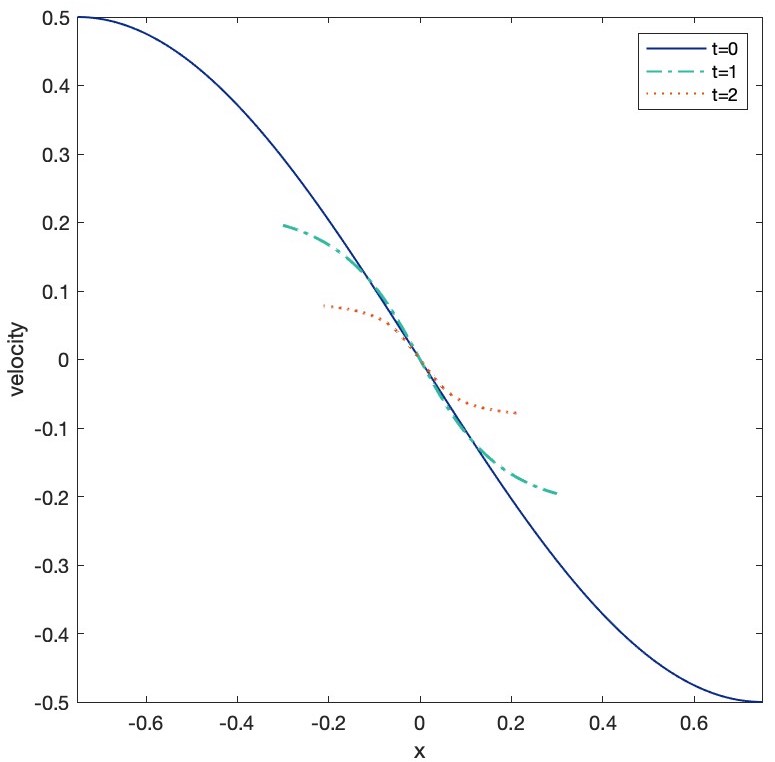}}\\	
	\caption{Test 5: Numerical simulation of \eqref{Euler1D_nochemo}-\eqref{vel_0} for $c_2=0.2$ (first line) and $c_2=0.5$ (second line). }
	\label{fig:Euler_nochemo}
\end{figure}

The results obtained are in good agreement with the ones presented in \cite{CarrilloEuler}, where a Lagrangian numerical scheme is adopted. 
Figure \ref{fig:Euler_nochemo} (first line) shows the density and velocity profile obtained for $c=0.2$ for three different time instants. The initial data lies in the subcritical region, and the convergence of the velocity to zero corresponds to the expected global consensus.
On the contrary, Figure \ref{fig:Euler_nochemo} (second line) shows the density and velocity profile obtained for $c=0.5$, hence initial data in the supercritical region. 
\REVISION{Plotting the velocity at different times, we observe that the derivative is getting larger negative at the origin as time goes. This implies high and rapidly increasing values of the density solution, which almost cannot be solved numerically after $t \approx 3$.
The obtained results represent our starting point for numerical investigations on the effect of the introduction of a chemotactic effect and a damping term.}

\COMMENT{To check: other papers which simulate Nonlocal Euler (with no chemo). Only this one of Carrillo?}

\subsubsection{with chemotaxis and/or damping}

Let now consider the coupled system in \eqref{Euler1D} adding a damping term:
\begin{eqnarray}\label{Euler1D_damping}
\begin{array}{cl}
\left\{
\begin{array}{l}
\partial_t\mu^t+\partial_x(u^t \mu^t) = 0\\ \\
\partial_t(\mu^tu^t)+
\partial_x(\mu^t(u^t)^{2}) = 
\mu^t\int\gamma(\cdot-y,u^t(\cdot)-u^t(y))\mu^t(y)dy\\ \\
\ \ \ \ \ \ \ \ \ \ \ \ \ \ \ \ \ \ \ \ \ \ \ \ \ \ \ \ \ \ \ \ +\eta\mu^t\partial_{x}\psi^t
-\alpha \mu^t u^t\\ \\
\partial_s\psi^s = D\partial^2_{x}\psi\REVISION{^s}-\kappa\psi\REVISION{^s} +
\displaystyle \int_{x-R}^{x+R}\mu^s(y)dy,\ s\in[0,t].
\end{array}
\right.
\end{array}
\end{eqnarray}
where $\eta$ and $\alpha$ are nonnegative constants regulating the chemotactic and the damping influence, respectively.
\REVISION{From a modelling perspective,} the introduction of damping is justified and necessary for biological applications, \REVISION{in particular when dealing with cells moving on a substrate}.

\begin{figure}[h!]
	\centering
\subfigure{\includegraphics[width=5cm, height=5cm]{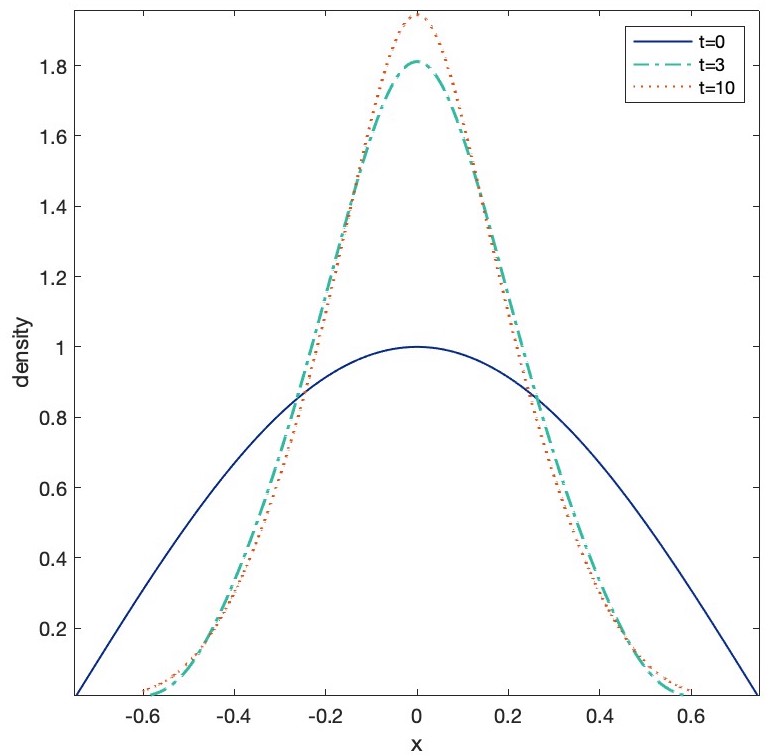}}
\subfigure{\includegraphics[width=5cm, height=5cm]{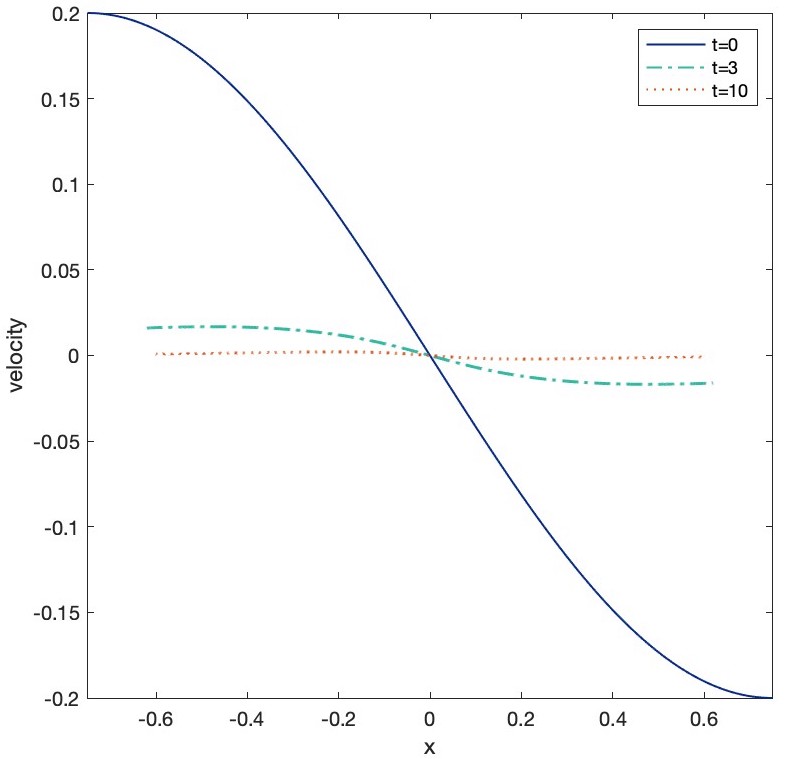}}\\	
\subfigure{\includegraphics[width=5cm, height=5cm]{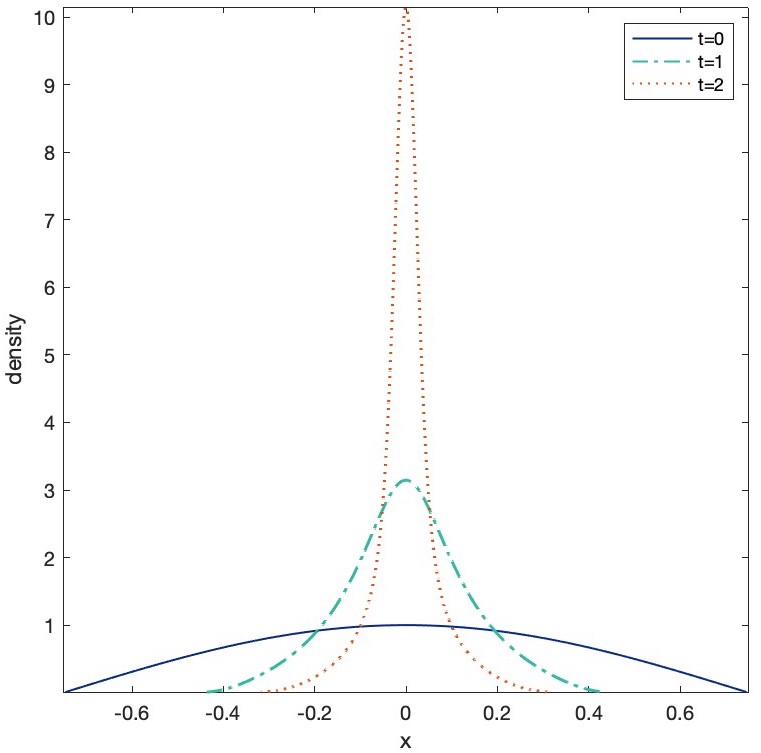}}
\subfigure{\includegraphics[width=5cm, height=5cm]{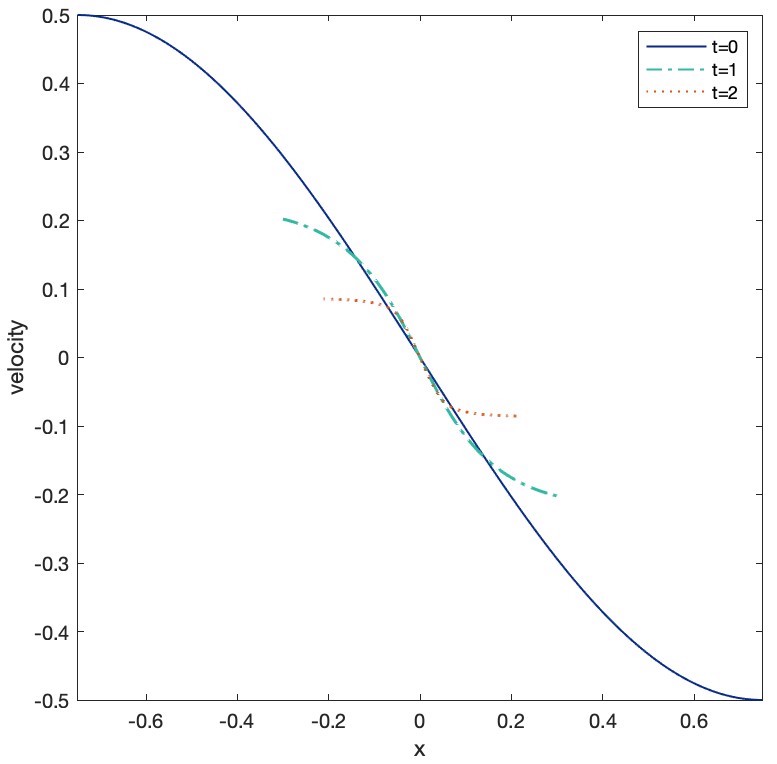}}\\	
\subfigure{\includegraphics[width=5cm, height=5cm]{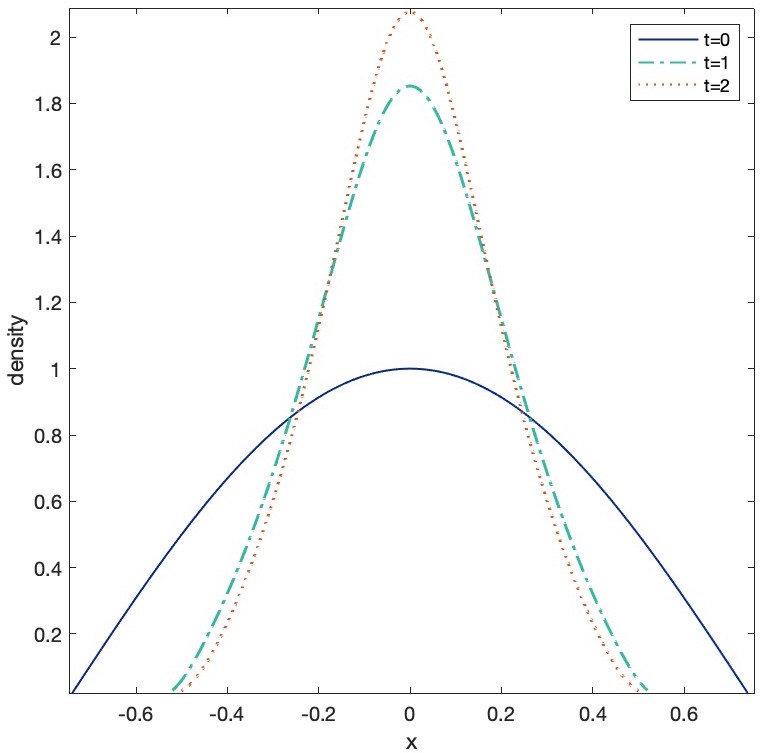}}
\subfigure{\includegraphics[width=5cm, height=5cm]{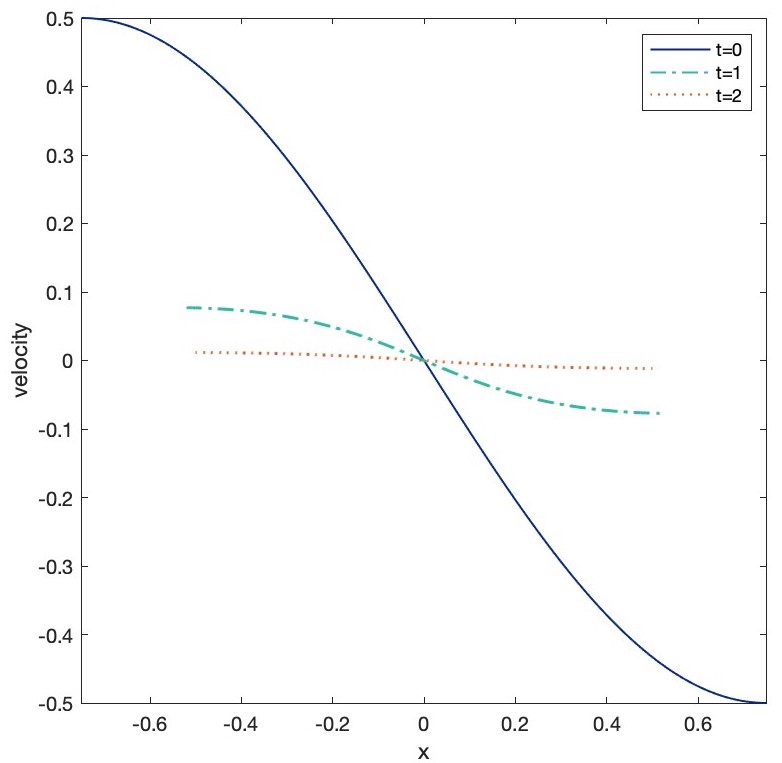}}
	\caption{Test 6: Numerical simulation of \eqref{Euler1D_damping} with initial condition \eqref{dens_0}-\eqref{vel_0}. First line: $c_2=0.2$, $\eta=1$, $\alpha=0$; second line: $c_2=0.5$, $\eta=1$, $\alpha=0$; third line $c_2=0.5$, $\eta=0$, $\alpha=1$. }
	\label{fig:Euler_chemo}
\end{figure}

Figure \ref{fig:Euler_chemo} shows the results of three different scenarios.
First, we simulate the setting with data in the subcritical region already seen ($c_2=0.2$, $\alpha=0$) adding the chemical influence ($\eta=1$). Numerical evidence shows that the global consensus is still preserved, and the presence of a chemotactic influence \REVISION{results in higher values of the density}.
In the other scenarios (Figure \ref{fig:Euler_chemo}, second and third lines) we investigate the role of chemotaxis and damping in suppressing the blow up phenomena ($c_2=0.5$). Numerical evidence (Figure \ref{fig:Euler_chemo}, second line) shows that a pure chemotactic effect ($\eta=1$, $\alpha=0$) enhances blow up phenomena: density value at the origin is higher than the case without chemotaxis, and the derivative of the velocity is sharper. 
The presence of a damping term seems to play the crucial role, as shown in (Figure \ref{fig:Euler_chemo}, third line). Here, the blow up exhibited in the same setting without damping (see Figure \ref{fig:Euler_nochemo}) is suppressed: velocity converges to zero and global consensus is reached.

\subsection{Numerical comparison between (V) and (E) }

In the following sections we present the main results of our work. 
One of the main aims is to get numerical insights on the behavior of the solutions of $(V)$ and $(E)$, in particular starting with non monokinetic initial data, being not covered by theoretical results.

First, we introduce the following notation for the zero-th and first order moments of the solution of $(V)$, respectively
\begin{equation}\label{nu0}
\nu^t_0(x)=\int \rho^t(x,\xi)d\xi,
\end{equation}
\begin{equation}\label{nu1}
\nu^t_1(x)=\int \xi \rho^t(x,\xi)d\xi,
\end{equation}
for any $t\ge0$.

In the following tests we compare the moments of the solution to Vlasov with the solution of the Euler system with initial data 

\be\label{mom_data}
\left\{
\begin{array}{l}
\mu^{0}(x)=\nu^0_0(x),\\ \\
Q^{0}(x)= \nu^0_1(x),
\end{array}
\right.
\ee
where $Q^t:=\mu^t u^t$ denotes the momentum, for any $t \ge 0$.
 Choosing different $\rho^{0}$ we simulate different scenarios, ranging from almost monokinetic to largely non-monokinetic initial data.
The aim is to compare the behavior of $\nu^t_0$ and $\nu^t_1$ with density $\mu^t$ and momentum $Q^t$ solving $(E)$ starting with \eqref{mom_data}.
 
\vspace{0.5cm}

\textit{Monokinetic initial data:}
We recall that analytical results obtained in \cite{rt2} ensure a complete correspondence between Vlasov moments and Euler solutions, in case of monokinetic initial data. 
As a preliminary test, we consider a numerical approximation of the monokinetic initial data, choosing a gaussian distribution with low value of $\sigma_v$ parameter. In greater detail, in Test 7 we assume

\begin{equation}\label{rho0_monokin}
\rho^0(x,v)= \frac{1}{2 \pi \sigma_x  \sigma_v} e^{\frac{-(x-x_0)^2}{2\sigma_x^2}+\frac{-(v-v_0)^2}{2\sigma_v^2}} 
\end{equation}
with $x_0=-2$, $v_0=1.5$, $\sigma_x=\sqrt{0.2}$, $\sigma_v=\sqrt{0.001}$. 
Clearly, even choosing a small value for $\sigma_v$, we are considering an approximation of the monokinetic case, and we do not expect a full correspondence with the analytical results on this case.

\begin{figure}[h!]
	\centering
\subfigure[]{\includegraphics[width=5cm, height=4.2cm]{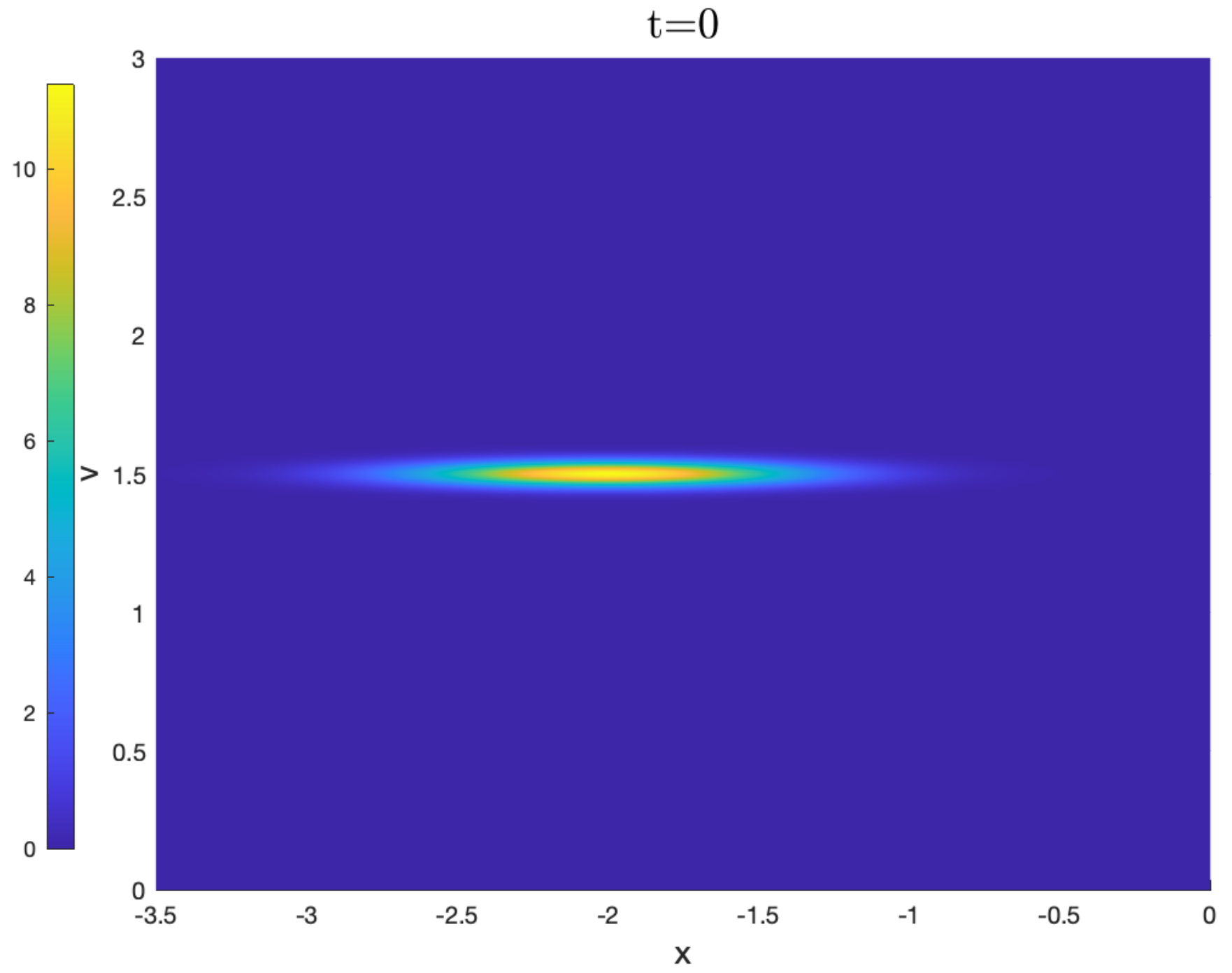}}
\subfigure[]{\includegraphics[width=5cm, height=4cm]{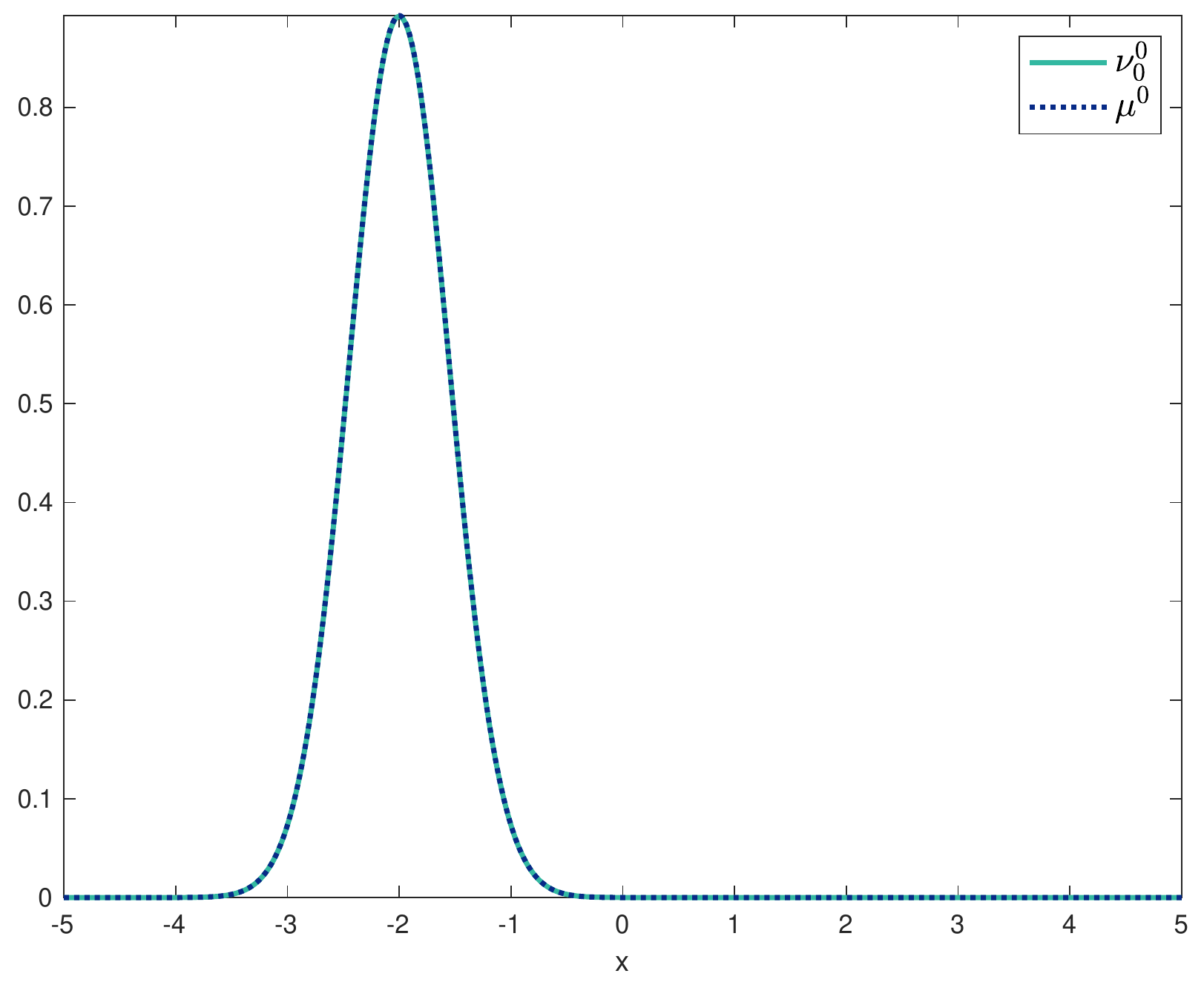}}
\subfigure[]{\includegraphics[width=5cm, height=4cm]{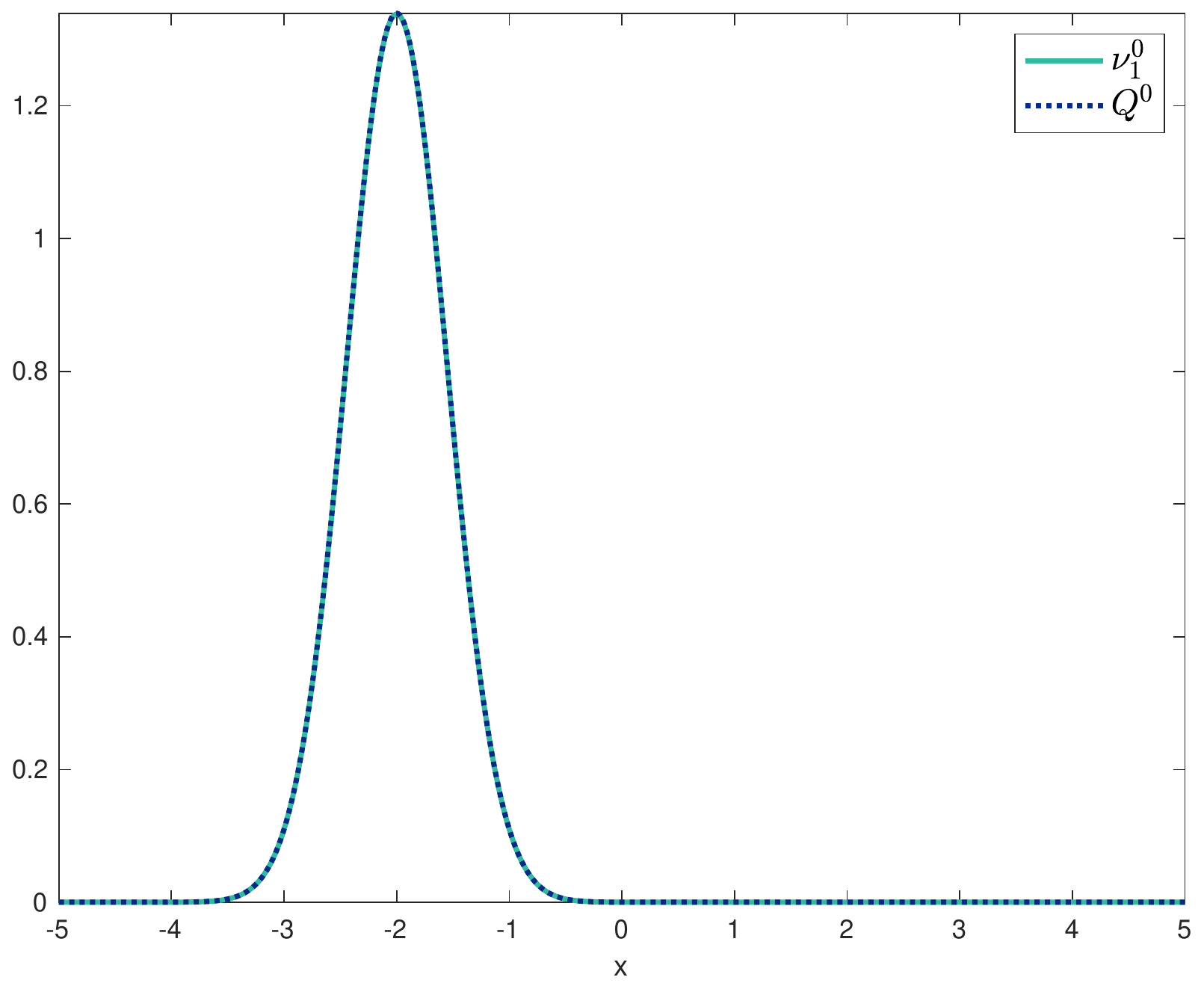}}
	\caption{Test 7, Test 11: Initial condition a) $\rho^0$ in the phase space, b) initial density $\mu^0$ and c) initial momentum $Q^0$.  }
	\label{fig:monokinetic_initialcondition }
\end{figure}

\begin{figure}[h!]
	\centering
\raisebox{10pt}{\parbox[b]{1.2\textwidth}{ \hspace{3cm} \small{Vlasov phase space} \hspace{2.5cm} $\nu_0^t$ and $\mu^t$ \hspace{3.5cm} $\nu_1^t$ and $Q^t$}}\\
\subfigure{\includegraphics[width=5cm, height=4.2cm]{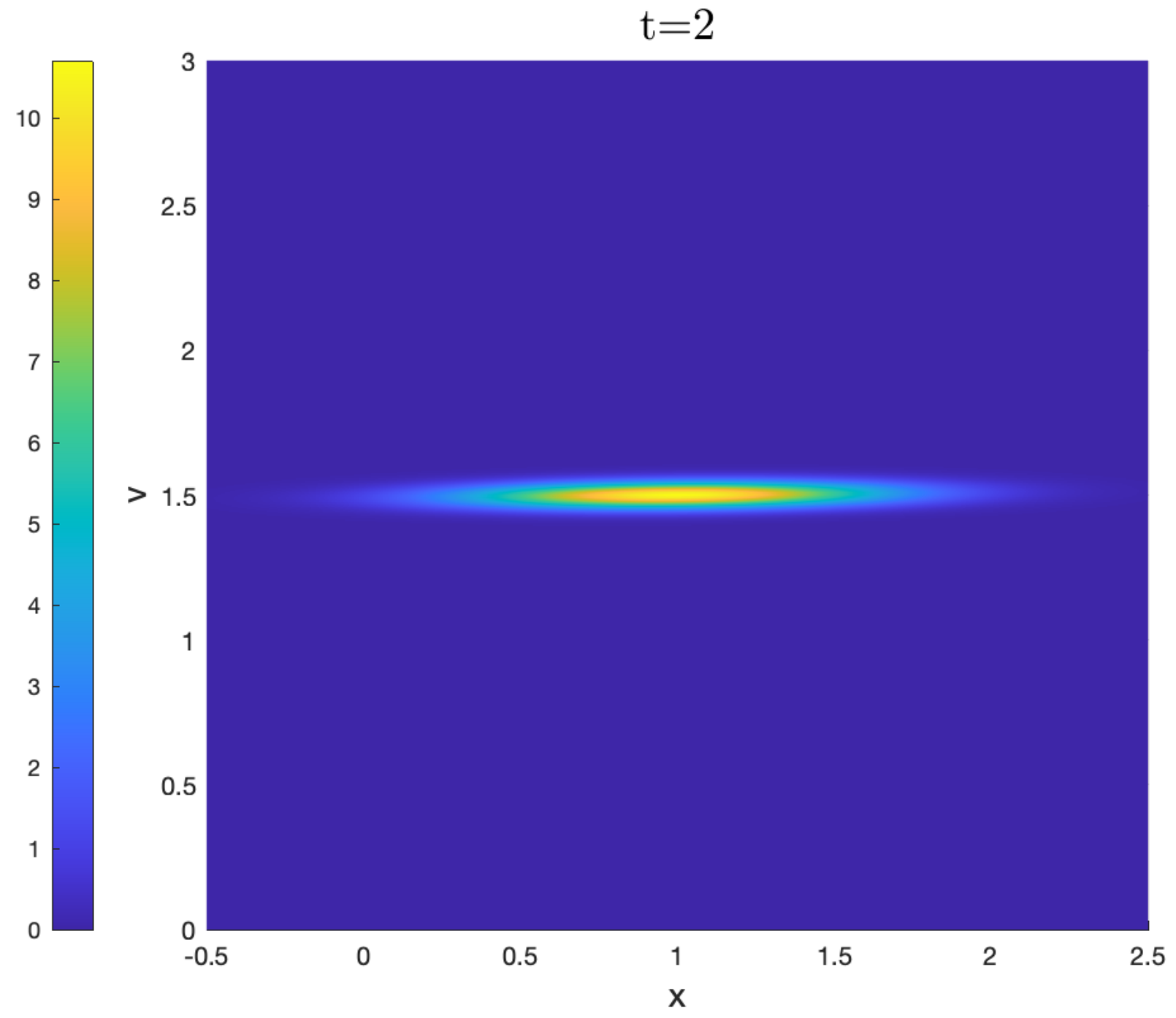}}
\subfigure{\includegraphics[width=5cm, height=4cm]{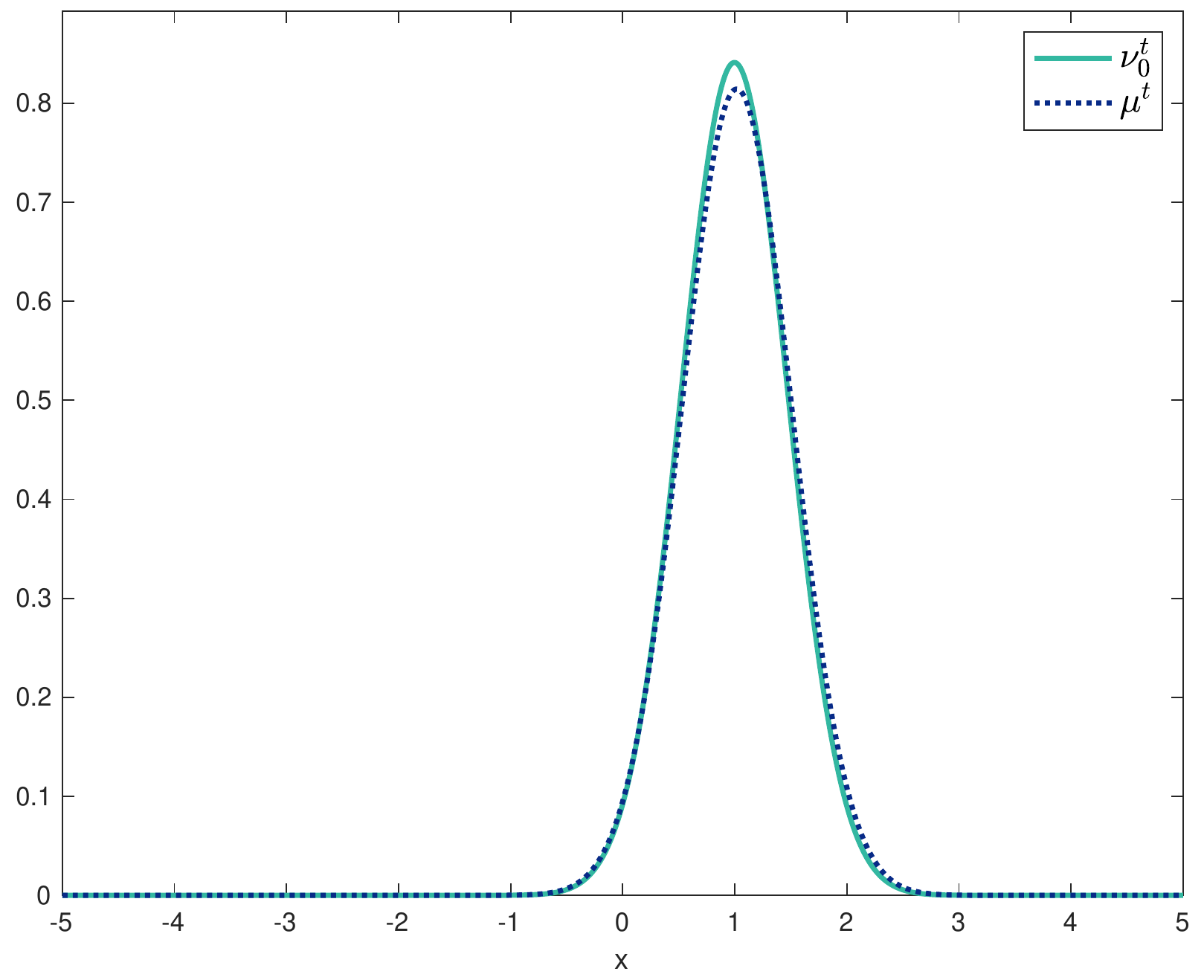}}
\subfigure{\includegraphics[width=5cm, height=4cm]{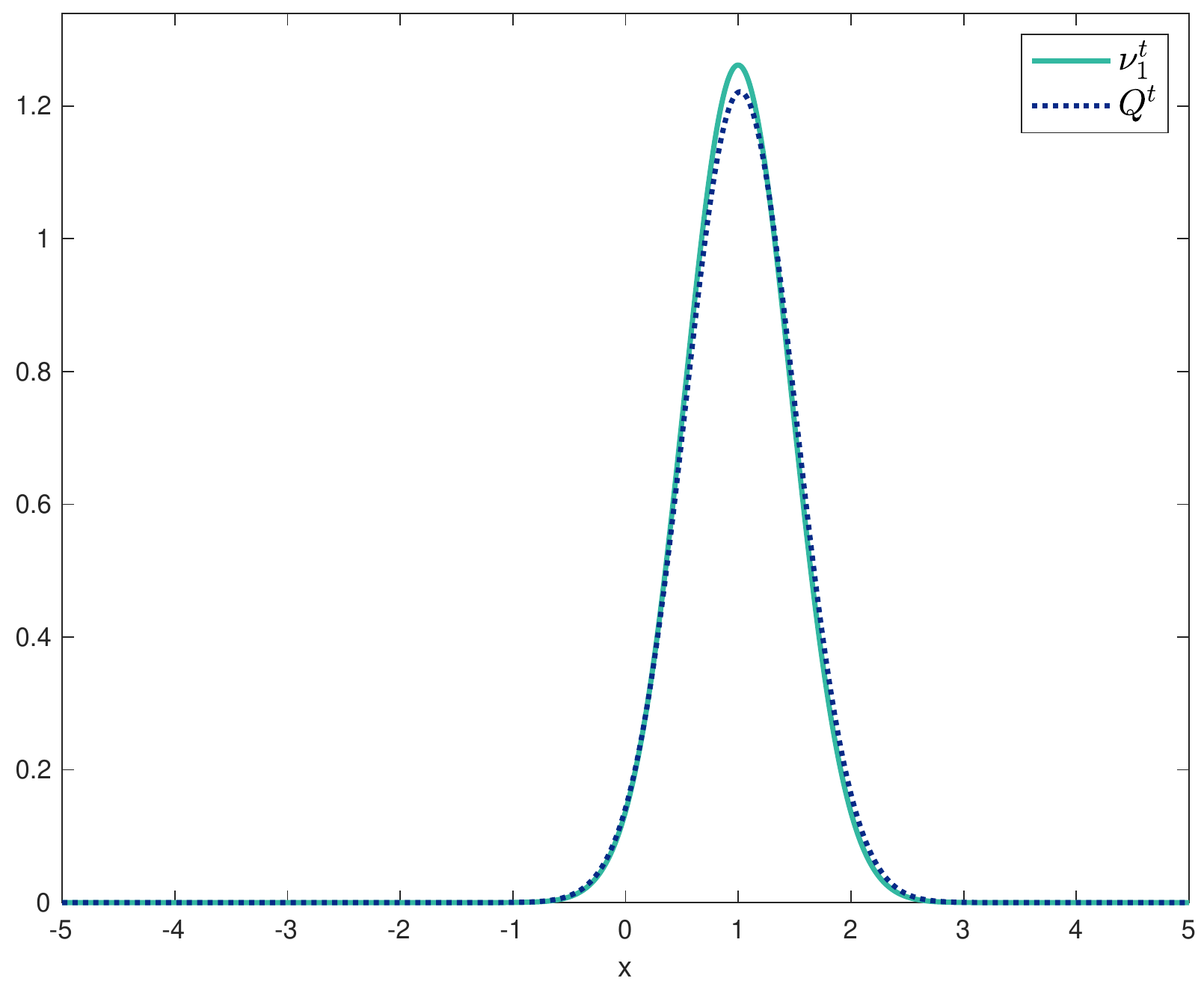}}\\
\hspace{0.5cm}
\subfigure{\includegraphics[width=4.5cm, height=4.2cm]{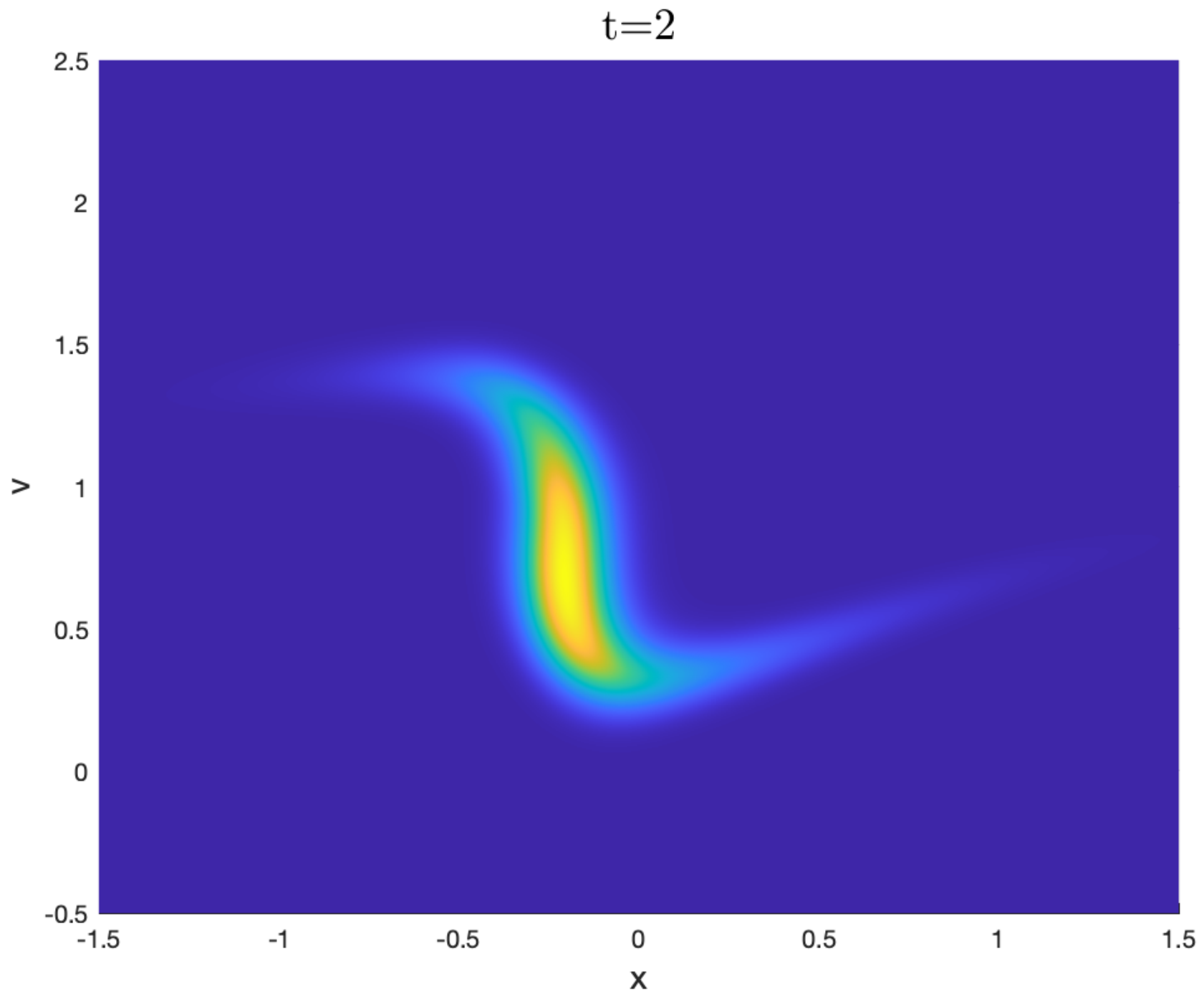}}
\subfigure{\includegraphics[width=5cm, height=4cm]{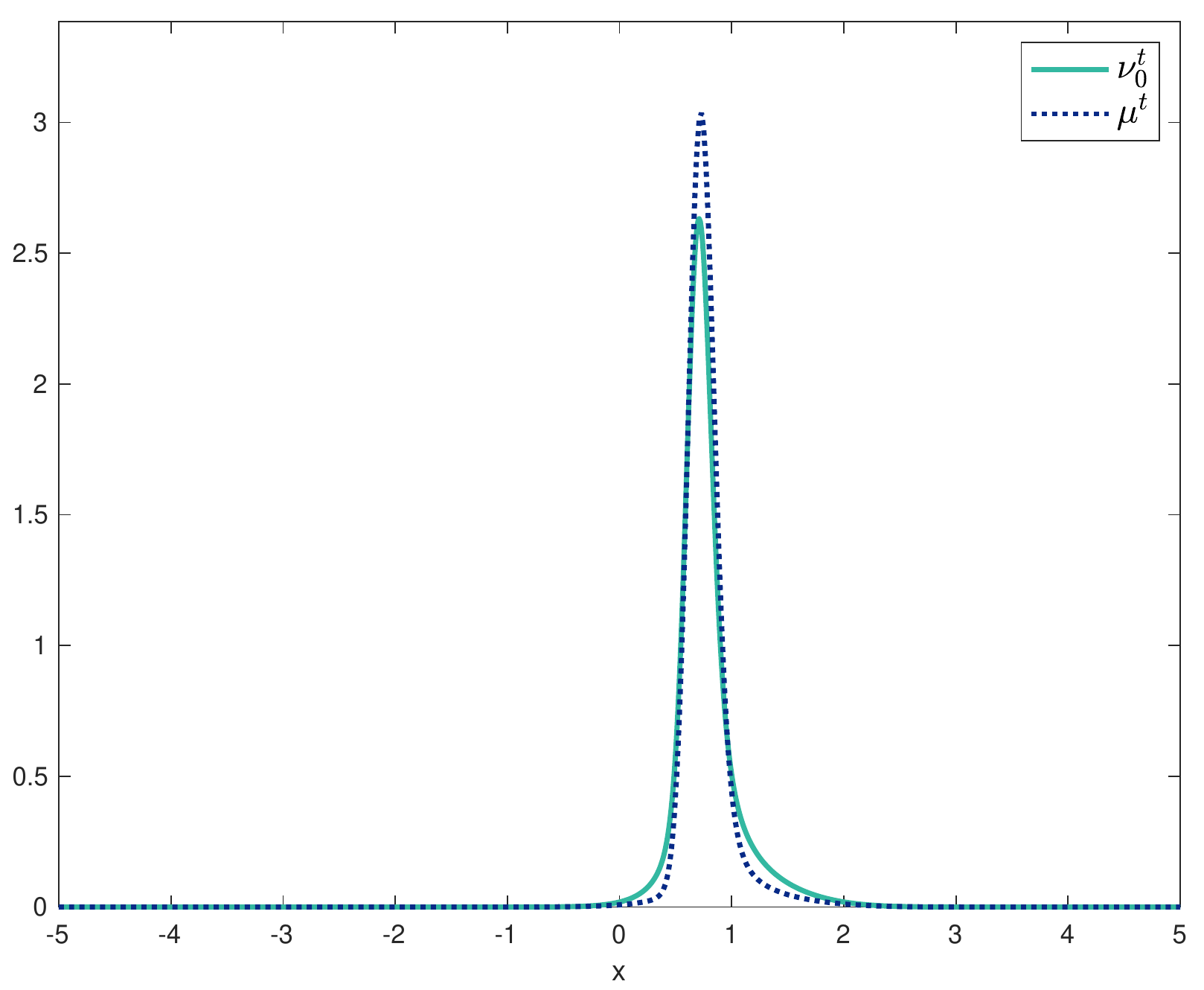}}
\subfigure{\includegraphics[width=5cm, height=4cm]{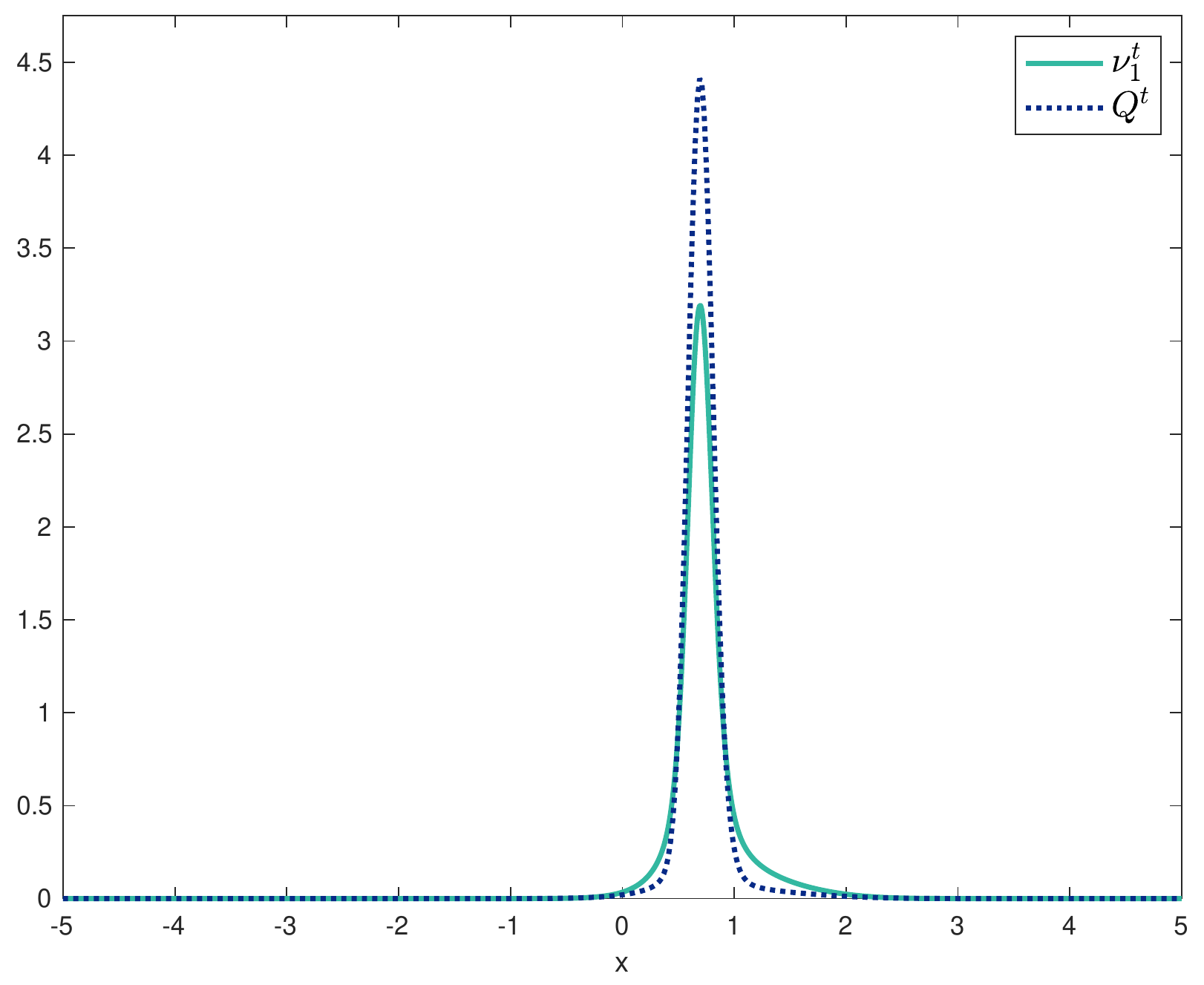}}	\\
\hspace{0.5cm}
\subfigure{\includegraphics[width=4.5cm, height=4.2cm]{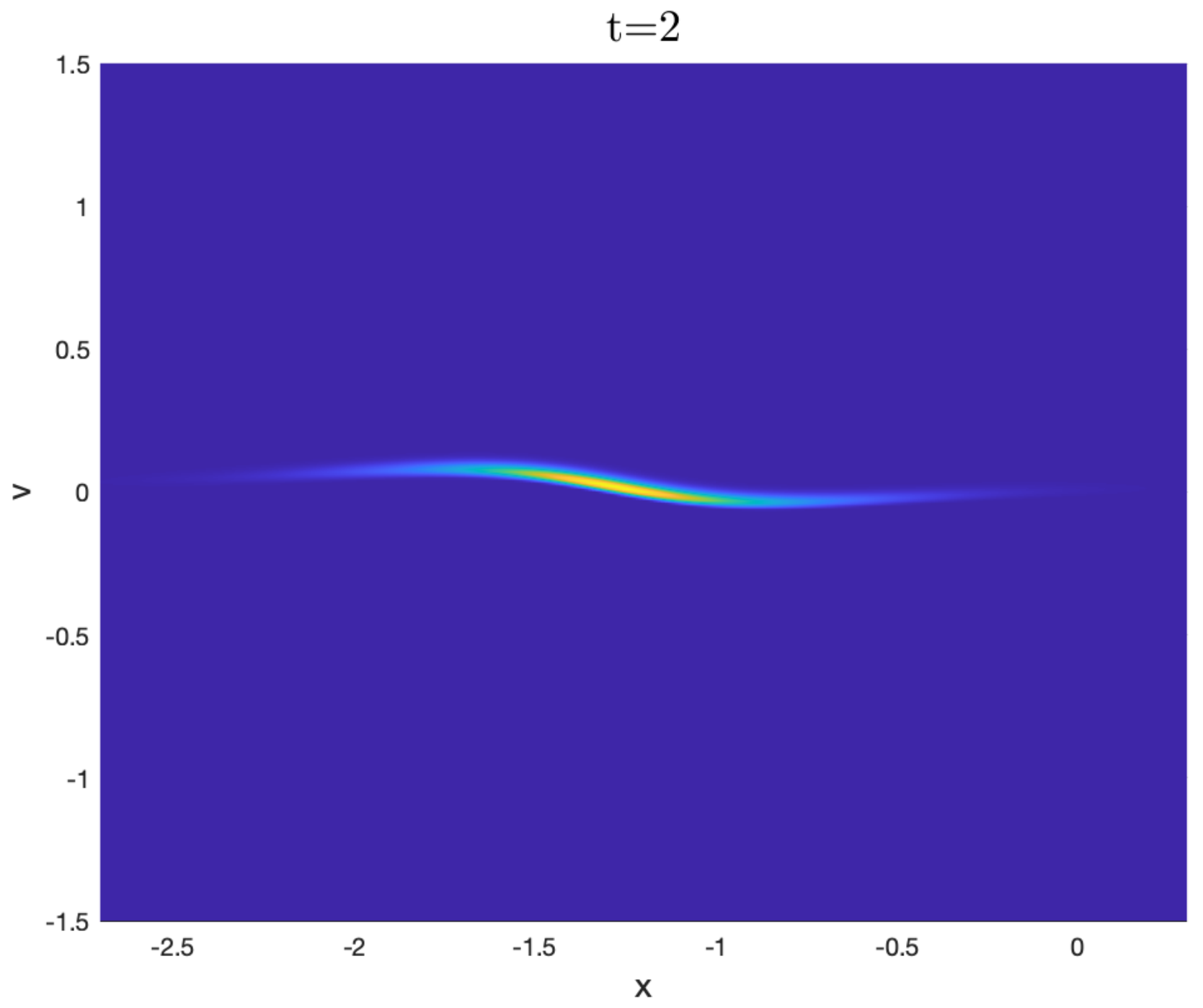}}
\subfigure{\includegraphics[width=5cm, height=4cm]{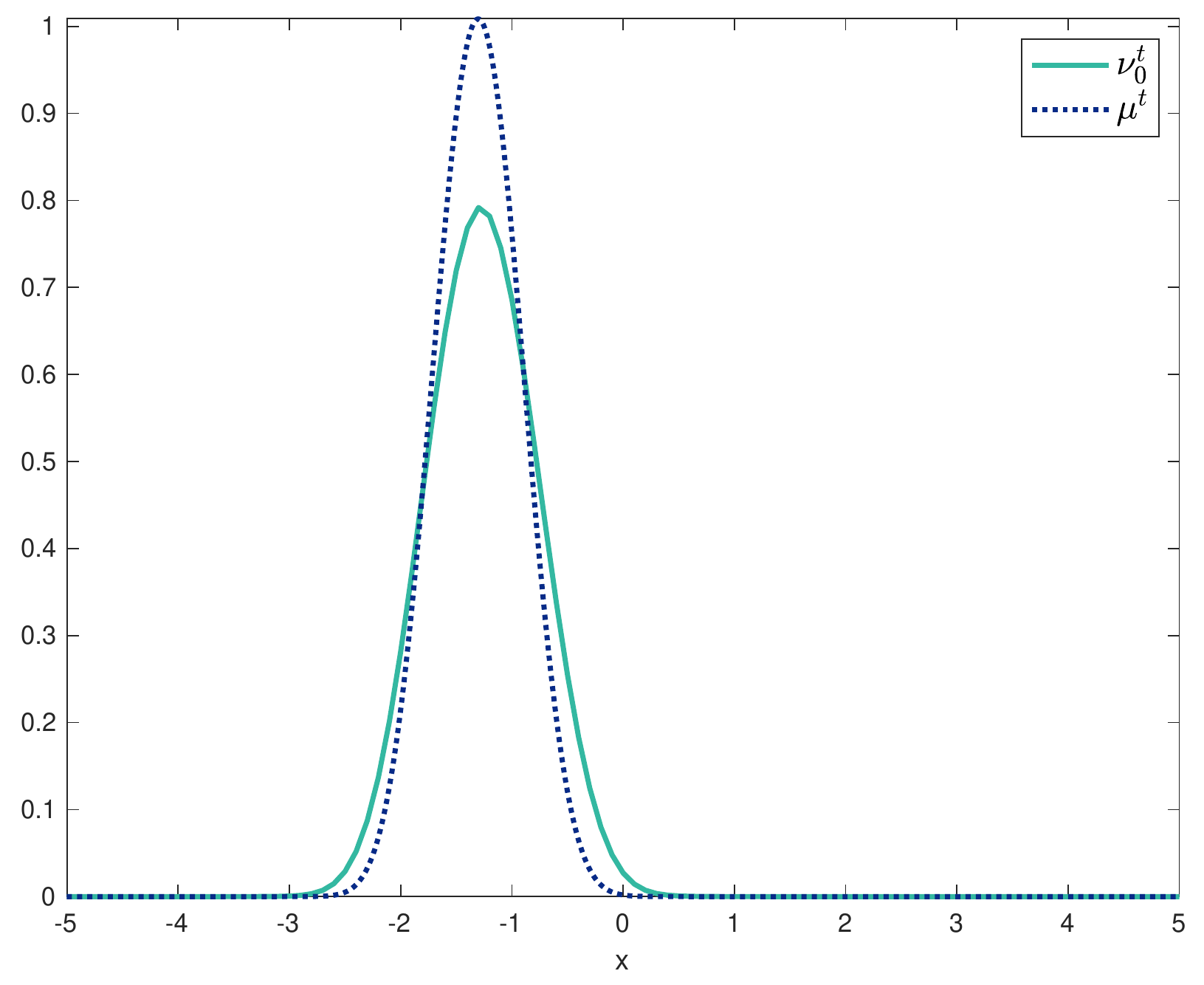}}
\subfigure{\includegraphics[width=5cm, height=4cm]{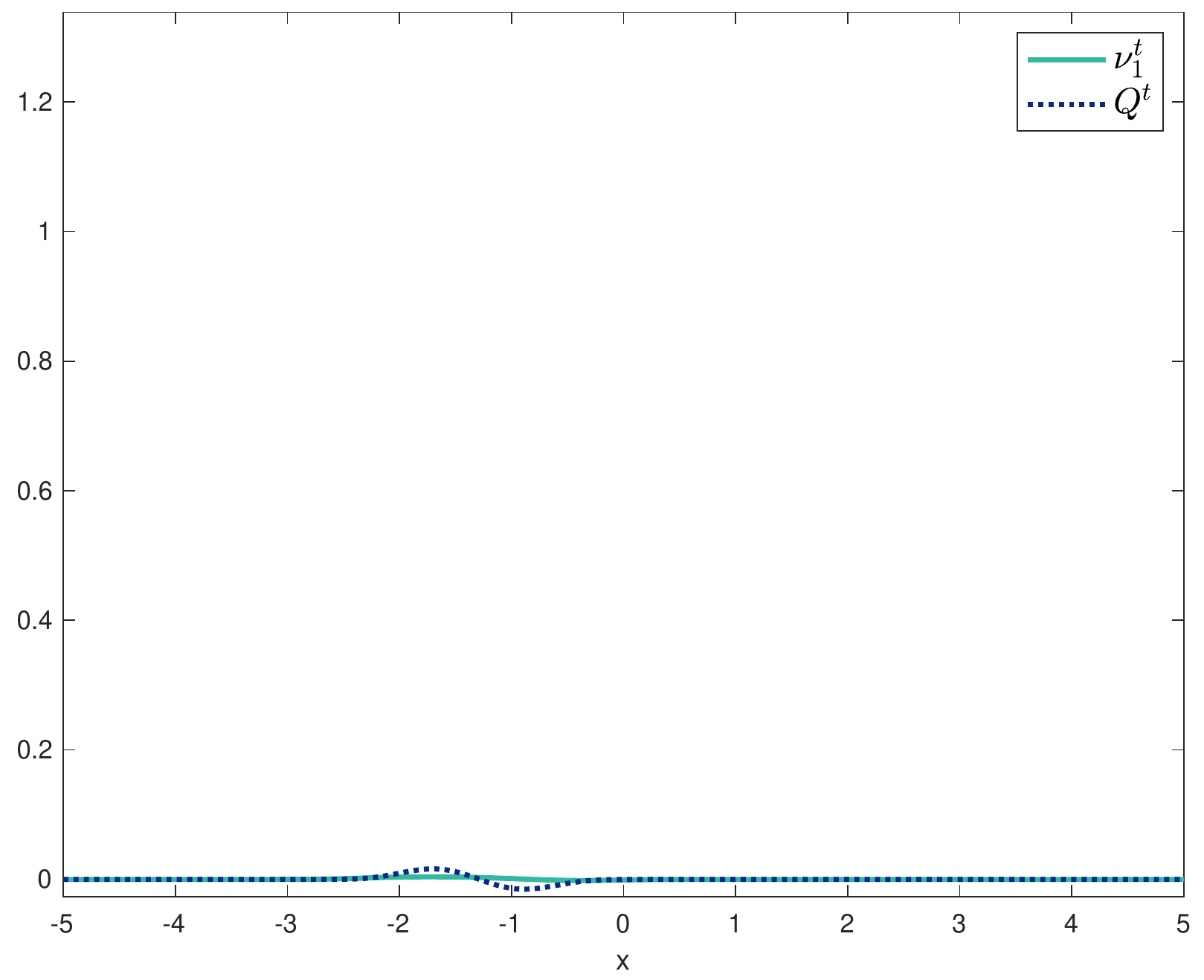}}	
	\caption{Test 7: approximation of monokinetic initial data. Comparison between zero-th order moment of the solution of (V) and $\mu^t$ (second column),  first order moment of the solution of (V) and $Q^t$ (third column): without chemotaxis (first line: $\eta=0$, $\alpha=0$ at $t=2$), with chemotaxis (second line: $\eta=0.2$, $\alpha=0$ at $t=2$), with chemotaxis and damping (third line: $\eta=0.2$, $\alpha=2$ at $t=2$).}
	\label{fig:monokinetic}
\end{figure}

Figure \ref{fig:monokinetic_initialcondition } shows the plot of $\rho^0$ in the phase space, and the profile of initial data for $(E)$ (blue dotted line), which are defined in \eqref{mom_data} and correspond to Vlasov moments (cyan line).
Figure \ref{fig:monokinetic} shows the behavior, at the same time instant $t=2$, of $\nu^t_0$ and $\nu^t_1$ together with $\mu^t$ and $Q^t$ respectively, both without (first line) and with chemotaxis (second and third line).
The evolution of Vlasov dynamics preserves the monokinetic (approximated) initial structure, and Vlasov moments and Euler solutions are in good agreement. 
Numerical evidence shows that adding chemotaxis, the approximated monokinetic-like structure is not preserved for long time.
As already observed at particle level, the presence of a chemical gradient influencing the dynamics tends to concentrate the density profile in the same point, and to asymptotically let the velocities converge to zero. This can be seen in the phase space (first column, second line). Moreover, introducing also a damping term, the momentum tends to zero faster (third line).
Focusing on the Vlasov moments, we observe that as soon as the monokinetic structure is lost, the agreement between Vlasov moments and Euler solution is worse. In particular, we observe a blow up of the density profile.
In conclusion, the monokinetic property seems to play a crucial role. For that reason we exploit with several tests the case of largely non monokinetic initial data.

\newpage

\begin{figure}[h!]
	\centering
		\subfigure[]{\includegraphics[width=0.30\textwidth]{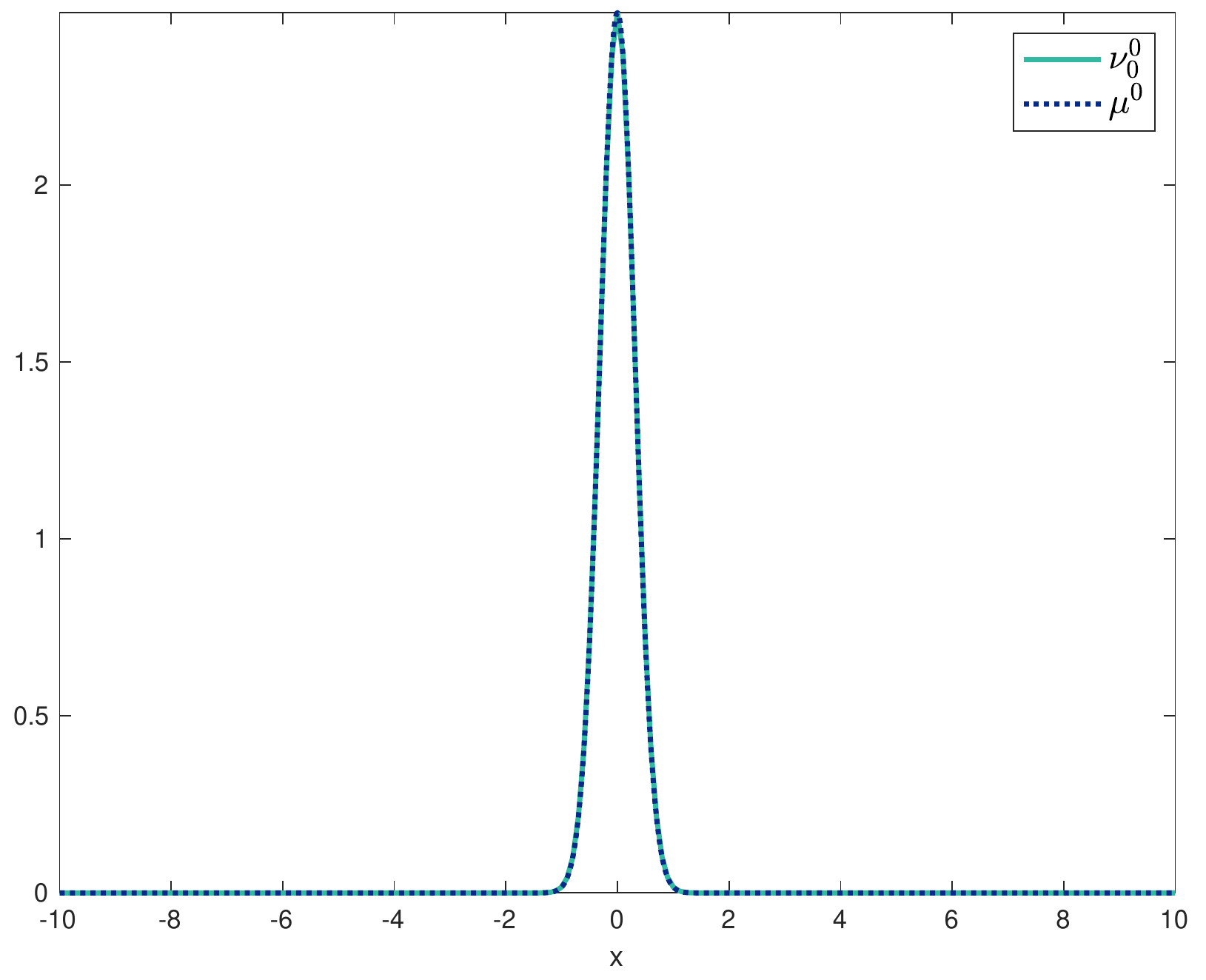}}
	\subfigure[]{\includegraphics[width=0.30\textwidth]{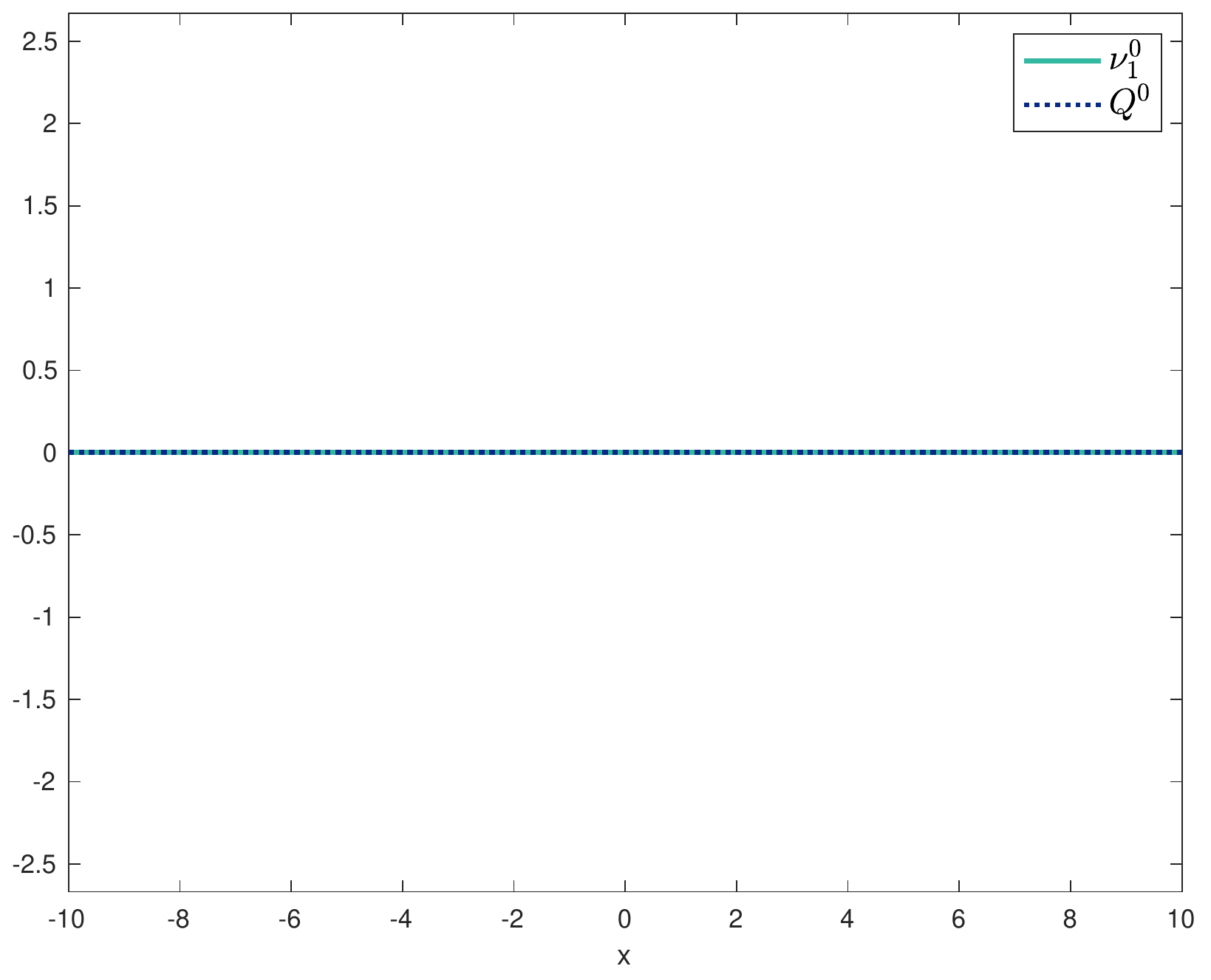}}
		\vspace{-0.3cm}
	\caption{Test8, Test9, Test10: a) initial density $\mu^0$ and b) initial momentum $Q^0$ (here null).  }
	\label{fig:VE_initial_condition}
\end{figure}

\textit{Non monokinetic initial data: }
As non monokinetic initial data, we consider the distribution function already introduced in \eqref{rho0}.
Figure \ref{fig:VE_initial_condition} shows initial conditions for the following Test 8 and Test 9, as defined in \eqref{mom_data}.

\begin{figure}[h!]
	\centering
		\subfigure[]{\includegraphics[width=0.30\textwidth]{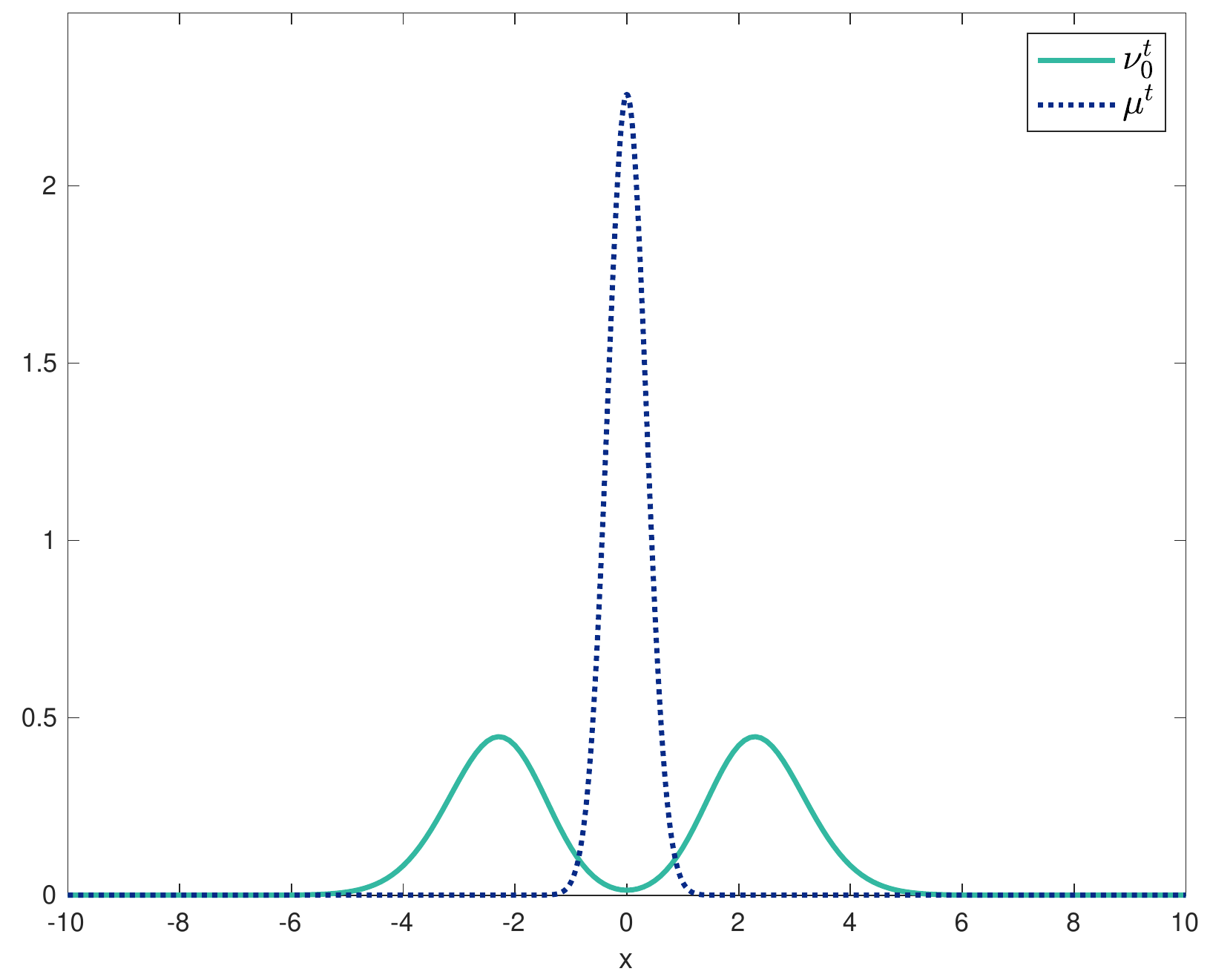}}
	\subfigure[]{\includegraphics[width=0.30\textwidth]{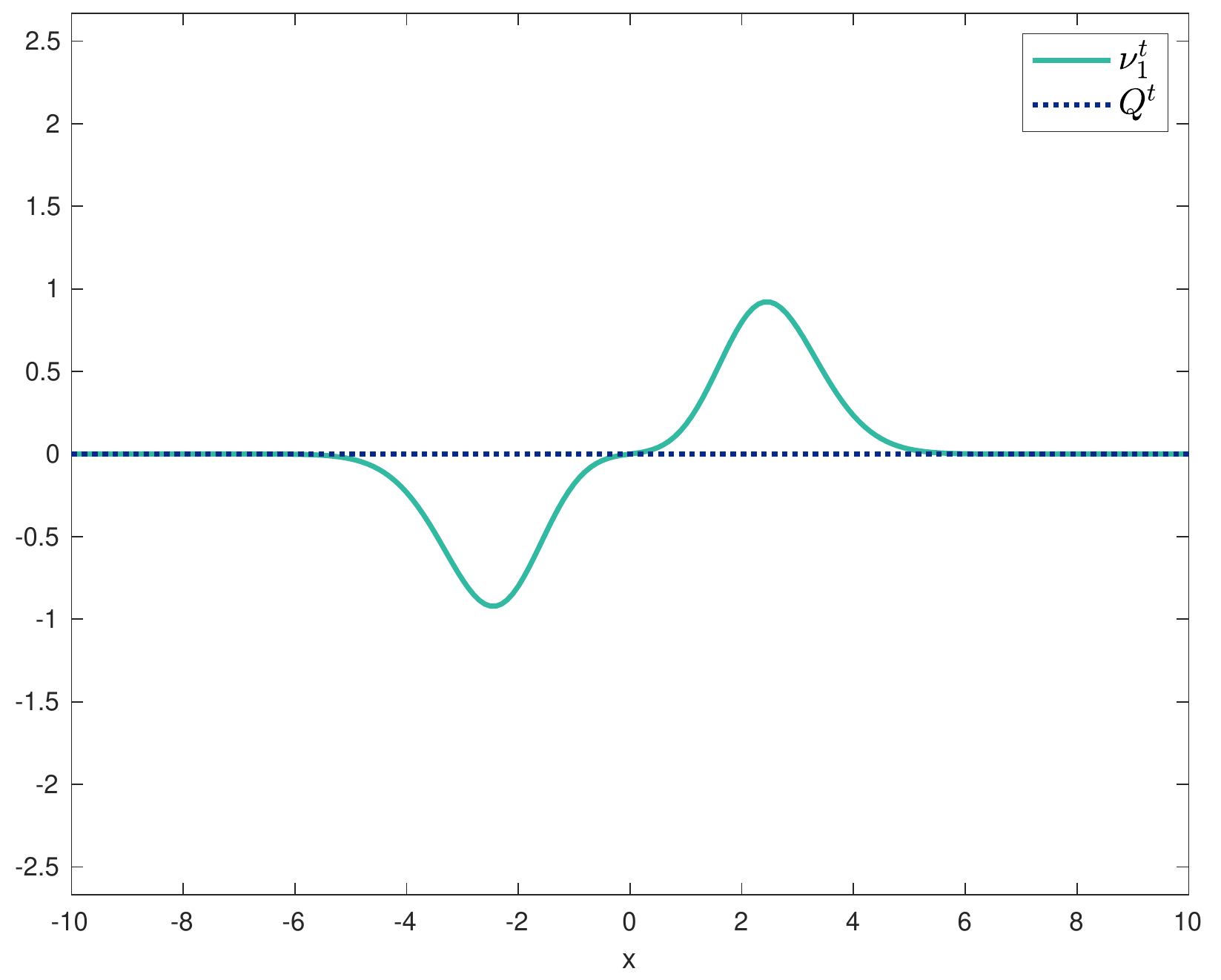}}\\
	\subfigure[]{\includegraphics[width=0.30\textwidth]{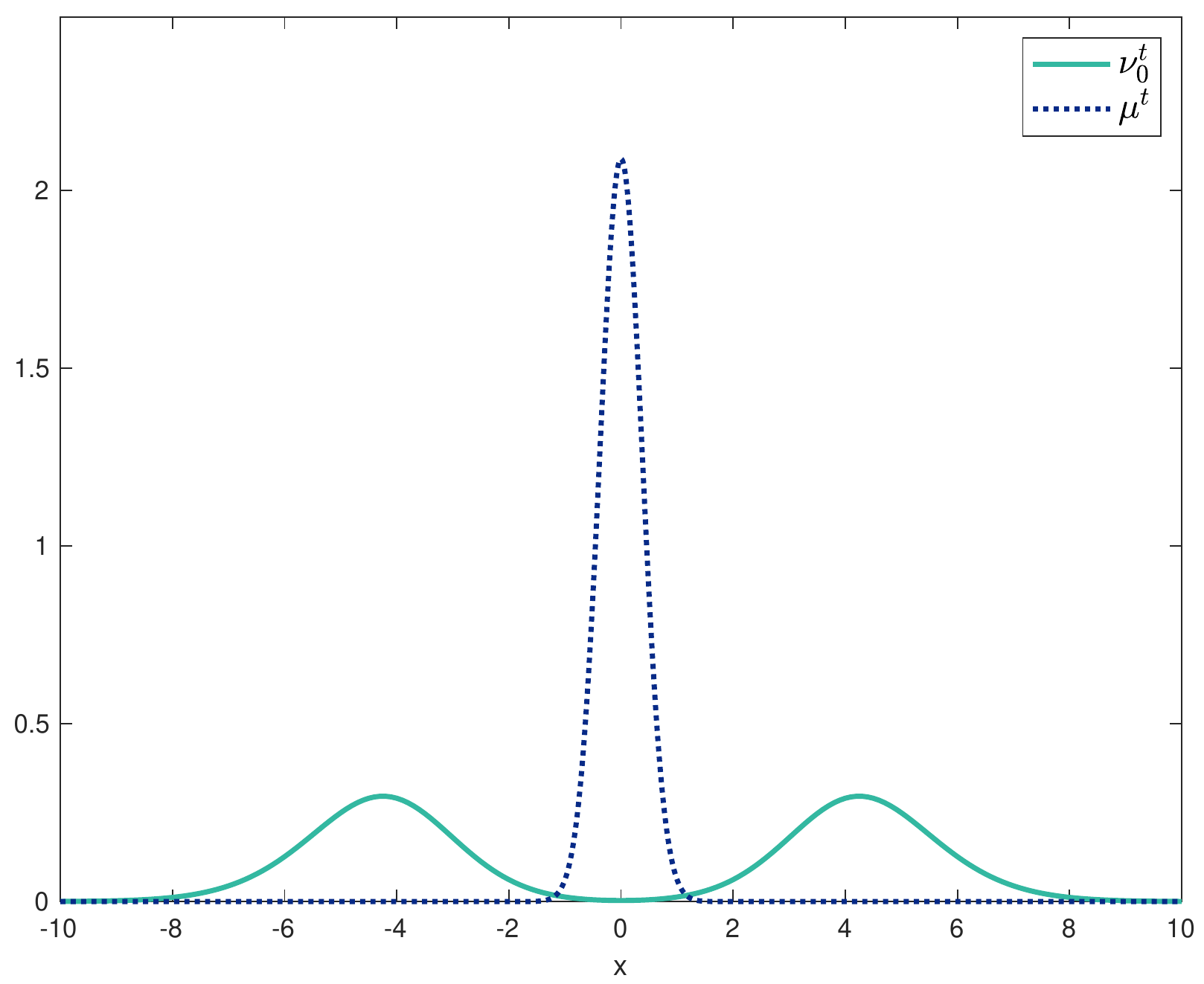}}
	\subfigure[]{\includegraphics[width=0.30\textwidth]{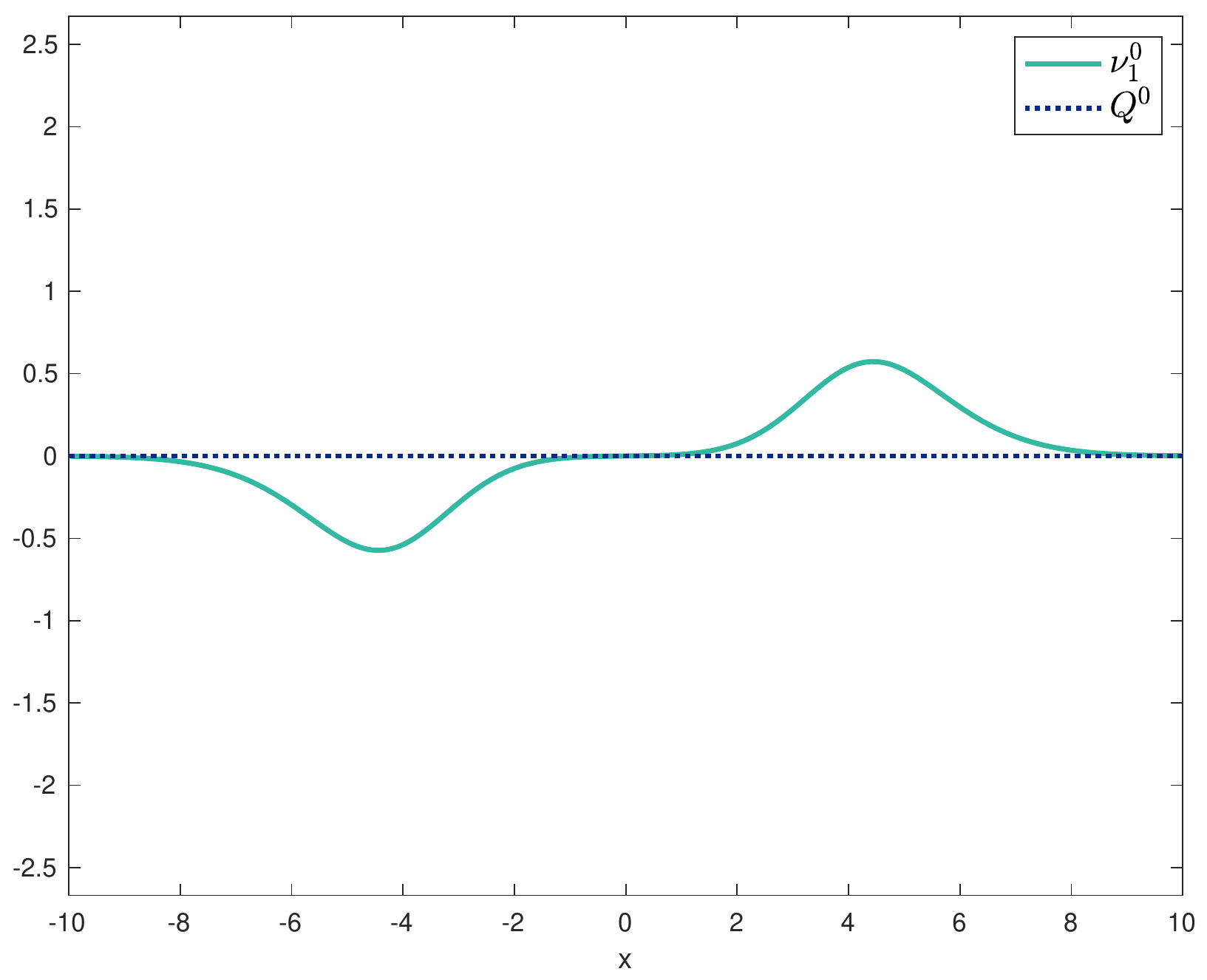}}\\
	\subfigure[]{\includegraphics[width=0.30\textwidth]{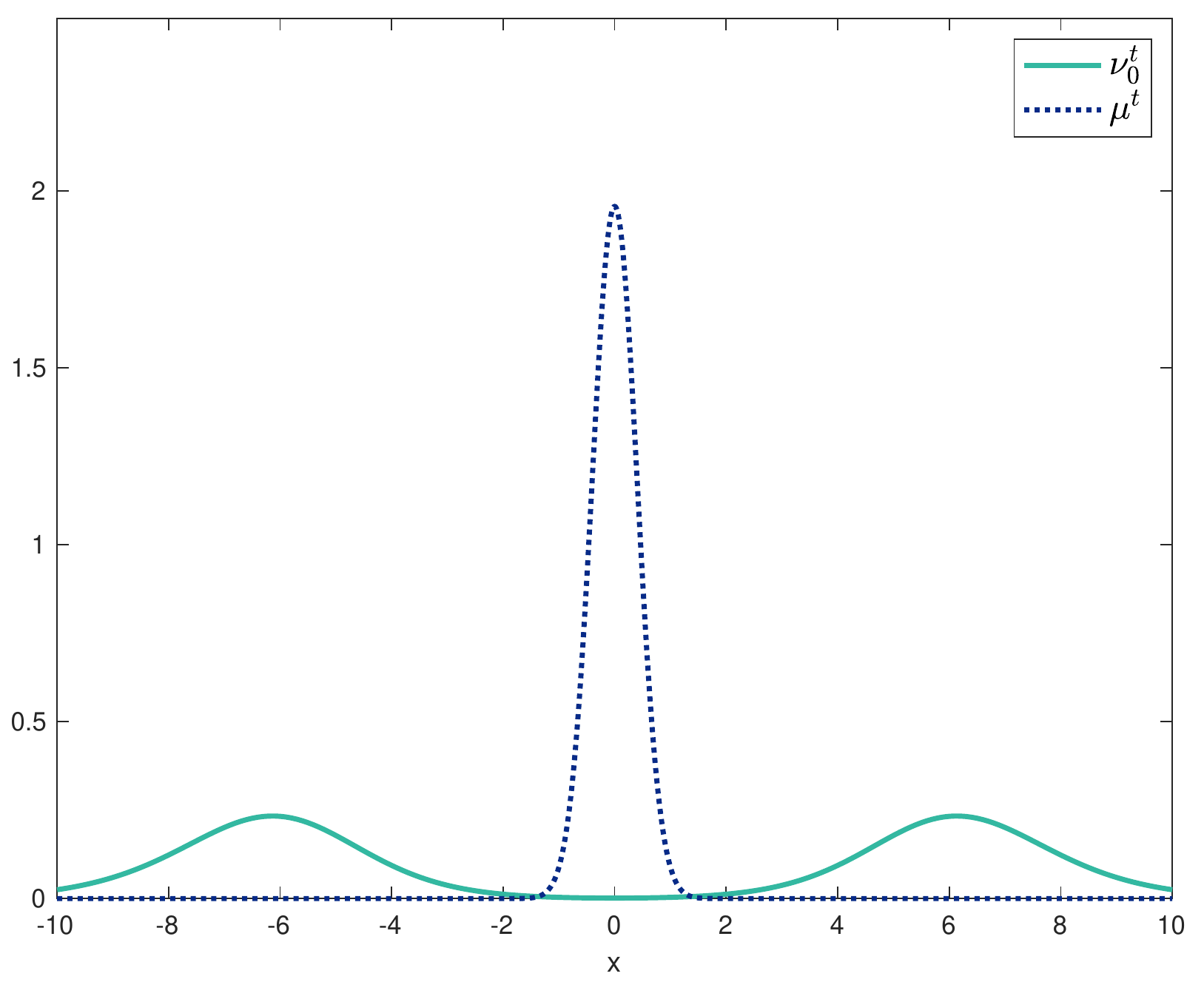}}
	\subfigure[]{\includegraphics[width=0.30\textwidth]{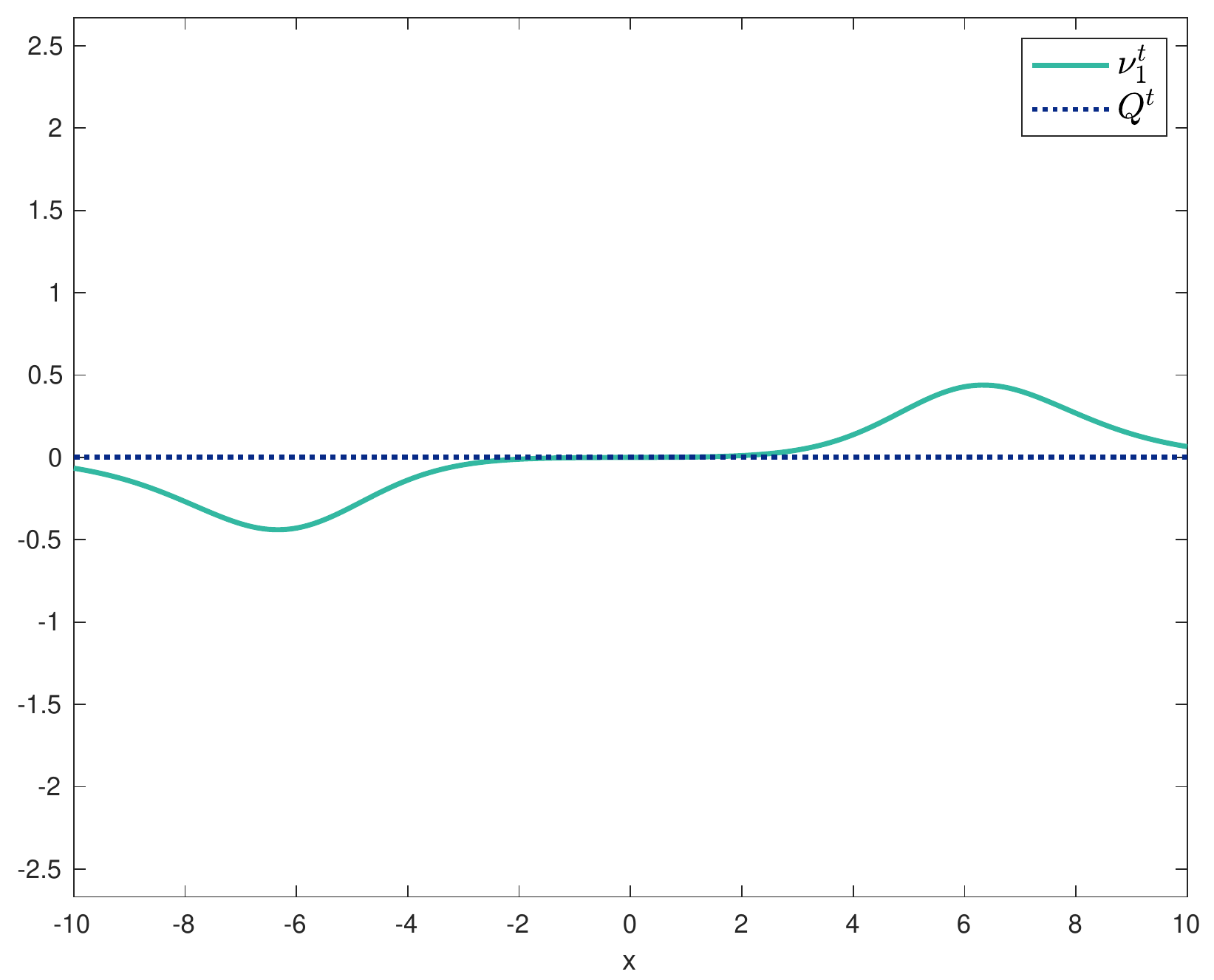}}
	\vspace{-0.3cm}
	\caption{Test 8: Comparison between zero-th order moment of the solution of (V) and $\mu^t$,  first order moment of the solution of (V) and $Q^t$  at a)-b) $t=1$, c)-d) $t=2$, e)-f) $t=3$. Here $\eta=0$ (no chemotaxis), $\alpha=0$.  }
	\label{fig:VE_nochemo}
\end{figure}

In Test 8 we compare the moments $\nu_0^t$ and $\nu_1^t$ of the solution of \eqref{Vlasov_nochemo}-\eqref{rho0} with the solution of system \eqref{Euler1D_nochemo} with initial conditions \eqref{mom_data}. 
Since $u^0(x)=0$, the velocity $u^t$ remains zero, being a stationary solution for the Euler system.
On the contrary, the first order moment profile shows a separation due to the effect of speeds.
Being far from the monokinetic case, in Vlasov dynamics $v$ is a variable which is not a priori related to $x$, whereas in Euler the velocity field depends on the position. 
As shown in Figure \ref{fig:VE_nochemo}, there is no agreement between Vlasov moments and Euler solutions starting with that initial condition.

\begin{figure}[h!]
	\centering
	\subfigure[]{\includegraphics[width=0.30\textwidth]{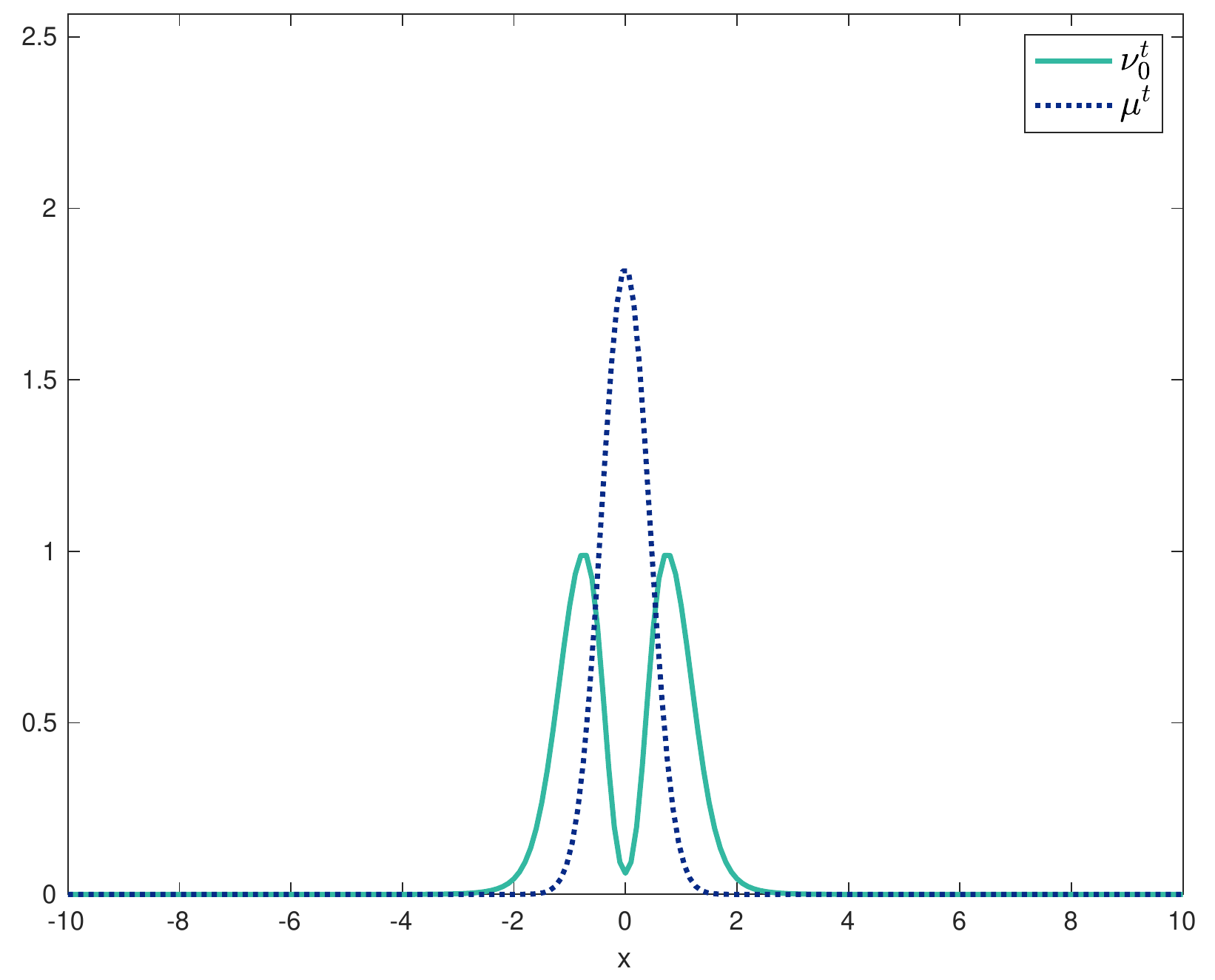}}
	\subfigure[]{\includegraphics[width=0.30\textwidth]{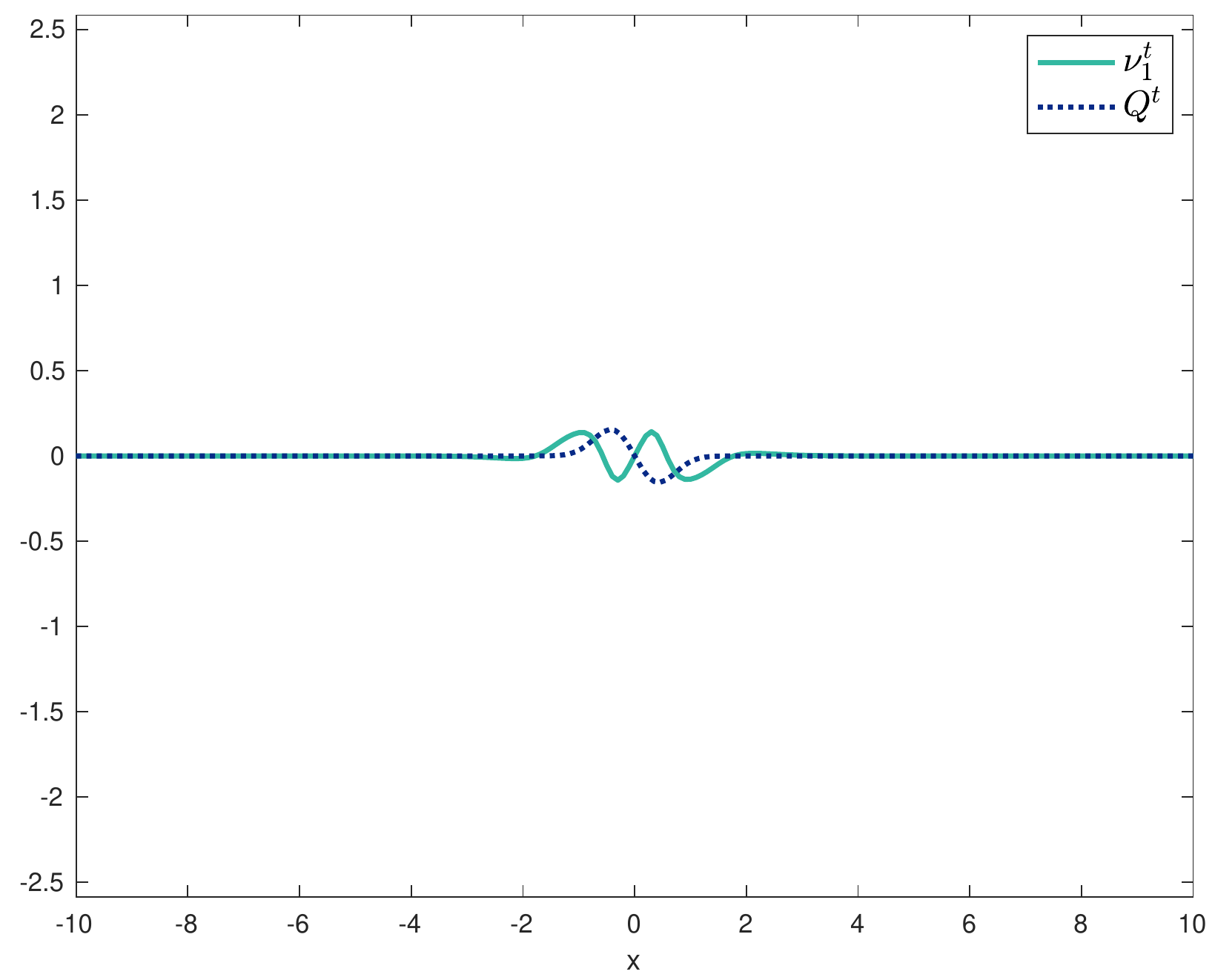}}\\
	\subfigure[]{\includegraphics[width=0.30\textwidth]{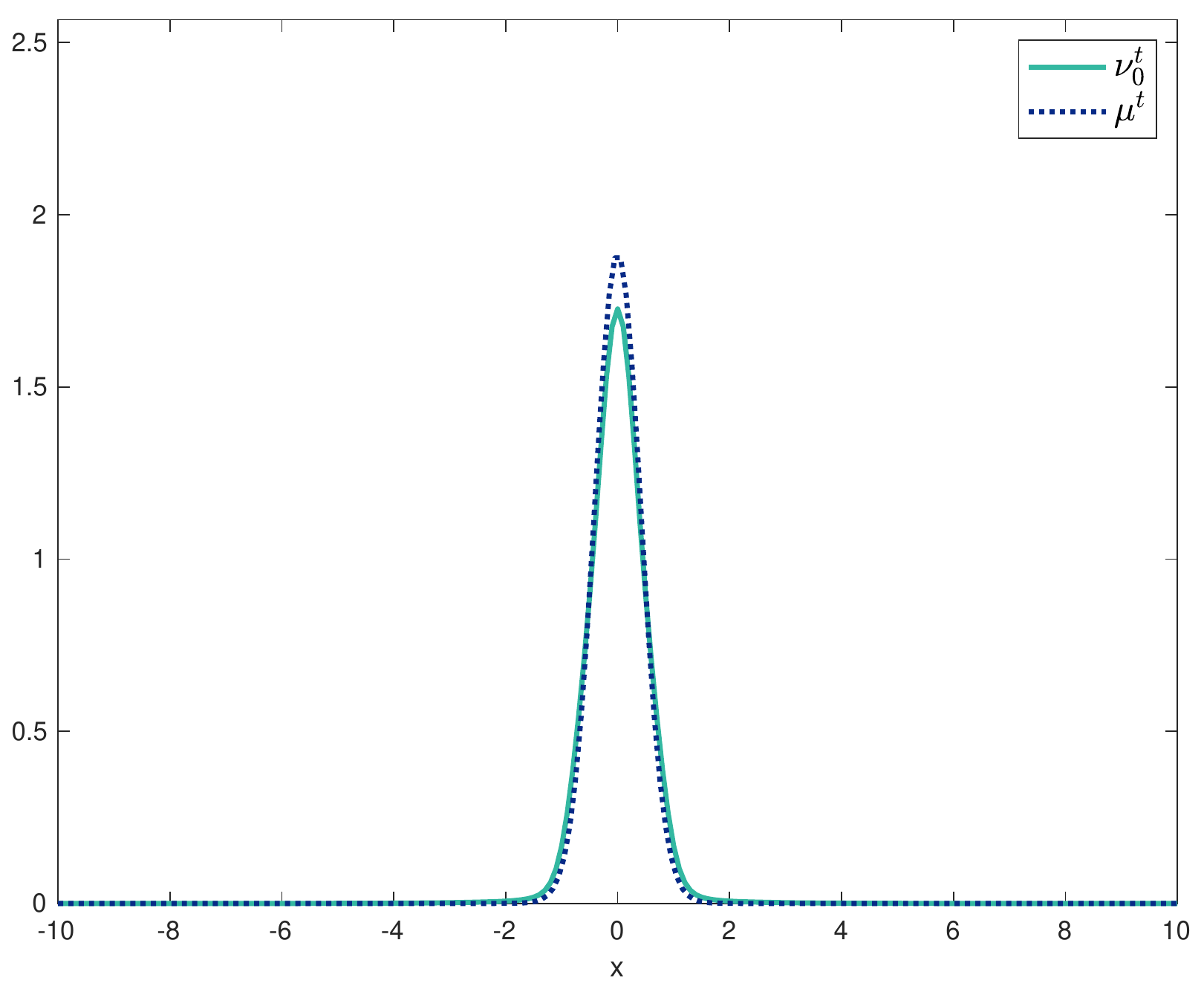}}
	\subfigure[]{\includegraphics[width=0.30\textwidth]{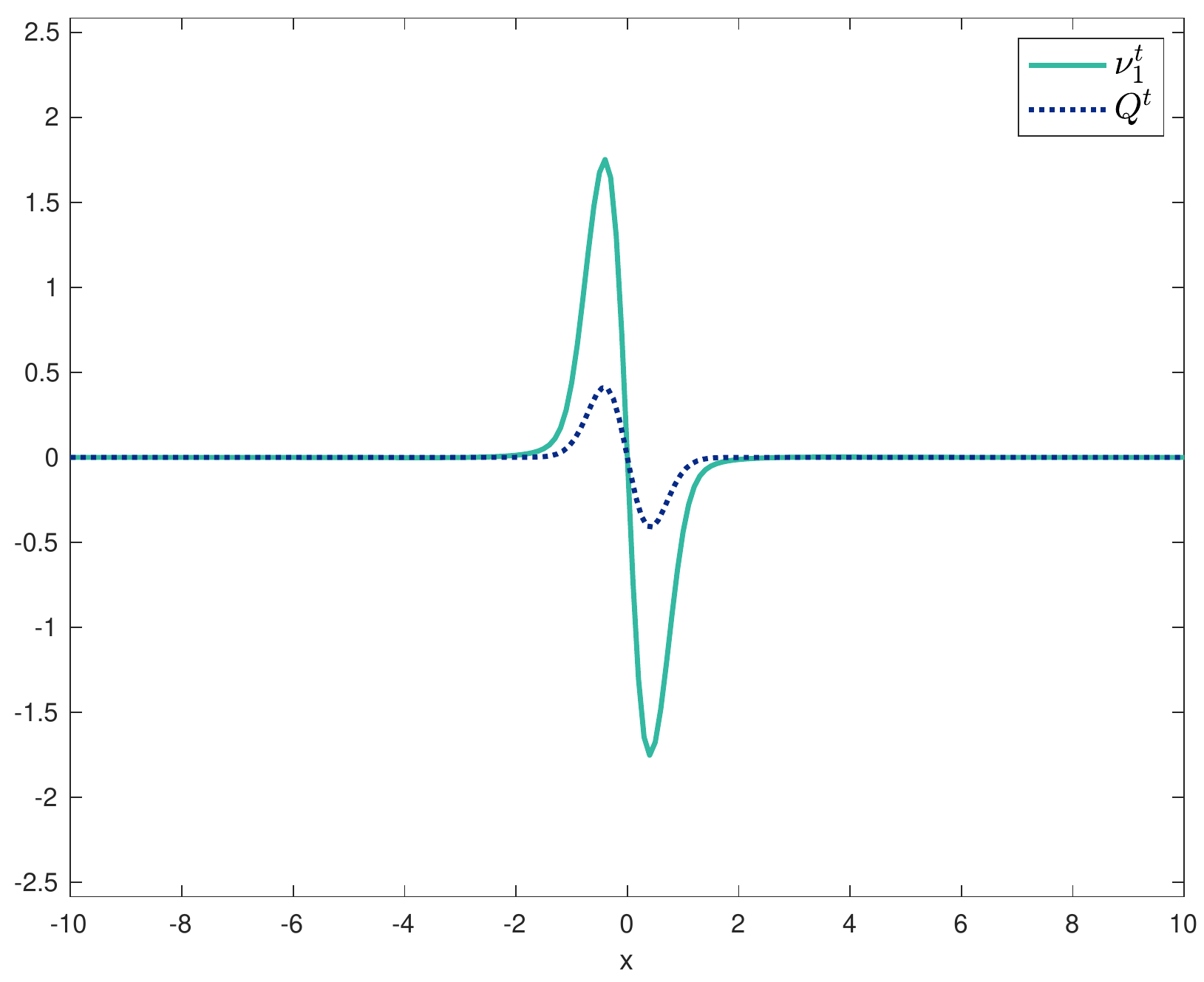}}\\
	\subfigure[]{\includegraphics[width=0.30\textwidth]{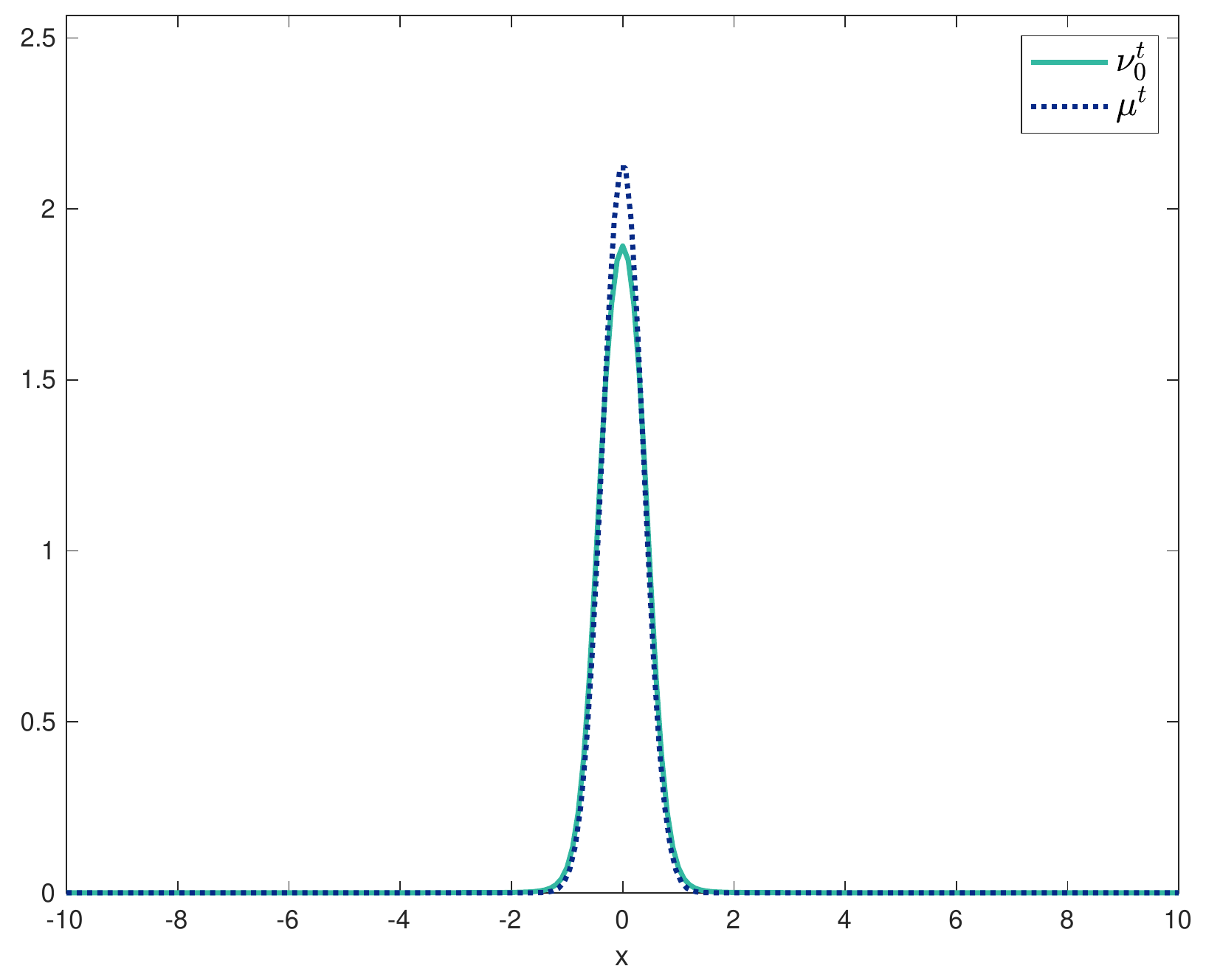}}
	\subfigure[]{\includegraphics[width=0.30\textwidth]{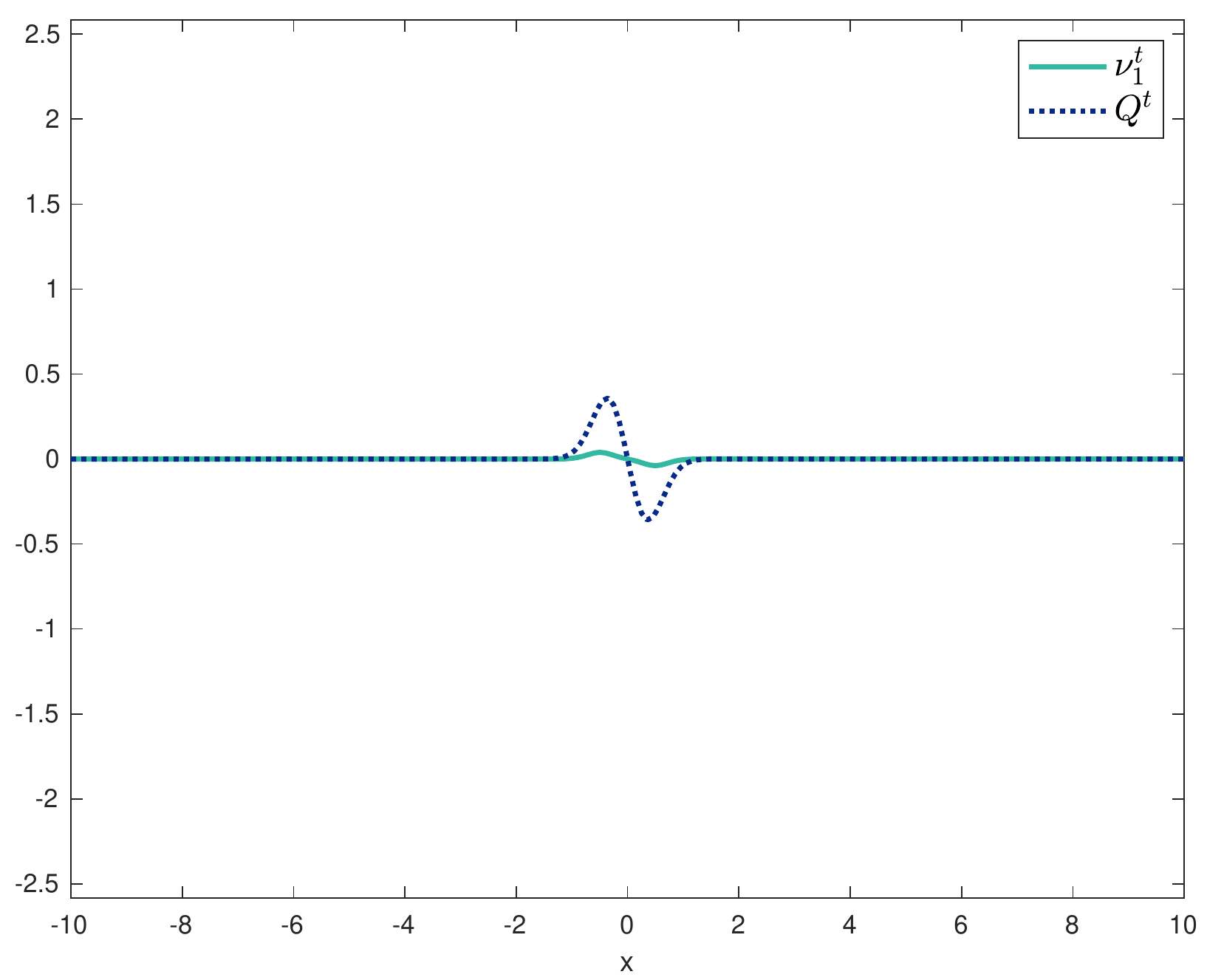}}
	\caption{Test 9: Comparison between zero-th order moment of the solution of (V) and $\mu^t$,  first order moment of the solution of (V) and $Q^t$   at a)-b) $t=0.5$, c)-d) $t=1$, e)-f) $t=7$.  Here $\eta=3$, $\alpha=0$.}
	\label{fig:VE_chemo}
\end{figure}

\begin{figure}[h!]
	\centering
	\subfigure[]{\includegraphics[width=0.30\textwidth]{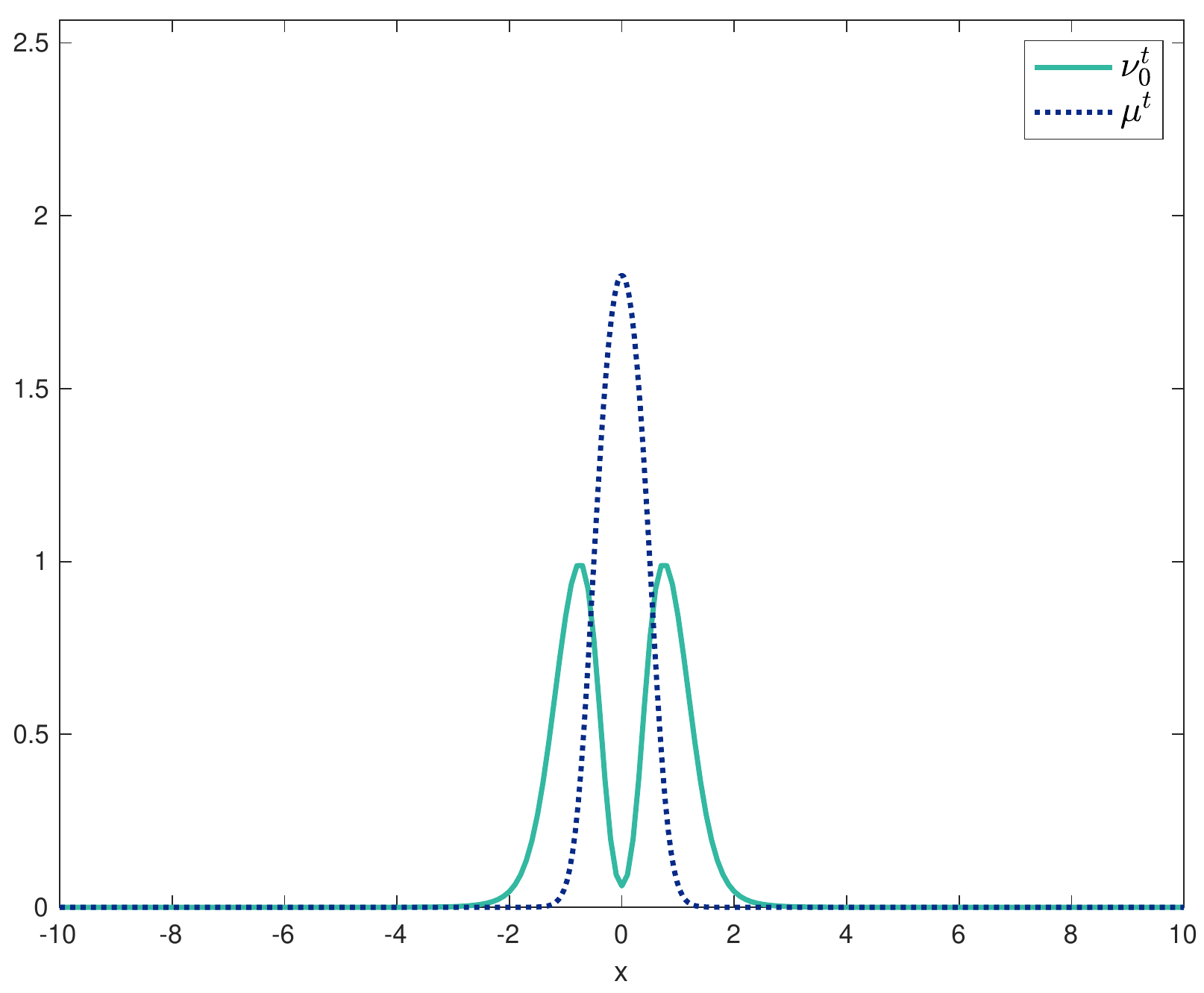}}
	\subfigure[]{\includegraphics[width=0.30\textwidth]{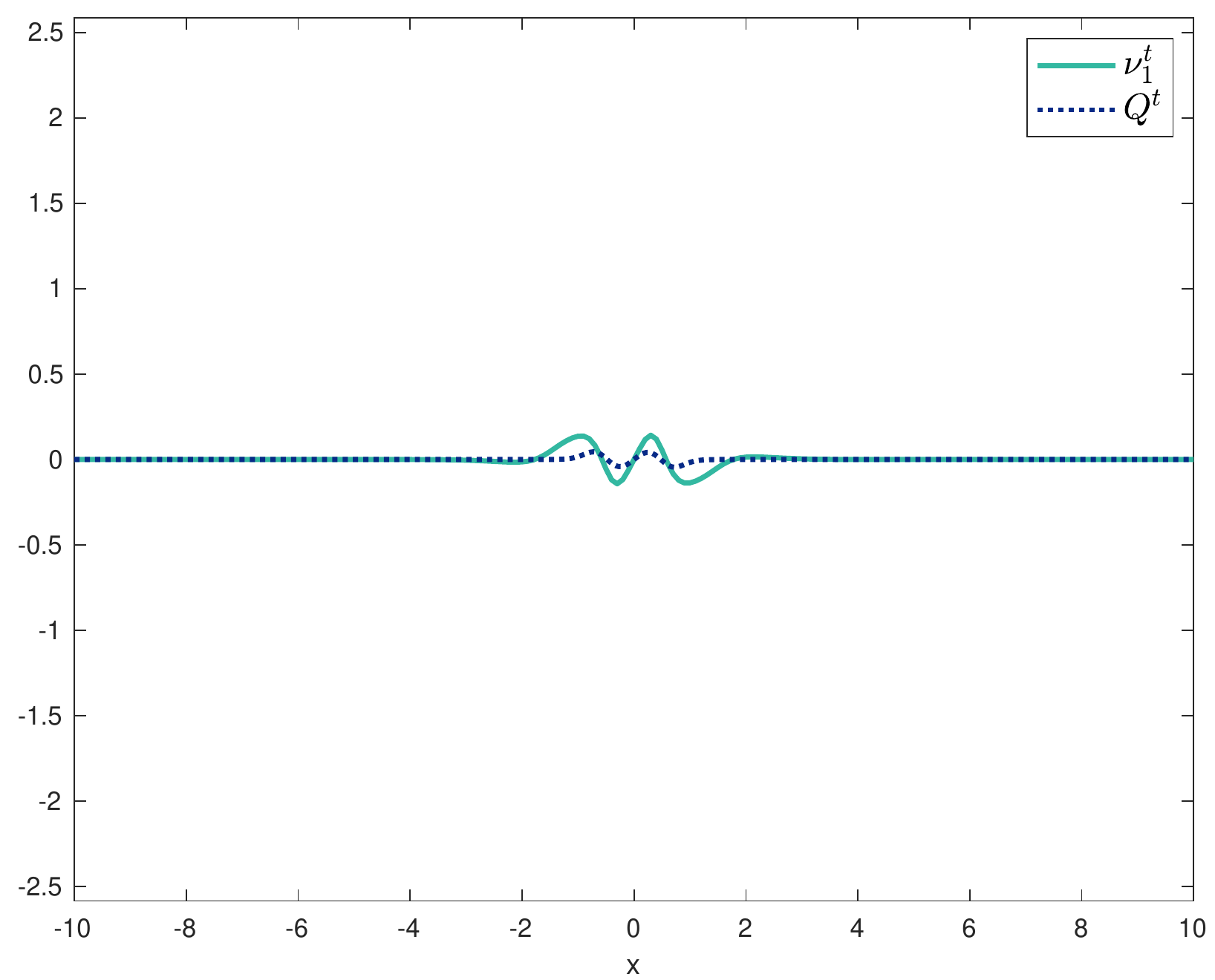}}\\
	\subfigure[]{\includegraphics[width=0.30\textwidth]{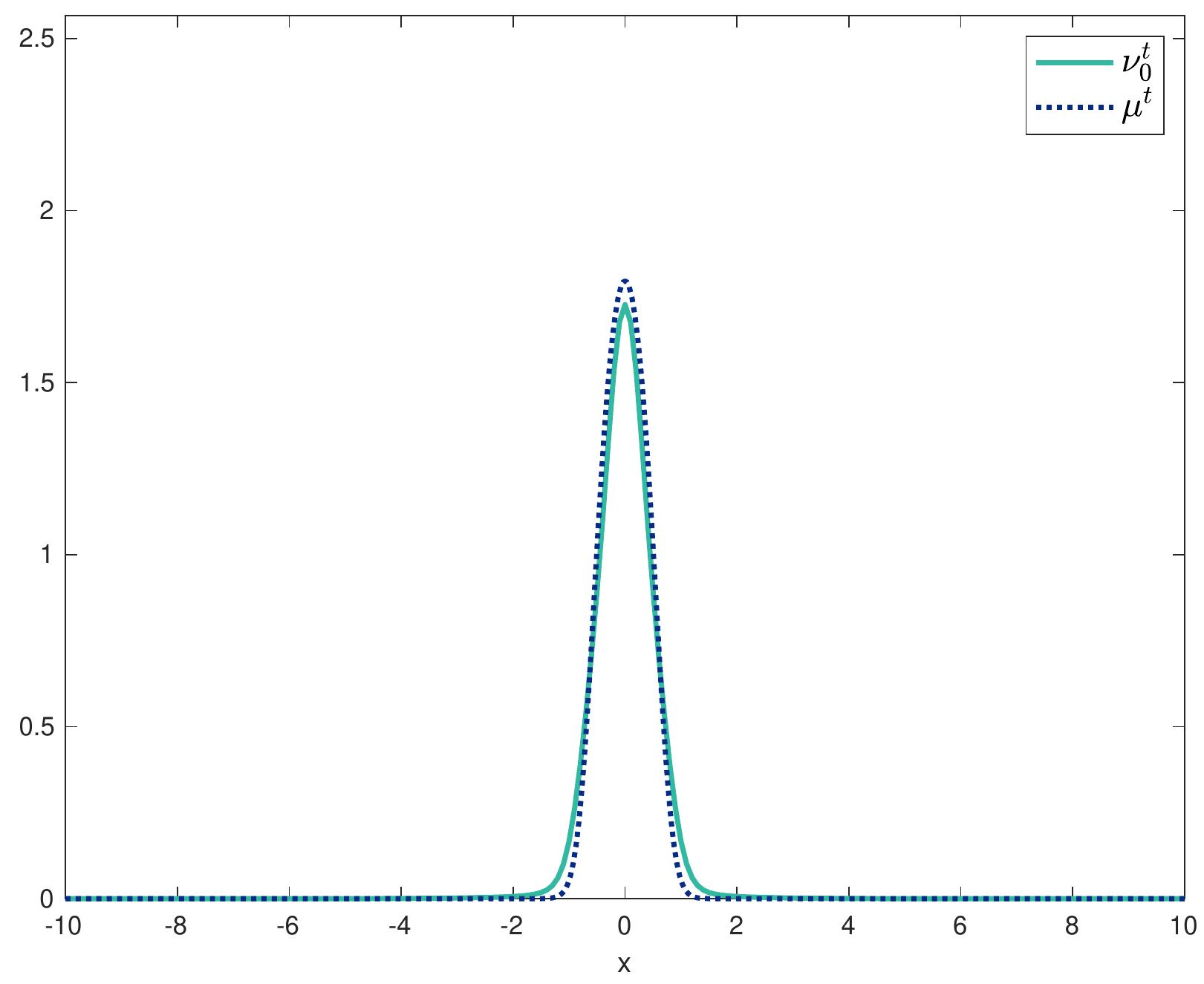}}
	\subfigure[]{\includegraphics[width=0.30\textwidth]{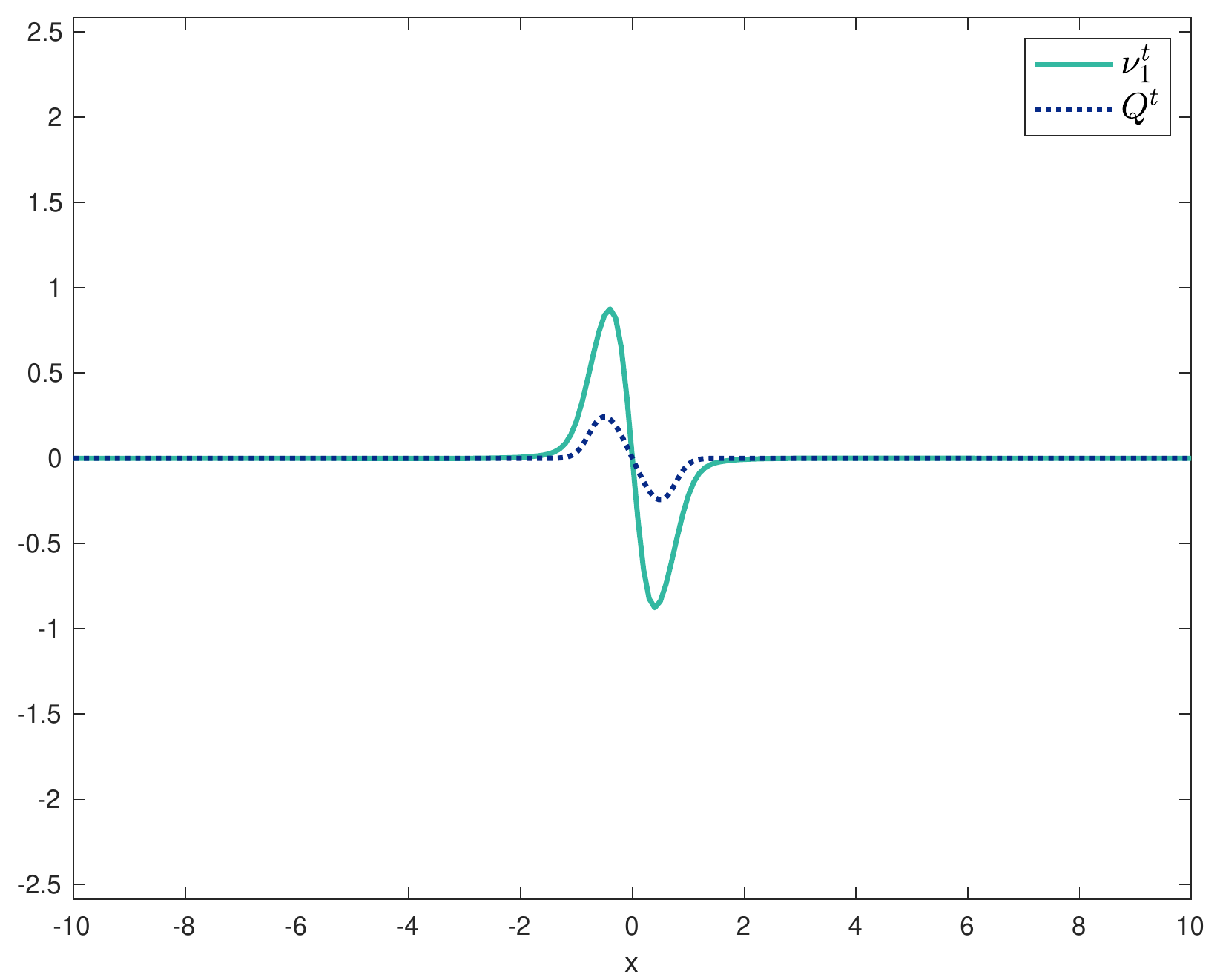}}\\
	\subfigure[]{\includegraphics[width=0.30\textwidth]{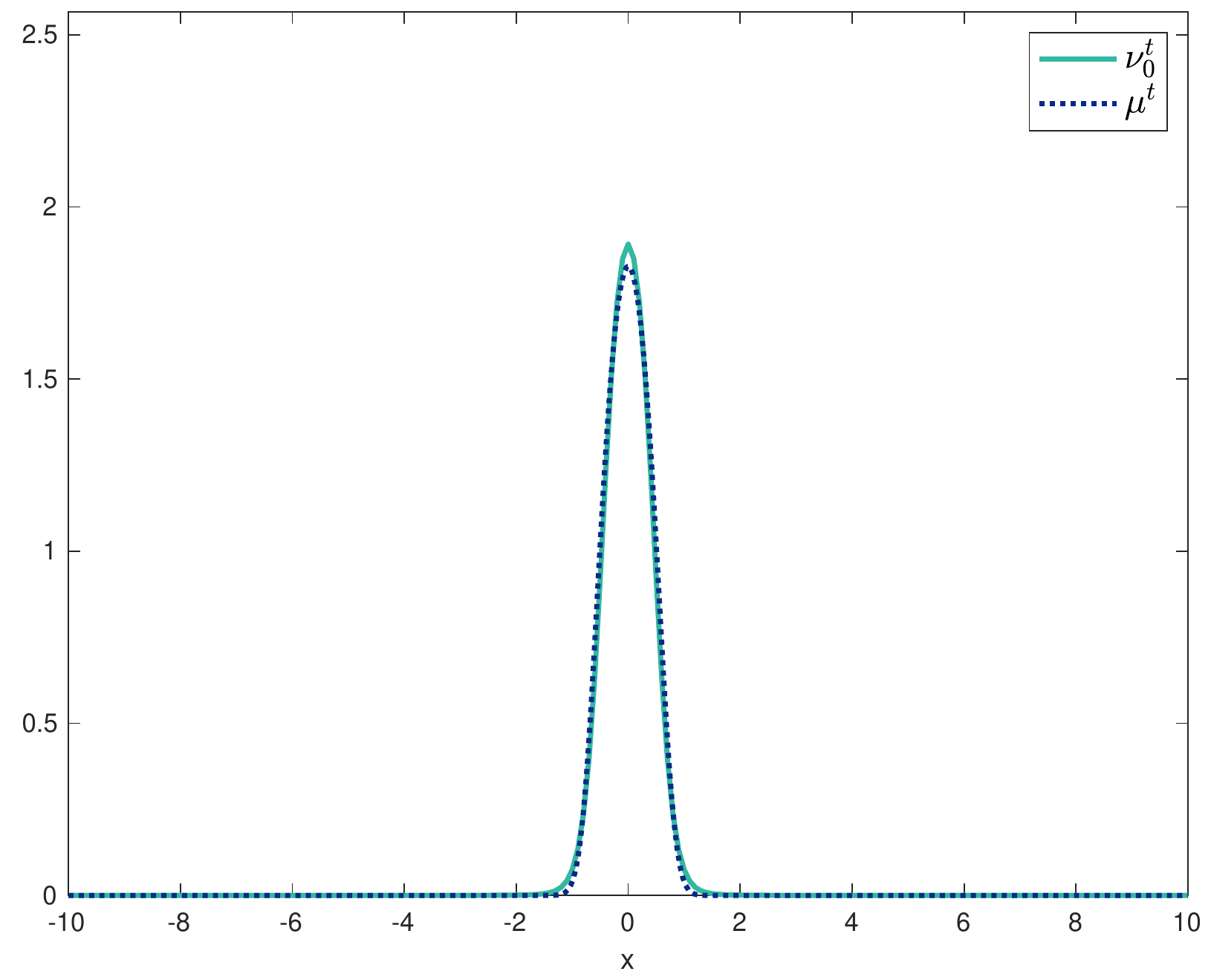}}
	\subfigure[]{\includegraphics[width=0.30\textwidth]{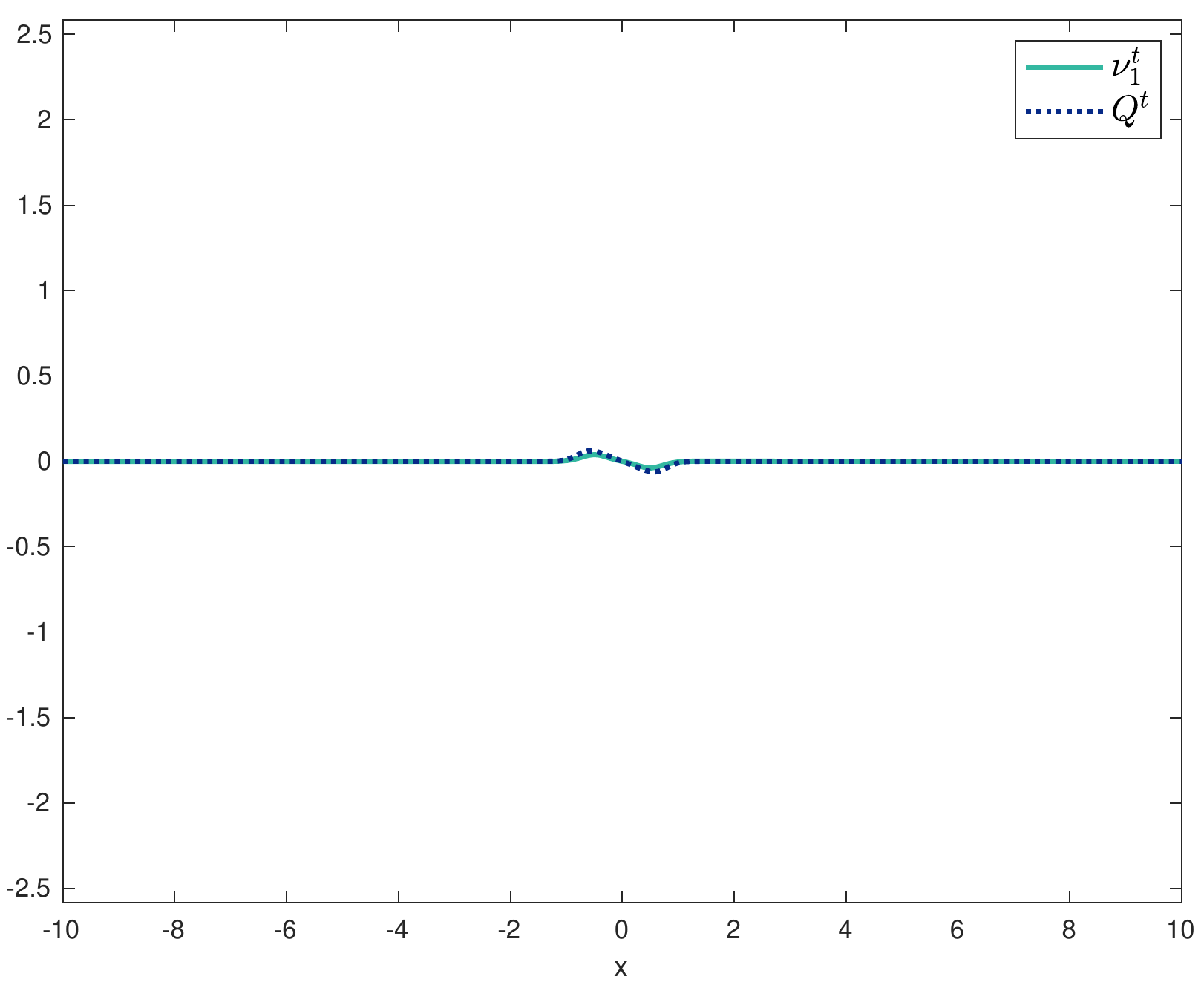}}
	\caption{Test 9 adding damping: Comparison between zero-th order moment of the solution of (V) and $\mu^t$,  first order moment of the solution of (V) and $Q^t$   at a)-b) $t=0.5$, c)-d) $t=1$, e)-f) $t=7$. Here $\eta=3$, $\alpha=0.8$.}
	\label{fig:chemo_damping}
\end{figure}

Let now solve \eqref{Euler1D_damping}, hence considering also a chemotactic effect.
Figure \ref{fig:VE_chemo} shows three different snapshots of the numerical simulation performed.
At the beginning, we observe that the spreading of the speed in the density, solution of Vlasov equation, is stronger than the chemistry, and two separate peaks appear. With respect to the previous case, the presence of the chemotactic gradient prevent the separation observed for longer time, leading to the convergence to a unique bump. 
Concerning the first order moment (Vlasov) and the momentum (Euler), at the beginning they are both null. Due to the presence of chemicals, $u^t(x)=0$ is no longer a stationary solution for Euler.
In particular in presence of a damping term ($\alpha=0.8$) first order moment and momentum asymptotically converge to zero, as shown in Figure \ref{fig:chemo_damping}. 
We observe that, as time grows, the agreement between Vlasov moments and Euler solution is almost recovered especially for the density. 
Numerical evidence hence shows that, far from the monokinetic case, the dissipation of energy induced by damping term helps in recovering a good agreement. In fact, the role of the second-order moment, that we are here neglecting, is reduced and almost negligible.

\begin{figure}[h!]
	\centering
	\subfigure[]{\includegraphics[width=0.30\textwidth]{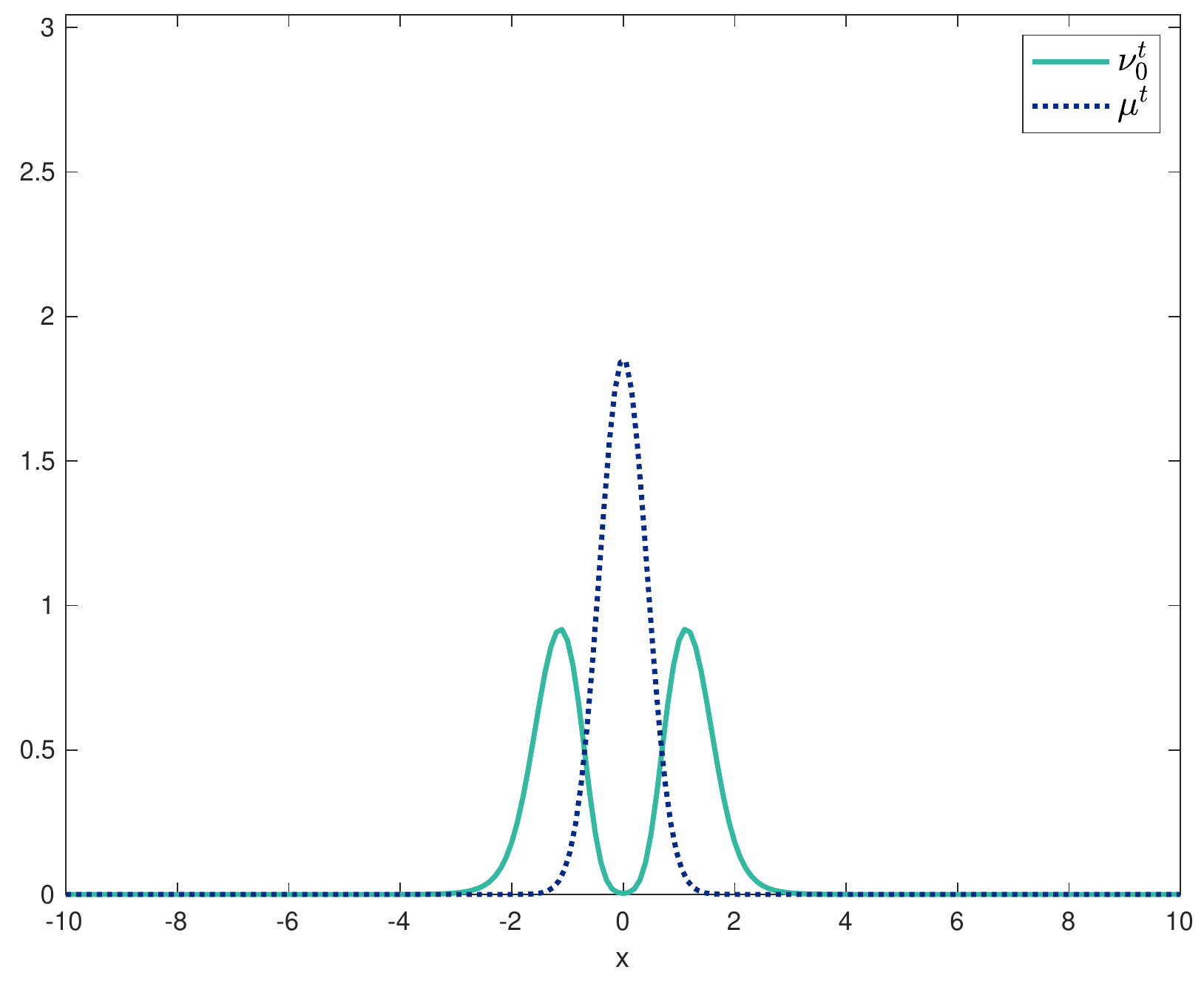}}
	\subfigure[]{\includegraphics[width=0.30\textwidth]{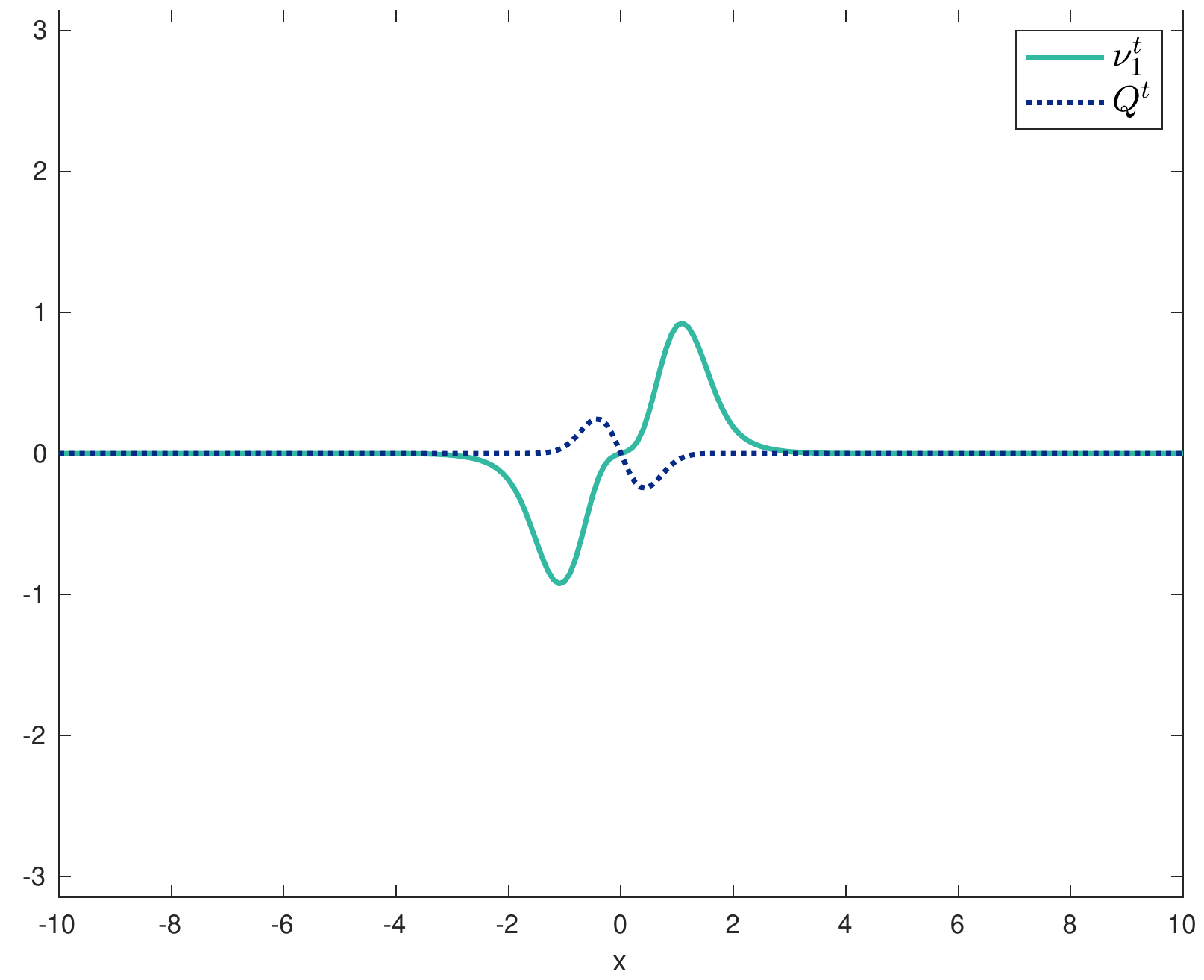}}\\
	\subfigure[]{\includegraphics[width=0.30\textwidth]{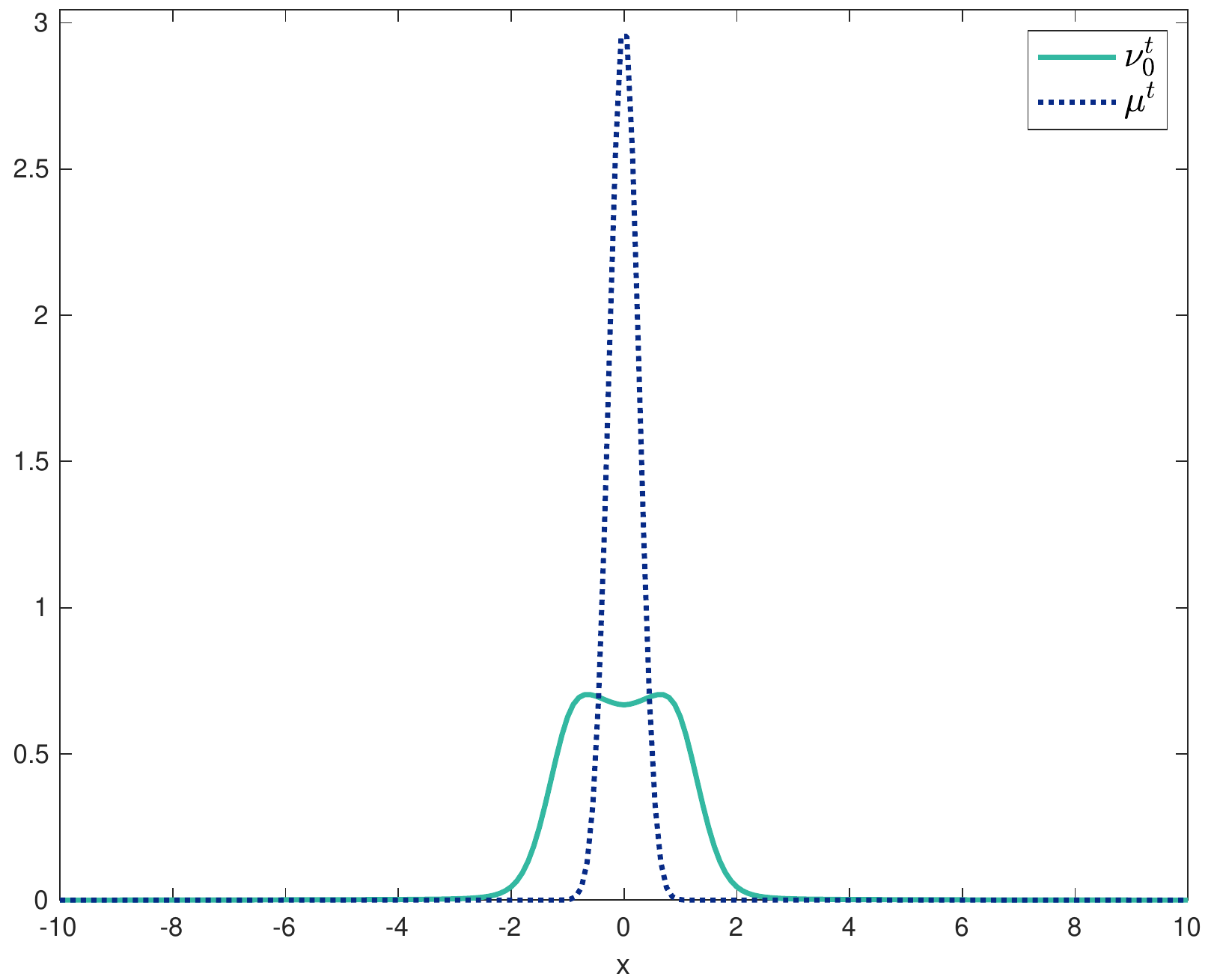}}
	\subfigure[]{\includegraphics[width=0.30\textwidth]{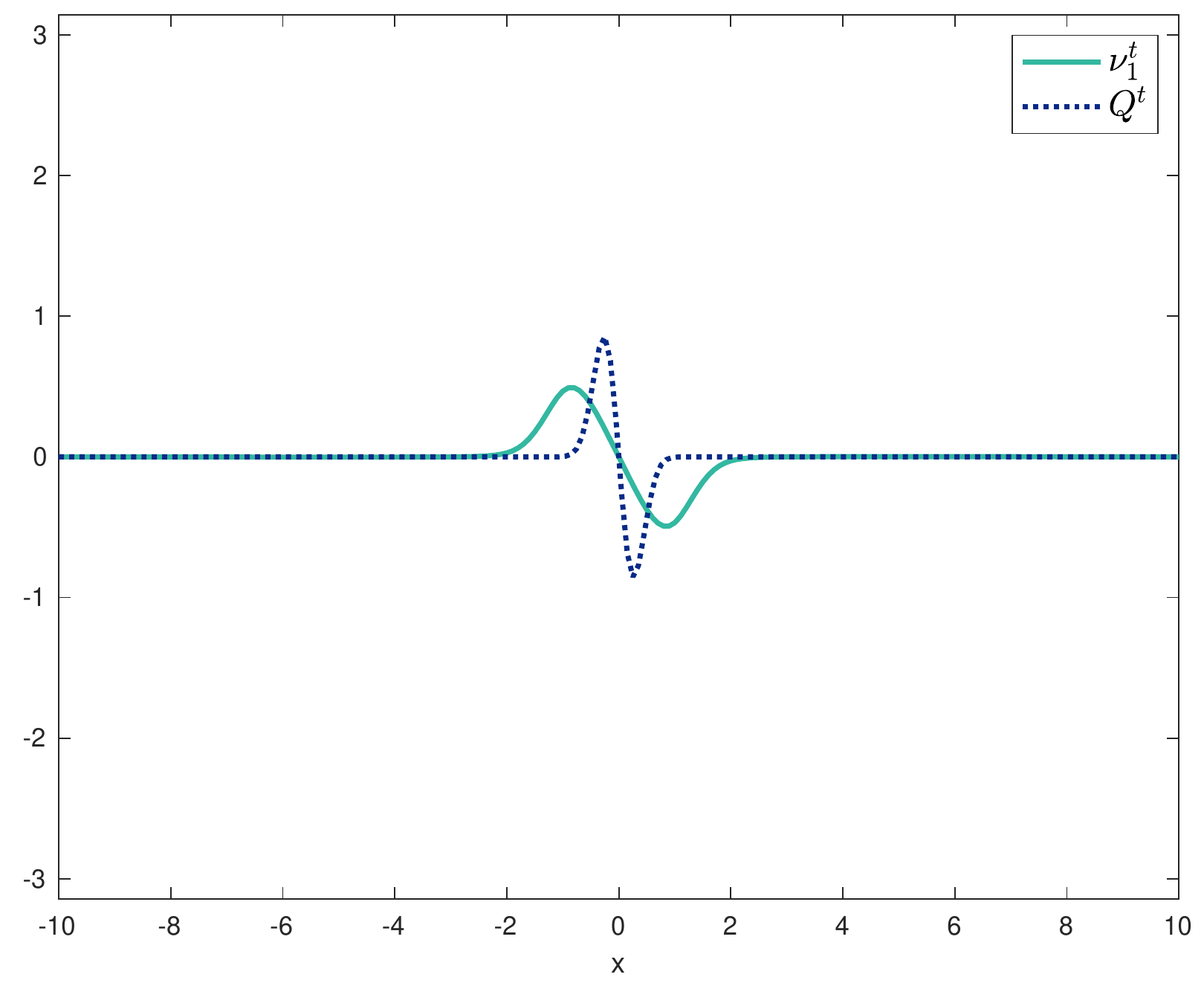}}\\
	\subfigure[]{\includegraphics[width=0.30\textwidth]{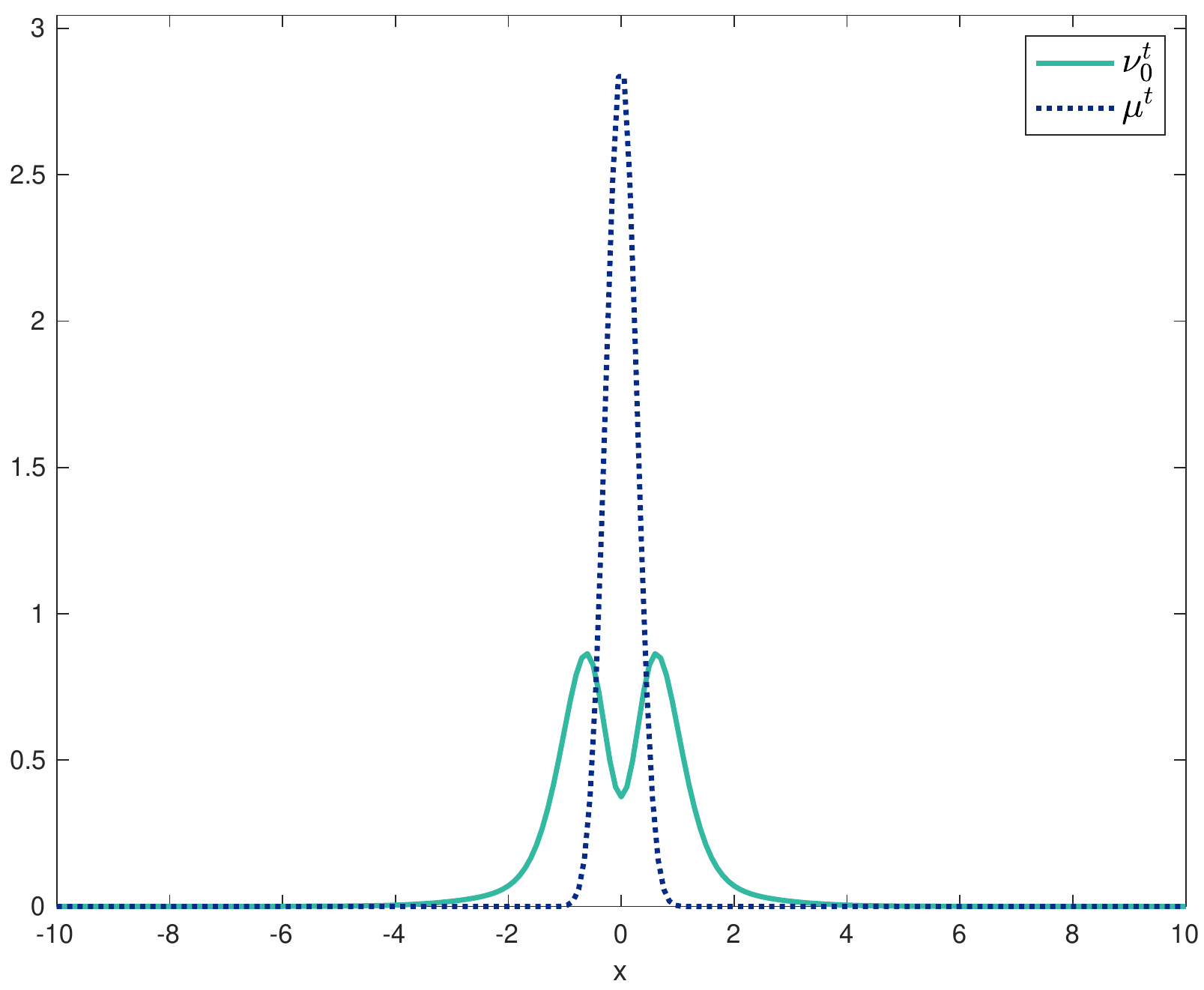}}
	\subfigure[]{\includegraphics[width=0.30\textwidth]{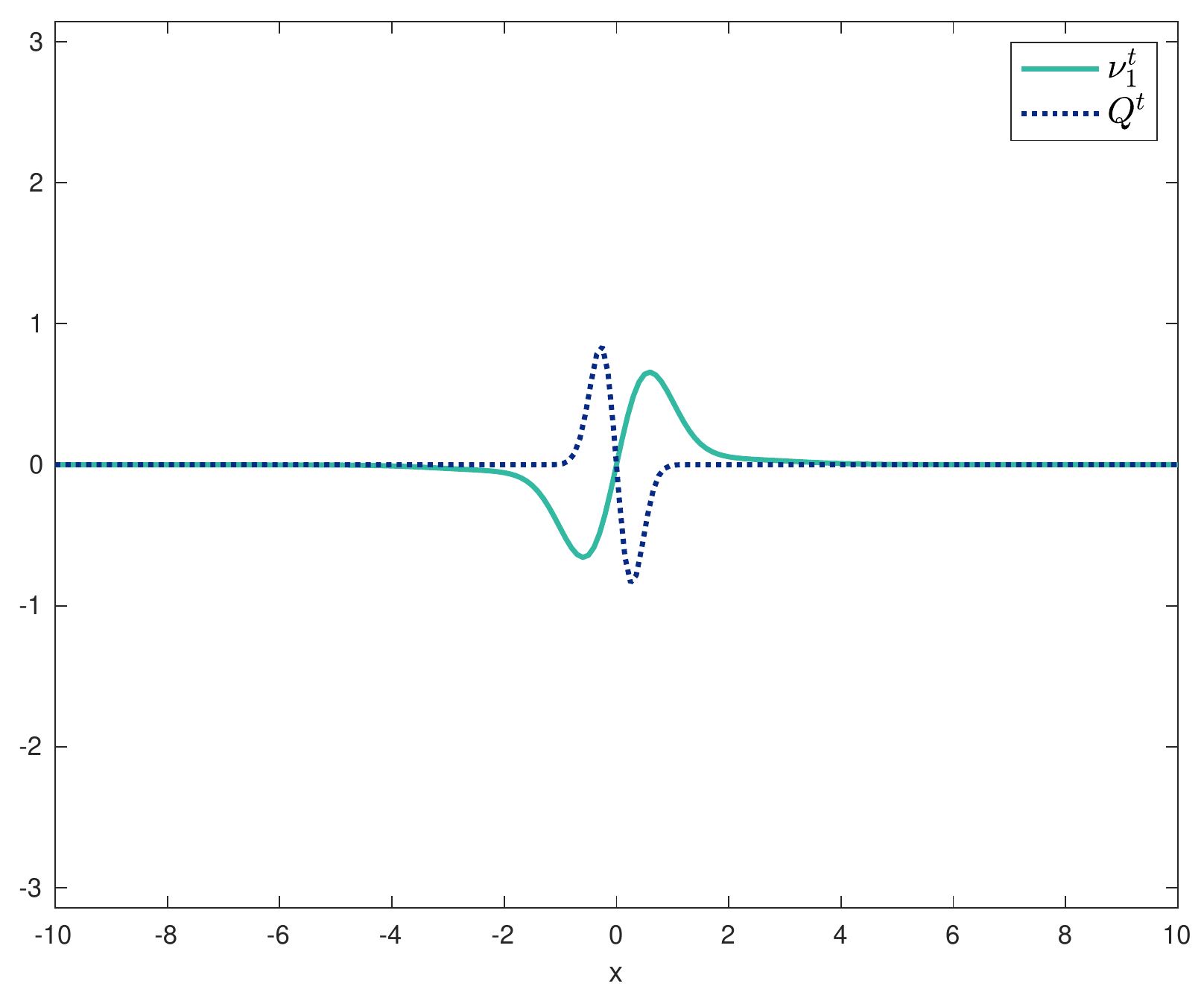}}
	\caption{Test 10:  Comparison between zero-th order moment of the solution of (V) and $\mu^t$,  first order moment of the solution of (V) and $Q^t$     a)-b) $t=0.5$, c)-d) $t=2$, e)-f) $t=3$.  Here $\eta=3.4$, $\alpha=0$ and the integral term in (E) is negletted.}
	\label{fig:VE_onlychemo}
\end{figure}

For the sake of completeness we consider the scenario without alignment, where only a chemotactic effect rules the dynamics. At kinetic level, we have already discussed the results (see Test 4 and Figure \ref{fig:Vlasov_onlychemo}). 
Figure \ref{fig:VE_onlychemo} shows the results of the comparison. We observe that the agreement between Vlasov moments and solutions of Euler is not recovered. The zero-th order moment shows the observed oscillations, with two bumps separating and merging, as reflecting by the sign of the first-order moment. The oscillating behavior is not reproduced in the solution of Euler. 
Regardless of the comparison between the two scale, numerical simulations in largely non monokinetic cases shed light on an important result.
As already stated, the main advantage of (E) is to clearly preserve the particle interactions in the passage from microscopic to the macroscopic scale, which are included in the nonlocal integral term on the right-hand side.
Neglecting the nonlocal integral term and damping, the system reduces to a pressureless Euler system, that arises in the modeling of sticky particles and has been investigated at theoretical level in several works of the literature (see \cite{sticky1,sticky2,bouchut} for seminal papers).
One of the main features of pressureless gas system is the development of delta-shocks.
Numerical simulations in Test 9 -Test 10 actually exhibit a blow-up in the density, which can be seen as a drawback of the pressureless structure of (E).
In particular we conclude that in a general setting, the presence of the integral term and/or the chemotactic gradient is not enough dissipative to prevent a delta-shocks phenomena.
In the next section we investigate the role of damping term, proposing a mixed pressure-nonlocal Euler system incorporating the system (E).

\newpage

\newpage

\newpage

\newpage

\subsection{Numerical comparison between (V) and (E) with pressure }\label{VEpressure}

We observed that, in the general case, the presence of the damping term helps in recovering agreement between Vlasov moments and the macroscopic quantities, density and momentum, solving the nonlocal pressureless Euler system. 
In absence of damping, or if the $\gamma$ function is not enough dissipative, a discrepancy between the two scales can be expected.
\REVISION{The structure of $(E)$ is actually due to the monokinetic property of the solution of $(V)$, but in the general case, also the second order moment should be involved when taking the momenta of the Vlasov solution. And the second moment would in turn involve the third one and so on. Therefore the hierarchy of equations satisfied by the moments would not be closed.}
In the following we extend the Euler system $(E)$ with an additional term approximating to a certain extent the role of the second order moment.
In the literature of hyperbolic chemotaxis model, the model in \cite{GambaPercolation}, originally conceived to reproduce the early stages of the vasculogenesis process, shares a common structure with $(E)$.
The main difference is that the nonlocal integral term in $(E)$ is replaced by a phenomenological density dependent pressure function, which does not have a clear connection with the microscopic interactions among cells, but it plays a crucial role in stabilizing the system.
Here we mix the original structure of $(E)$ adding a pressure term inspired by the model in \cite{GambaPercolation}. 
Note that the nonlocal-pressure structure has been investigated in the context of hydrodynamic limit of agent-based dynamics (without chemotaxis) in a recent paper \cite{tadmor2022}, proposing rigorous results on the long-time behavior of alignment Euler systems for a general class of entropy pressure.
Our aim is to investigate, from a numerical perspective, the effect of the additional pressure term, and to gain some insights on possible agreement between the scales, far from the monokinetic assumption on initial data. 
The hydrodynamic alignment Euler system we here consider reads: 

\begin{eqnarray}\label{Euler1D_damping_press}
\begin{array}{cl}
\left\{
\begin{array}{l}
\partial_t\mu^t+\partial_x(u^t \mu^t) = 0\\ \\
\partial_t(\mu^tu^t)+
\partial_x\left(\mu^t(u^t)^{2}+\varepsilon P(\mu^t)\right) = 
\mu^t\int\gamma(\cdot-y,u^t(\cdot)-u^t(y))\mu^t(y)dy\\ \\
\ \ \ \ \ \ \ \ \ \ \ \ \ \ \ \ \ \ \ \ \ \ \ \ \ \ \ \ \ \ \ \ \ \ \ \ \ \ \ \ \ \ \  +\eta\mu^t\partial_{x}\psi^t
-\alpha \mu^t u^t\\ \\
\partial_s\psi^s = D\partial^2_{x}\psi\REVISION{^s}-\kappa\psi\REVISION{^s} +
\displaystyle \int_{x-R}^{x+R}\mu^s(y)dy,\ s\in[0,t].
\end{array}
\right.
\end{array}
\end{eqnarray}

We assume a pressure law for isentropic gases $P(\mu)=\mu^{\REVISION{p}}$,
choosing $p=2$ in our simulations. 
Parameter $\varepsilon \ge 0$ tunes the effect of the pressure term. Choosing $\varepsilon =0$, we recover the pressureless system \eqref{Euler1D_damping}. In the following we present the results of numerical simulations performed for different values of $\varepsilon$. As in the previous case, the effect of the damping term is investigated.

\begin{figure}[h!]
	\centering
	\subfigure[]{\includegraphics[width=0.30\textwidth]{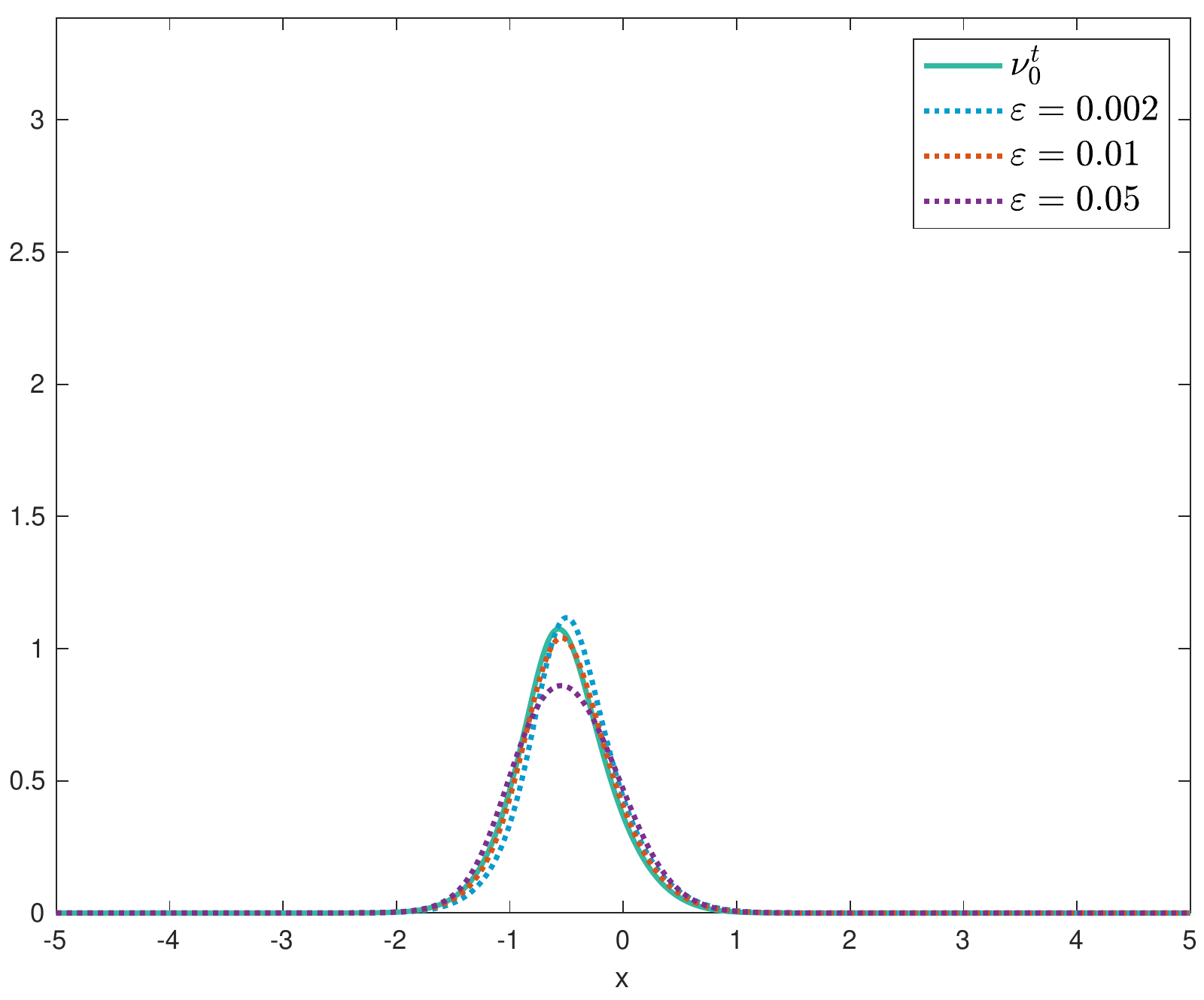}}
	\subfigure[]{\includegraphics[width=0.30\textwidth]{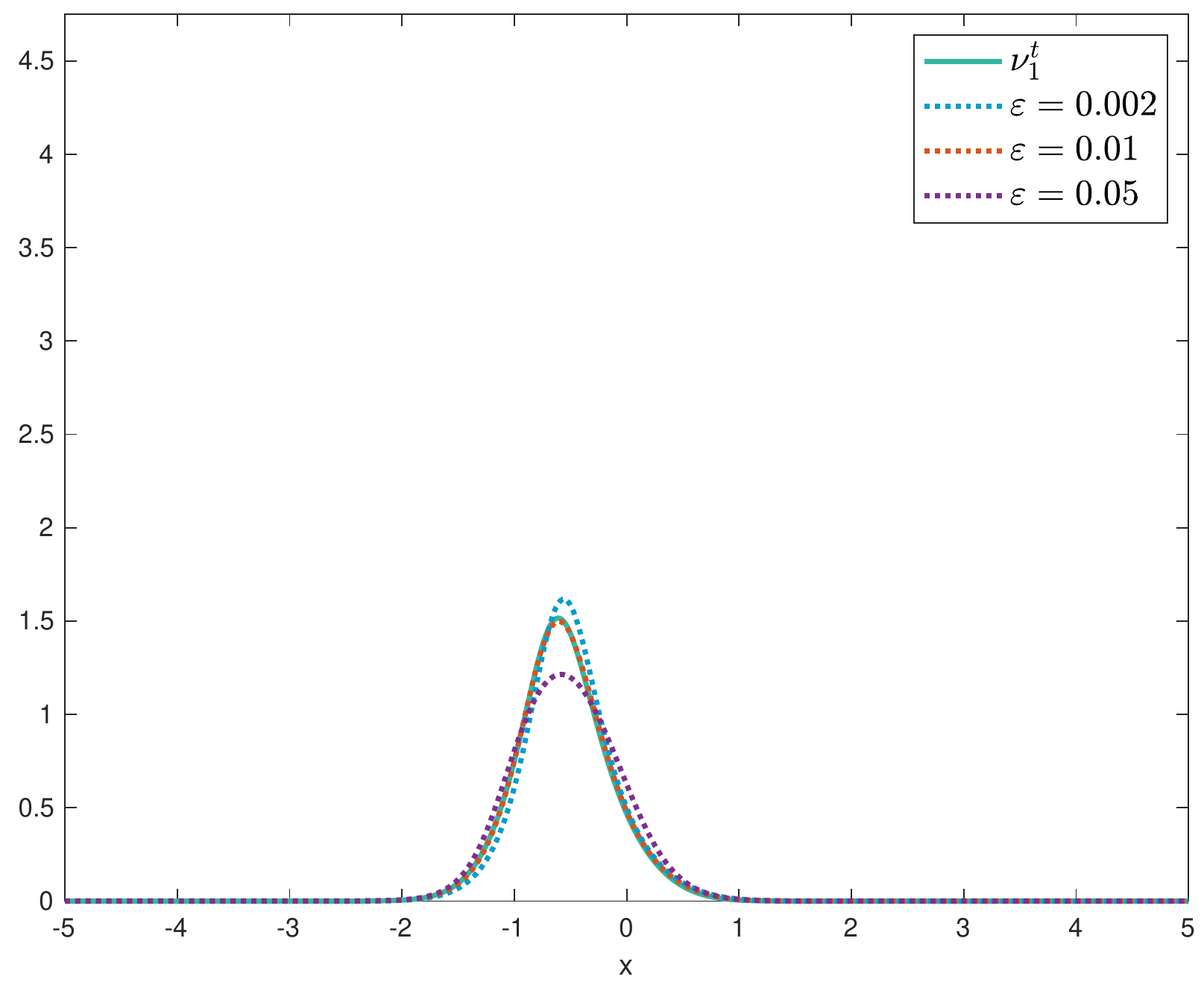}}\\
	\subfigure[]{\includegraphics[width=0.30\textwidth]{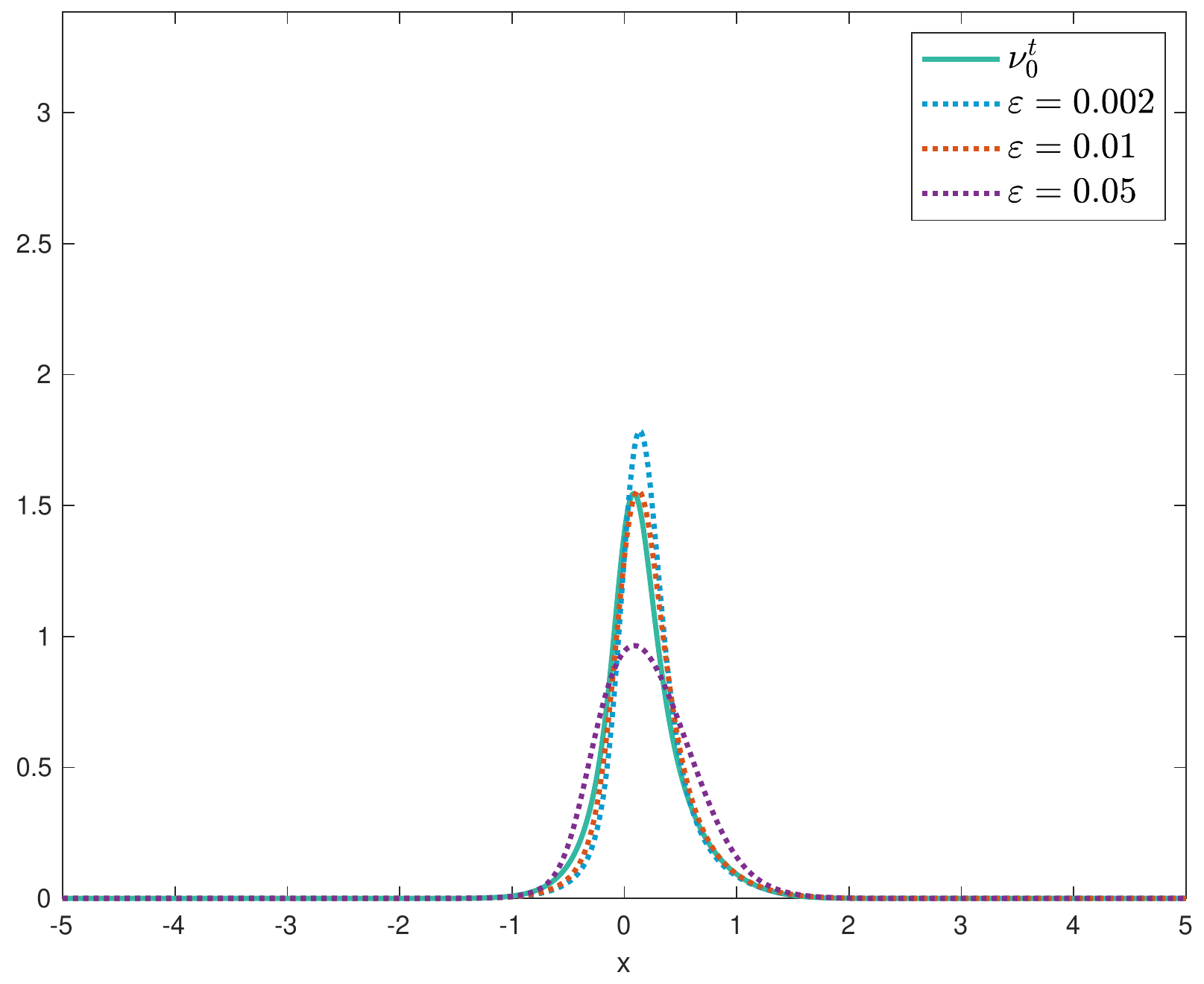}}
	\subfigure[]{\includegraphics[width=0.30\textwidth]{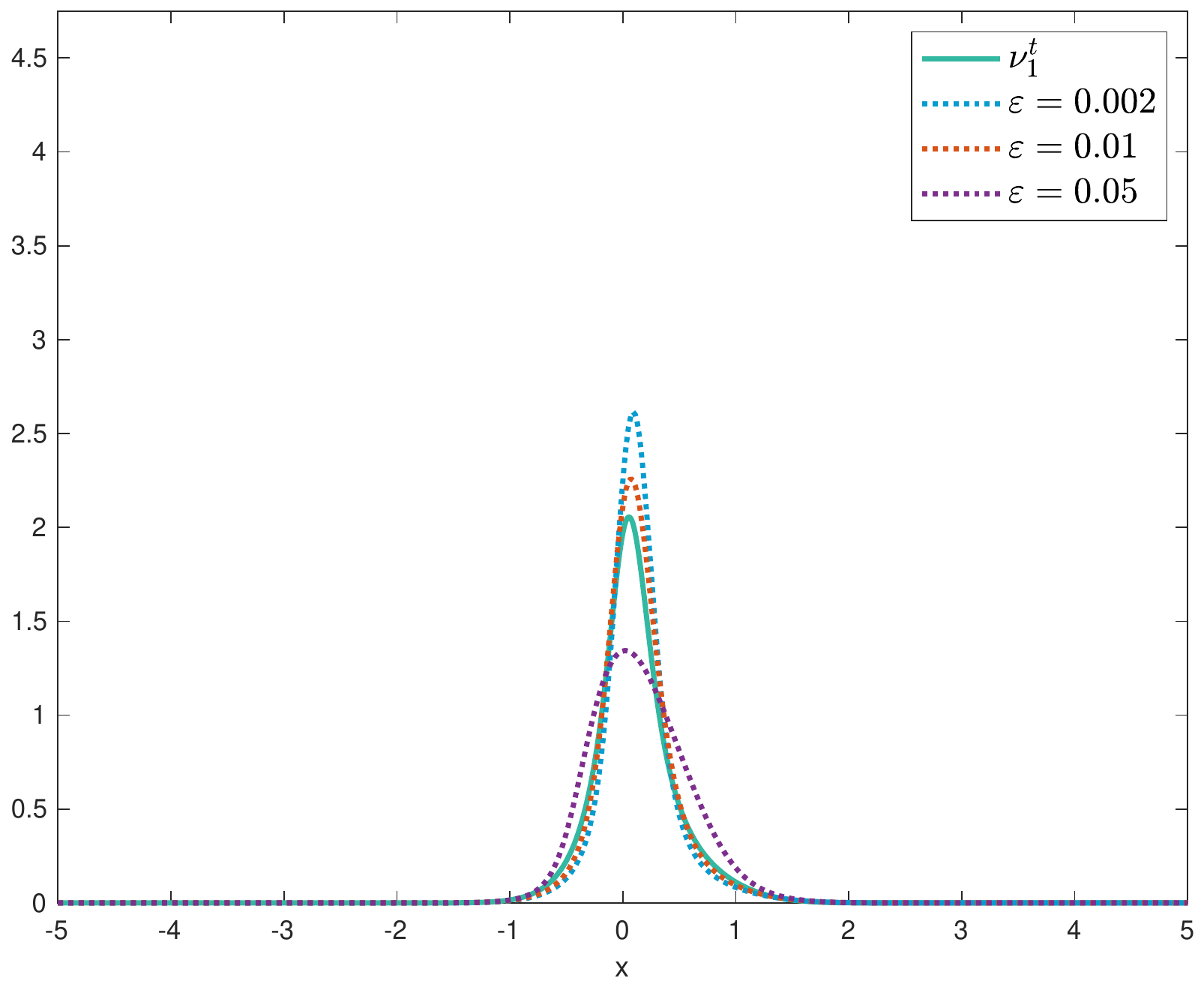}}\\
	\subfigure[]{\includegraphics[width=0.30\textwidth]{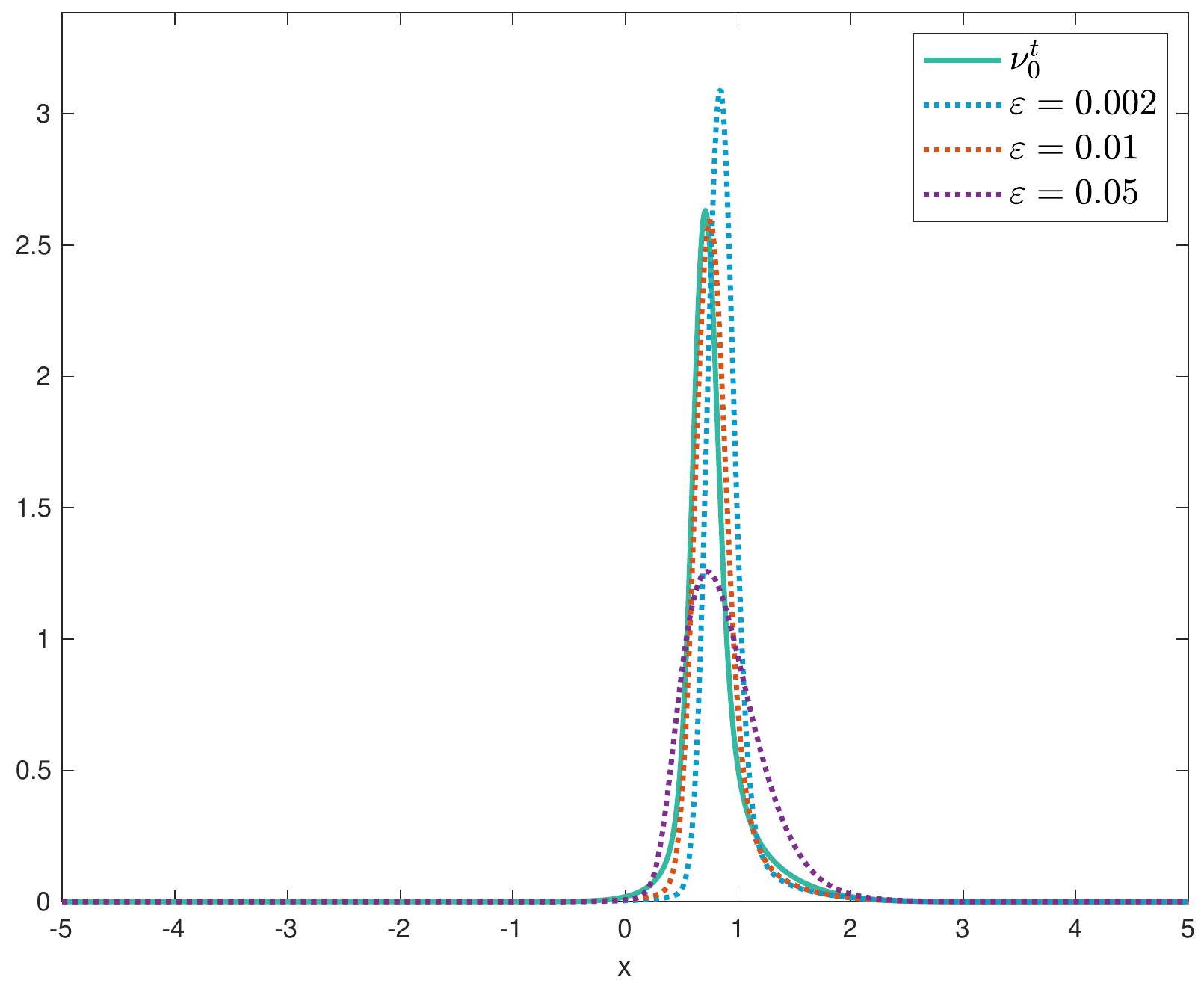}}
	\subfigure[]{\includegraphics[width=0.30\textwidth]{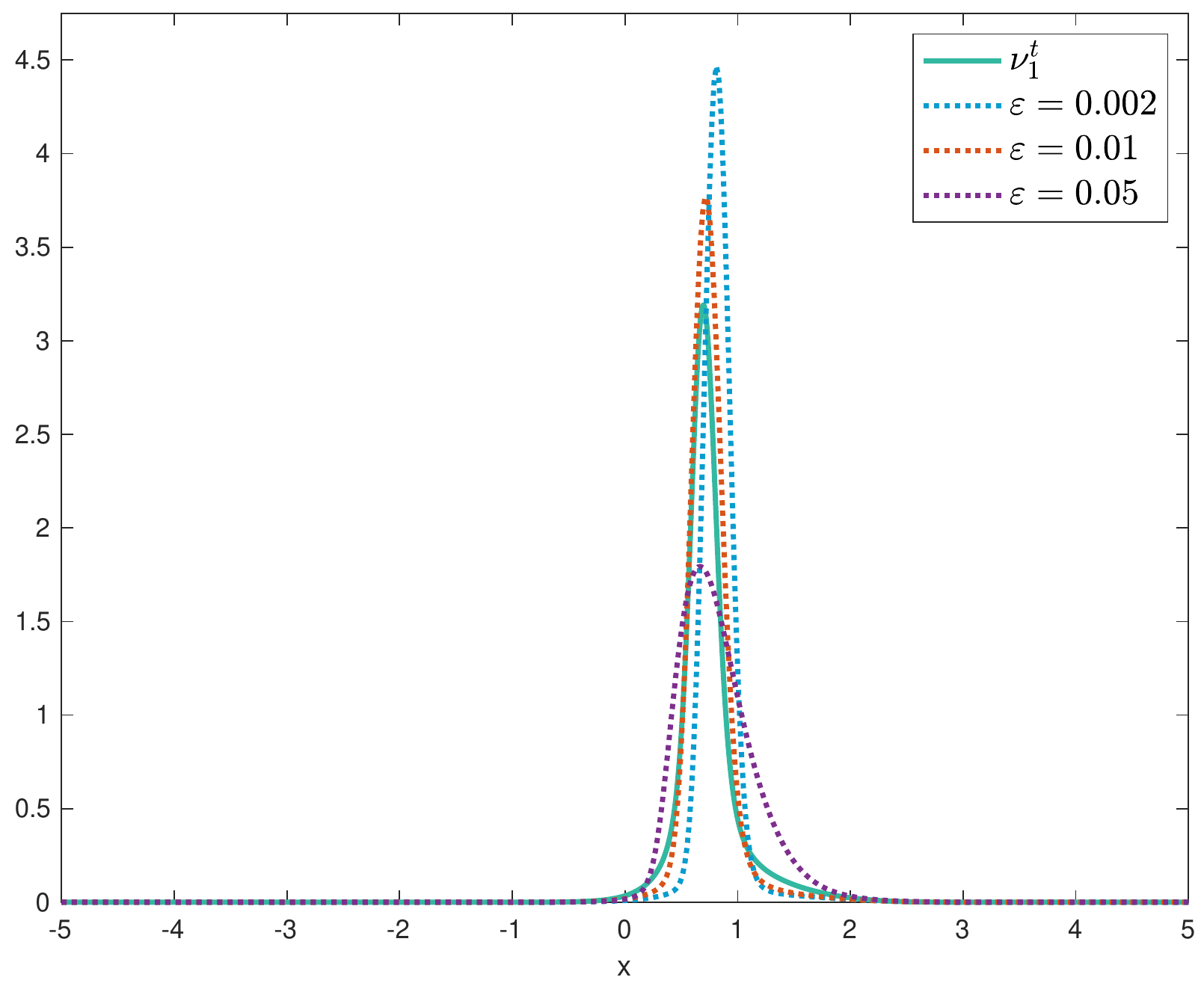}}
	\caption{Test 11: Screenshots of numerical simulation of \eqref{Euler1D_damping_press} for different values of $\varepsilon$. Comparison between zero-th order moment of the solution of (V) and $\mu^t$,  first order moment of the solution of (V) and $Q^t$  at a)-b) $t=1$, c)-d) $t=1.5$, e)-f) $t=2$.  Initial condition and parameters are the same of Test 7. }
	\label{fig:monokinetic_press}
\end{figure}

As a preliminary test, we run the system starting with the (approximated) monokinetic initial data in \eqref{rho0_monokin} following the same approach already considered for Test 7.
Figure \ref{fig:monokinetic_press} shows the comparison between the first two moments of the solution of \eqref{VlasovChemo} and density-momentum solutions of \eqref{Euler1D_damping_press}
 at three different time-steps, for different values of $\varepsilon \ge 0$.

Varying the value of $\varepsilon$ we get different solutions of system \eqref{Euler1D_damping_press}.
As expected, the presence of a pressure term introduces a dissipative behavior, which increases with the value of $\varepsilon$. Tuning the value of the parameter starting with $\varepsilon=0$, we observe that, at least from a qualitative perspective, there exists a critical value $\varepsilon^*>0$ such that
 $\varepsilon < \varepsilon^*$ the pressure term is not enough to counterbalance the density blow up, whereas if $\varepsilon >\varepsilon^*$ the effect of the pressure term is predominant and the dynamic is completely different from the kinetic one. 
 The results exhibited in Figure \ref{fig:monokinetic_press} show that $\epsilon^*=0.01$ represents a good candidate as optimal value minimizing the distance between $\nu^t_0$ and $\mu^t$.
 
In order to quantitatively estimate the agreement between $\nu^t_0$ and $\mu^t$, and $\nu^t_1$ and $Q^t$, we compute, at each iteration $k=1,2,...$  the value of distances $E_0^{k}(\varepsilon)$, $E_1^{k}(\varepsilon)$ defined by
\begin{equation}
E_0^k (\varepsilon):= \| \nu_0^k-\mu^k(\varepsilon)   \|_{L^1(\Omega)} \ \ \ \ \ \
E_1^k(\varepsilon) := \| \nu_1^k-Q^k(\varepsilon)   \|_{L^1(\Omega)}
\end{equation}
where  $ \mu^k(\varepsilon), Q^k(\varepsilon)$ denote the solution of \eqref{Euler1D_damping_press} obtained for a fixed value of $\varepsilon$ \REVISION{at iteration $k$}.

Table \ref{tab:norm_E0} exhibits the obtained values for the three time iterations considered in Figure  \ref{fig:monokinetic_press}. 
Computing the quantities at any iteration $k$, we observe that the values of $E_0^k$ and $E_1^k$ decrease as $\varepsilon$ increases from 0 to a certain value $\varepsilon^*$, whereas their value increase for $\varepsilon > \varepsilon^*$, giving insights for the existence of a minimum point.
To this end, we run an optimization procedure, aiming at finding, at each iteration, the value $\epsilon^*_k \in [0, 0.05] $  minimizing $E_0^k$. The same can be done for $E_1^k$ (data not shown).
Table \ref{tab:optimazed_epsilon} shows the results obtained for the selected time instants. Considering an accuracy of three decimal places, we observe that the optimal values are slightly higher than our qualitative guess.

\begin{table}[h!]
		\begin{center}
	\begin{tabular}{|l|c|c|c|c|c|c|c|c|c|}
		\hline
		$E_0^k$ &
		$\varepsilon_1=0.002$&
		$\varepsilon_2=0.01$ &
		$\varepsilon_3=0.05$ 
		\\ \noalign{\hrule height2pt}
		$t_k=0$  & 0  & 0  & 0    \\ \hline
  		$t_k=1$  & 0.09  & 0.08   & 0.16      \\ \hline
                 $t_k=1.5$  & 0.23  & 0.18   & 0.32      \\ \hline
        		$t_k=2$  & 0.75  & 0.67   & 0.74     \\ 
		\hline
	\end{tabular}
	
	\vspace{1cm}
	
	\begin{tabular}{|l|c|c|c|c|c|c|c|c|c|}
		\hline
		$E_1^k$ &
		$\varepsilon_1=0.002$&
		$\varepsilon_2=0.01$ &
		$\varepsilon_3=0.05$ 
		\\ \noalign{\hrule height2pt}
		$t_k=0$  & 0  & 0  & 0    \\ \hline
  		$t_k=1$  & 0.06  & 0.05   & 0.23      \\ \hline
                 $t_k=1.5$  & 0.32  & 0.20   & 0.43      \\ \hline
        		$t_k=2$  & 0.99  & 0.85   & 0.95     \\ 
		\hline
	\end{tabular}

	\end{center}
\caption{Value of $E_0^k$ and $E_1^k$ for different values of $\varepsilon$.
As representative value we choose  $\varepsilon_1=\varepsilon_2/5$ and $\varepsilon_3=\varepsilon_2 * 5$, where $\varepsilon_2$ represents the numerical (qualitative) candidate for the optimal value of parameter $\varepsilon$. Time instants correspond to the ones selected for snapshots in Figure \ref{fig:monokinetic_press}.  }
\label{tab:norm_E0}
\end{table}

\begin{table}[h!]
		\begin{center}
	\begin{tabular}{|l|c|c|c|c|c|c|c|c|c|}
		\hline
		 &
		$\varepsilon^*_k$&
		$E_0^k(\varepsilon^*_k)$ 
		\\ \noalign{\hrule height2pt}
		$t_k=0$  & 0  & 0     \\ \hline 
  		$t_k=1$  & 0.013  & 0.07       \\ \hline 
                 $t_k=1.5$  & 0.015  & 0.14        \\ \hline 
        		$t_k=2$  & 0.017  & 0.55        \\ 
		\hline
	\end{tabular}
		\end{center}
\caption{Computation of $\epsilon^*_k$ and corresponding value of $E^k_0$. Time instants correspond to the ones selected for snapshots in Figure \ref{fig:monokinetic_press}.   }
\label{tab:optimazed_epsilon}
\end{table}

\newpage

We recall that the correspondence between Vlasov moments and Euler solutions, in the monokinetic case, is analytically proved in \cite{rt2} for the pressureless case, and there is no a priori reason to expect some sort of correspondence adding the density-dependent pressure term.
The results shown in Figure \ref{fig:monokinetic_press} suggest that the moments of the solution of Vlasov system \eqref{VlasovChemo}, starting with something that is like-monokinetic, look like solutions of Euler with pressure, at least for some values of $\varepsilon$.

\begin{figure}[h!]
	\centering
	\subfigure[]{\includegraphics[width=0.30\textwidth]{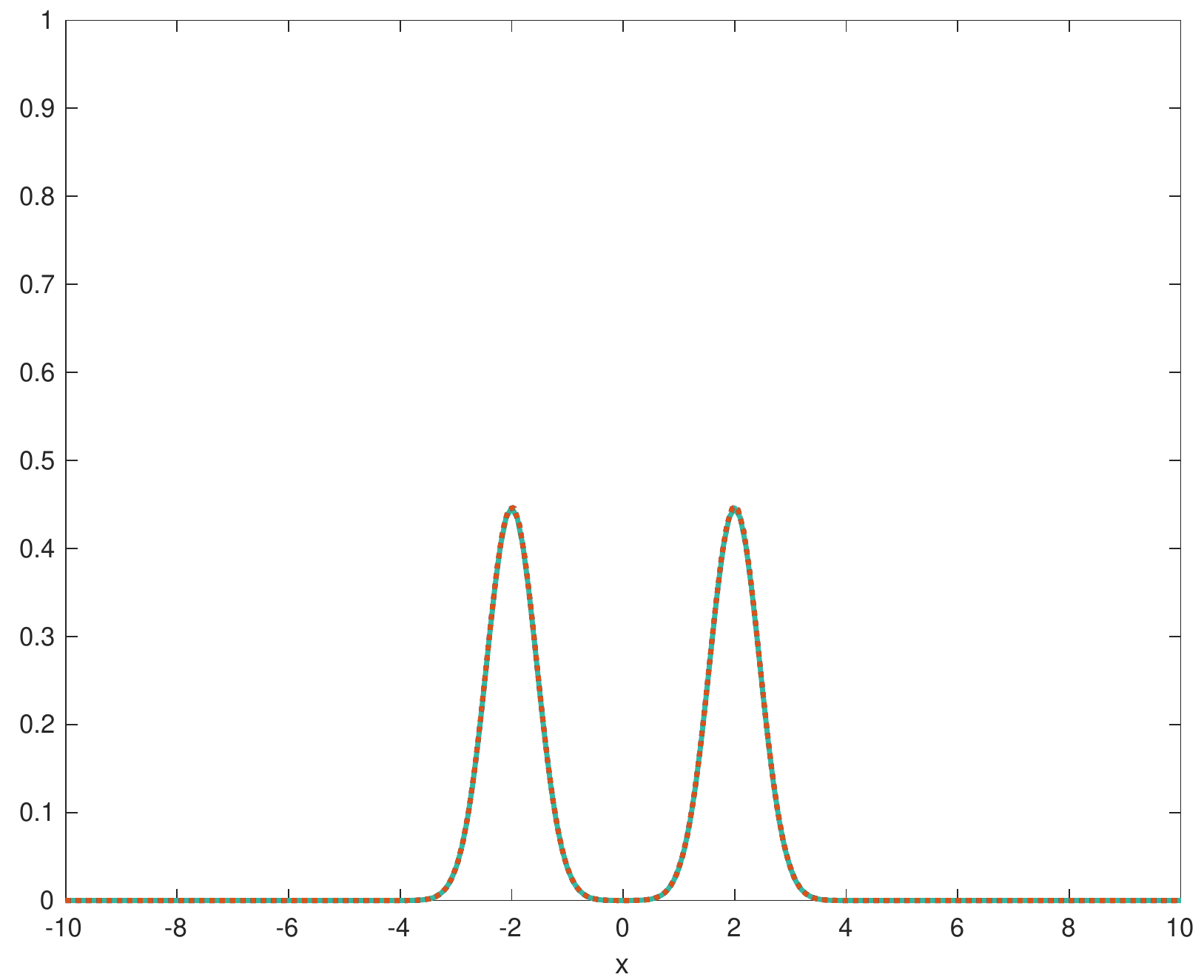}}
	\subfigure[]{\includegraphics[width=0.30\textwidth]{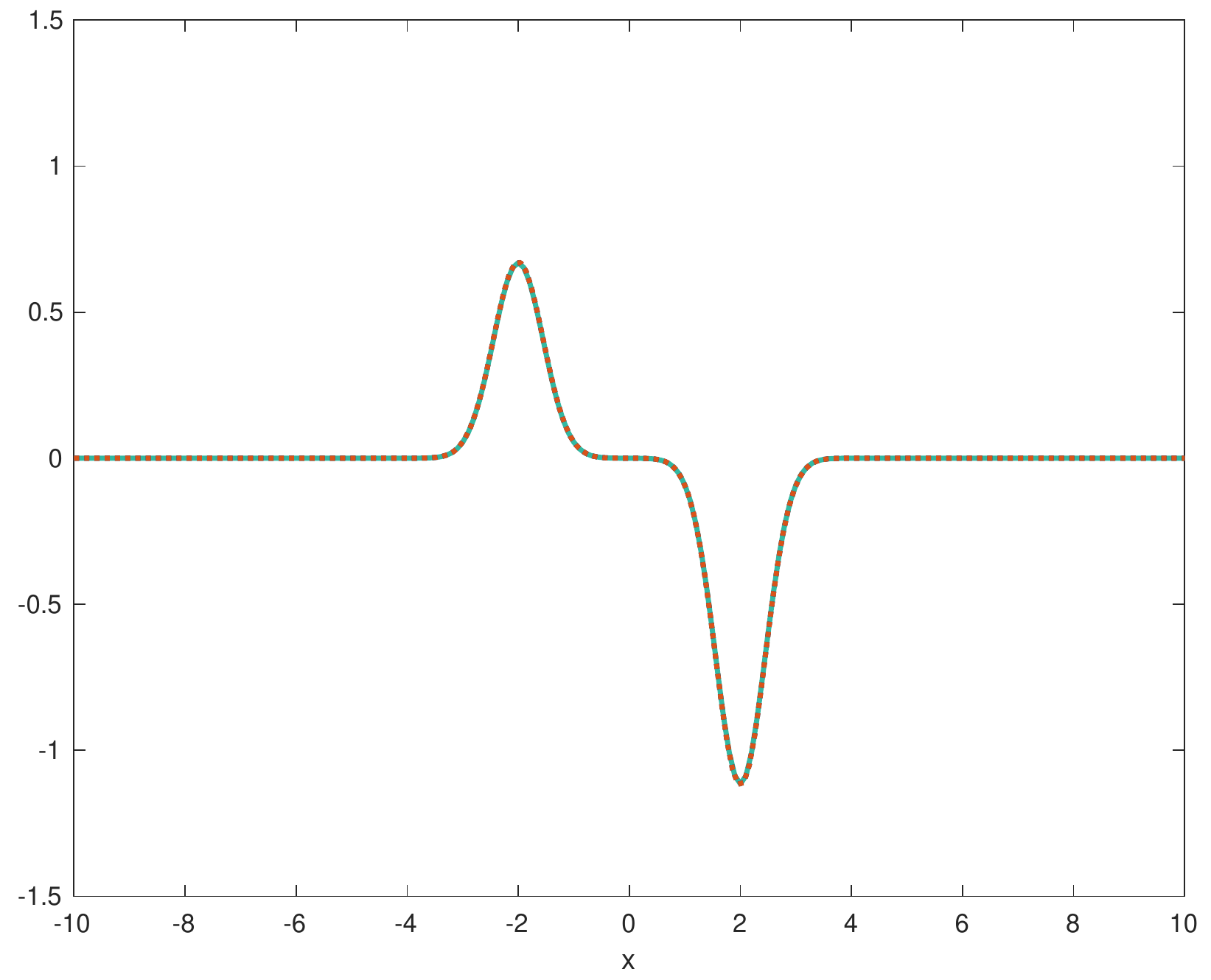}}
	\caption{Test 12, Test 13: a) initial density $\mu^0$ and b) initial momentum $Q^0$.  }
	\label{fig:2bumps_initialcondition}
\end{figure}

 In a second test we consider a general largely non monokinetic scenario.
Initial data are shown in Figure  \ref{fig:2bumps_initialcondition} . 
As defined in \eqref{mom_data}, $\mu^0$ and $Q^0$ are choosen as zero and first moment of the following initial distribution

\begin{equation}\label{rho0_2bumps}
\rho^0(x,v)= \frac{1}{\sqrt{2 \pi \sigma_x  \sigma_v}} 
\left(e^{\frac{-(x^2-x_1)}{2\sigma_x^2}}
e^{\frac{-(v-v_1)^2}{2\sigma_v^2}} + 
e^{\frac{-(x^2-x_2)}{2\sigma_x^2}}
e^{\frac{-(v-v_2)^2}{2\sigma_v^2}} \right)
\end{equation}
with $x_1=-2$, $v_1=1.5$, $x_2=2$, $v_2=-2.5$, $\sigma_x=\sqrt{0.2}$, $\sigma_v=\sqrt{0.5}$.

\begin{figure}[h!]
	\centering
		\subfigure[]{\includegraphics[width=0.30\textwidth]{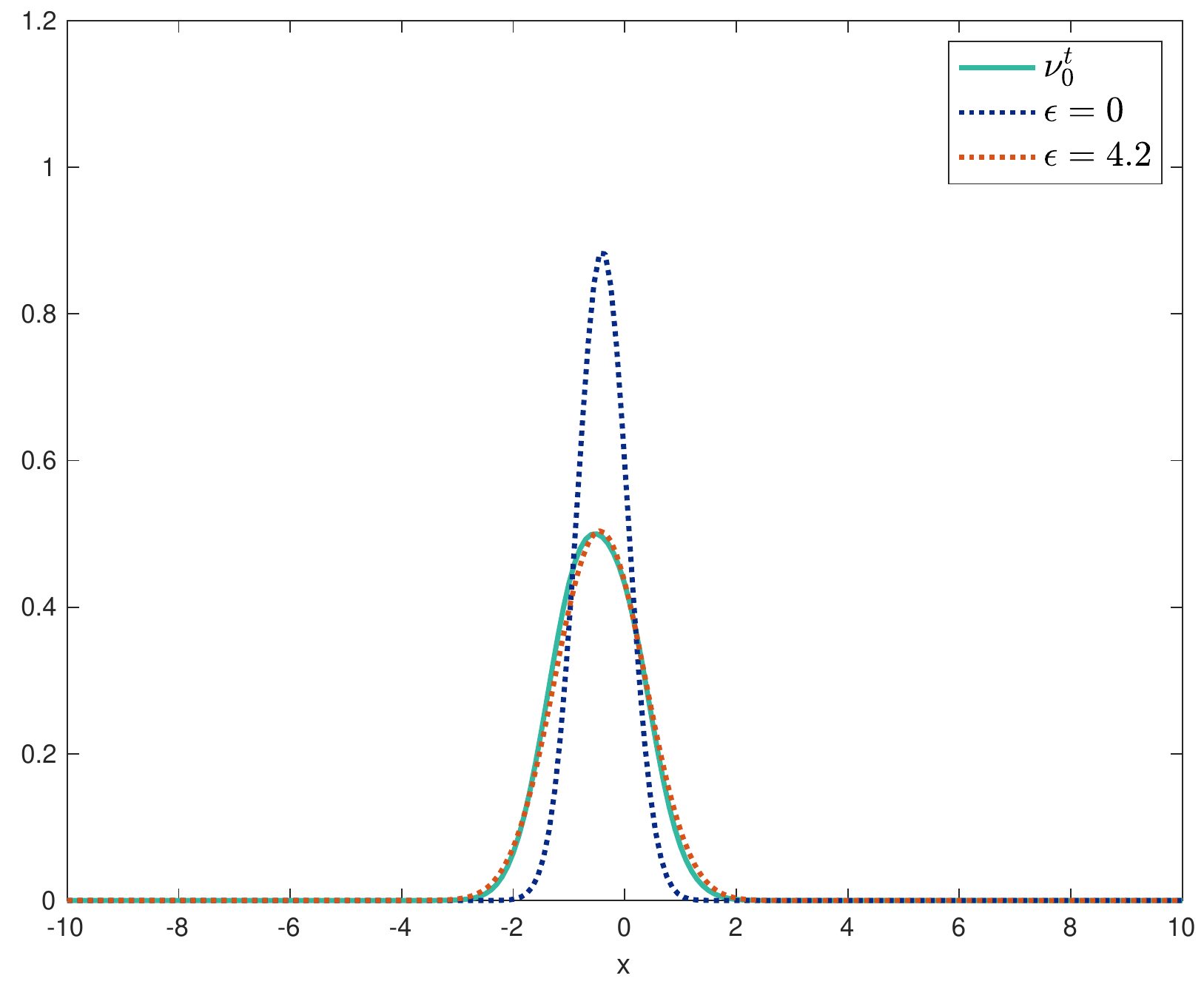}}
	\subfigure[]{\includegraphics[width=0.30\textwidth]{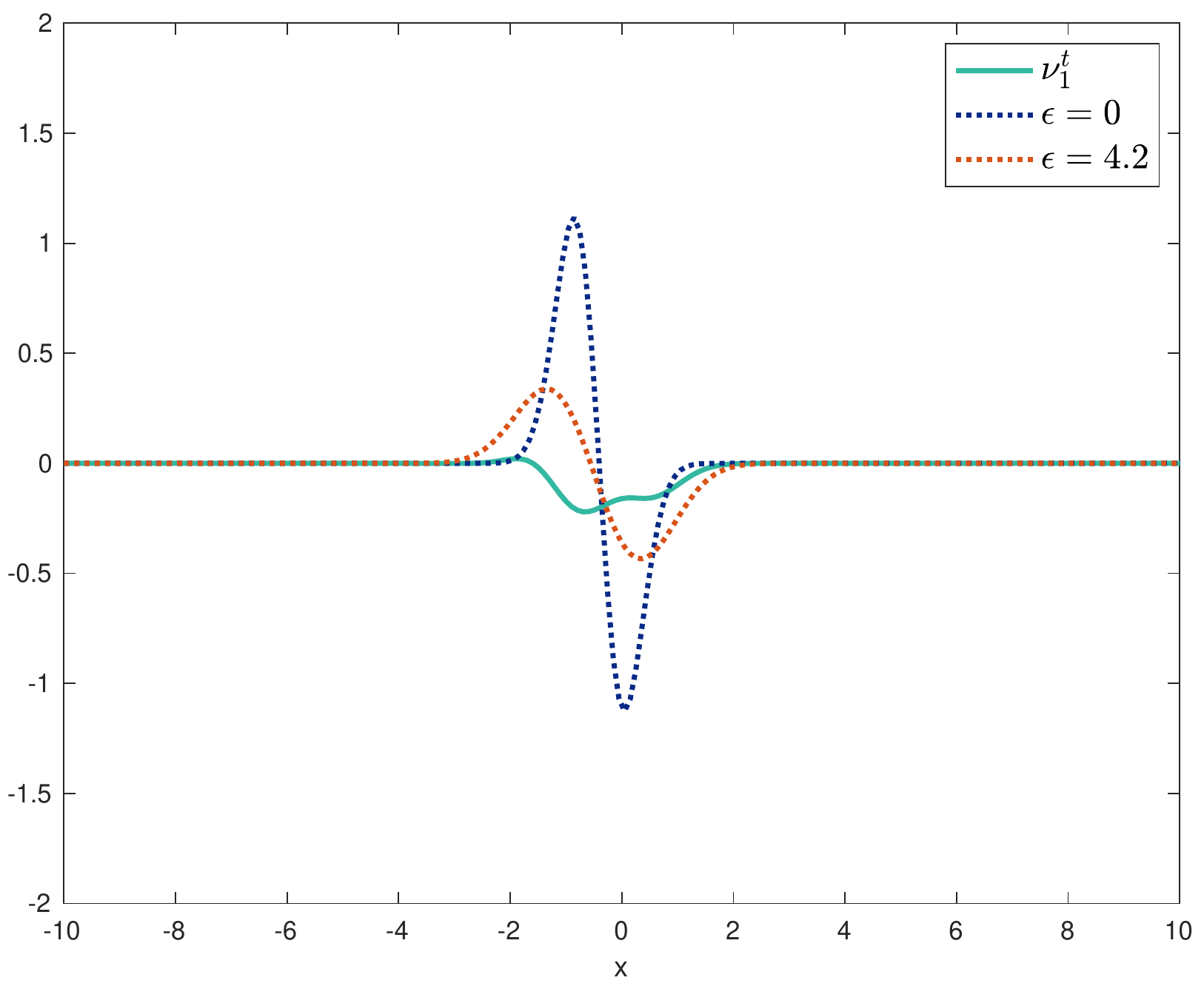}}\\
	\subfigure[]{\includegraphics[width=0.30\textwidth]{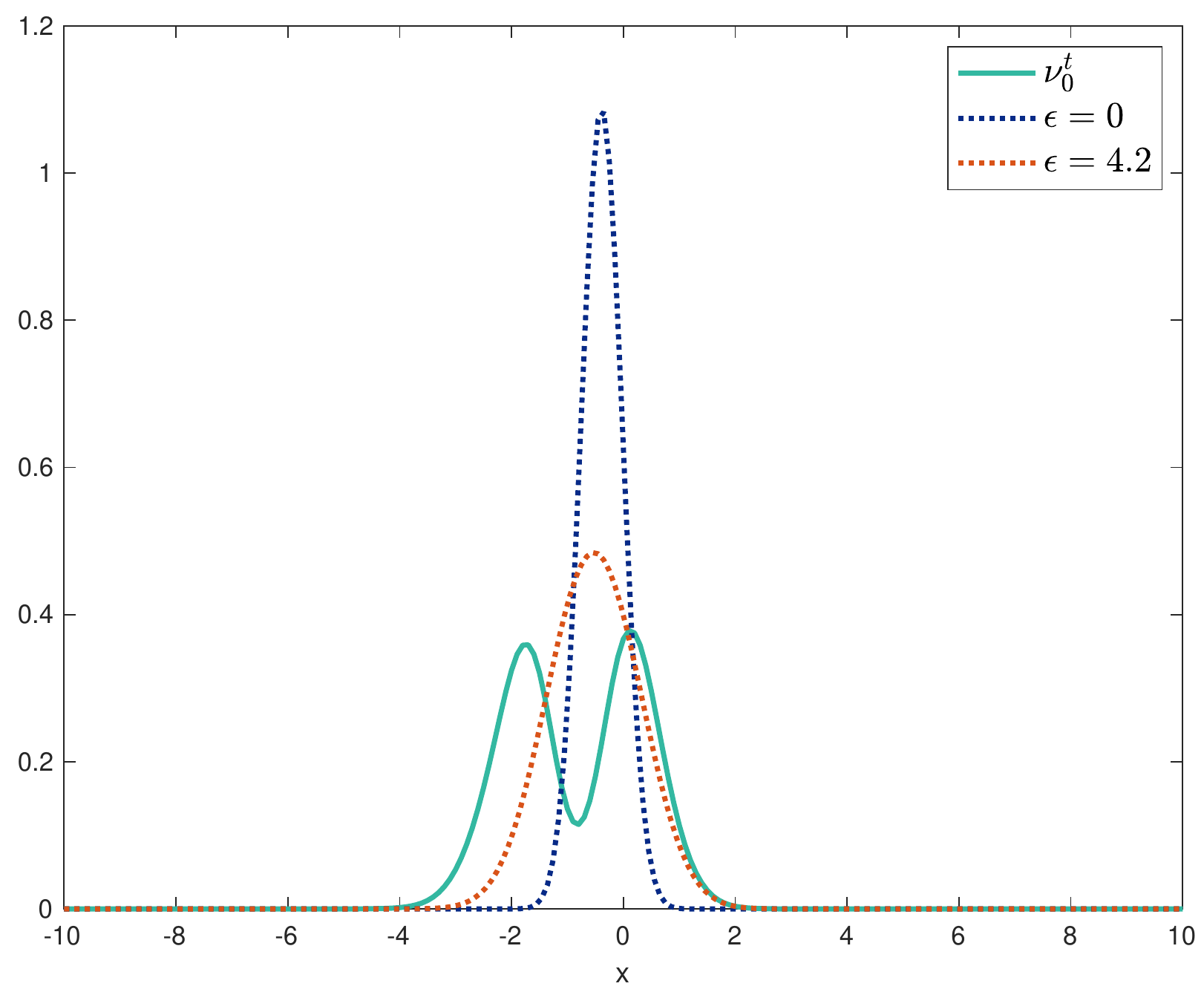}}
	\subfigure[]{\includegraphics[width=0.30\textwidth]{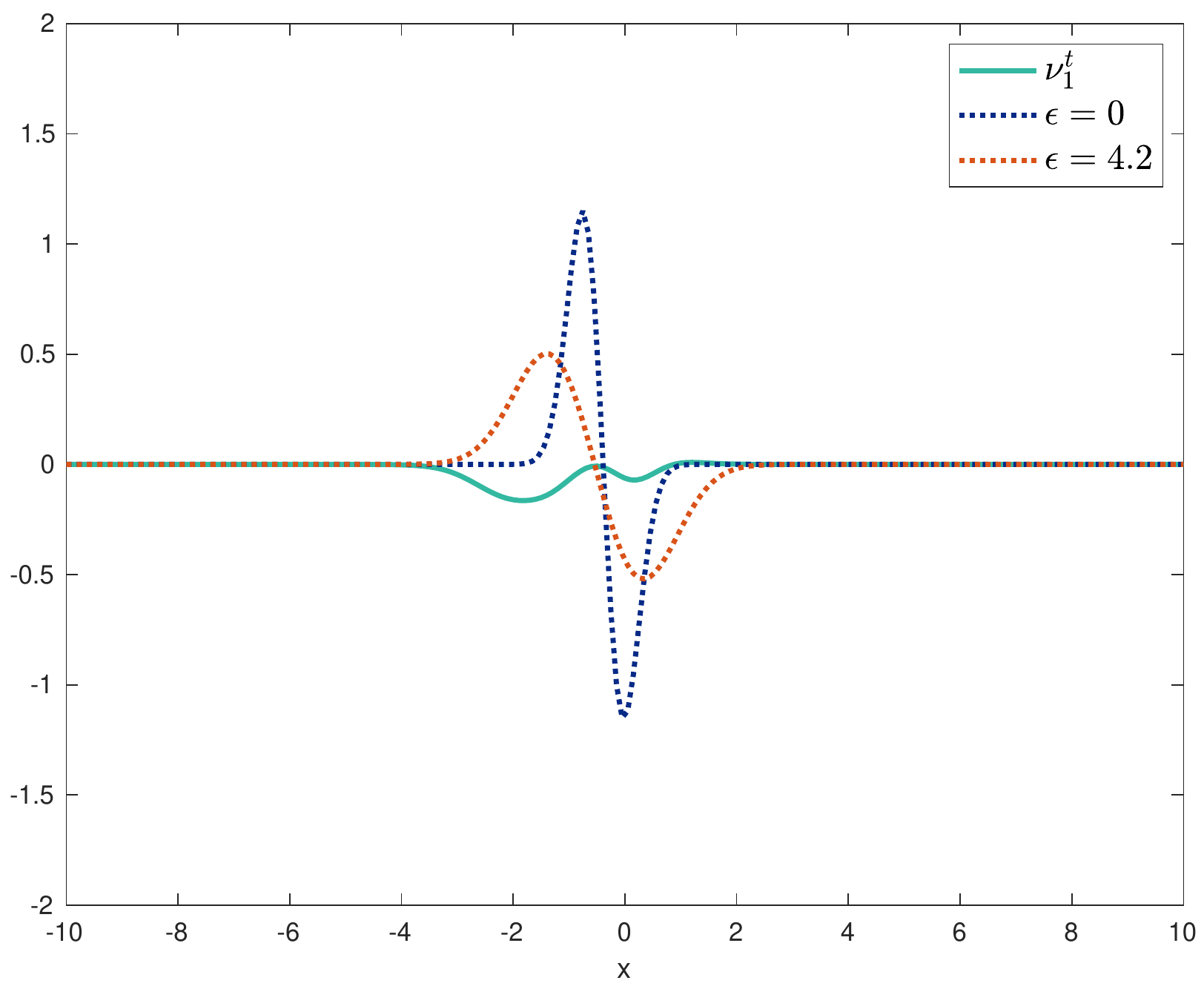}}\\
	\subfigure[]{\includegraphics[width=0.30\textwidth]{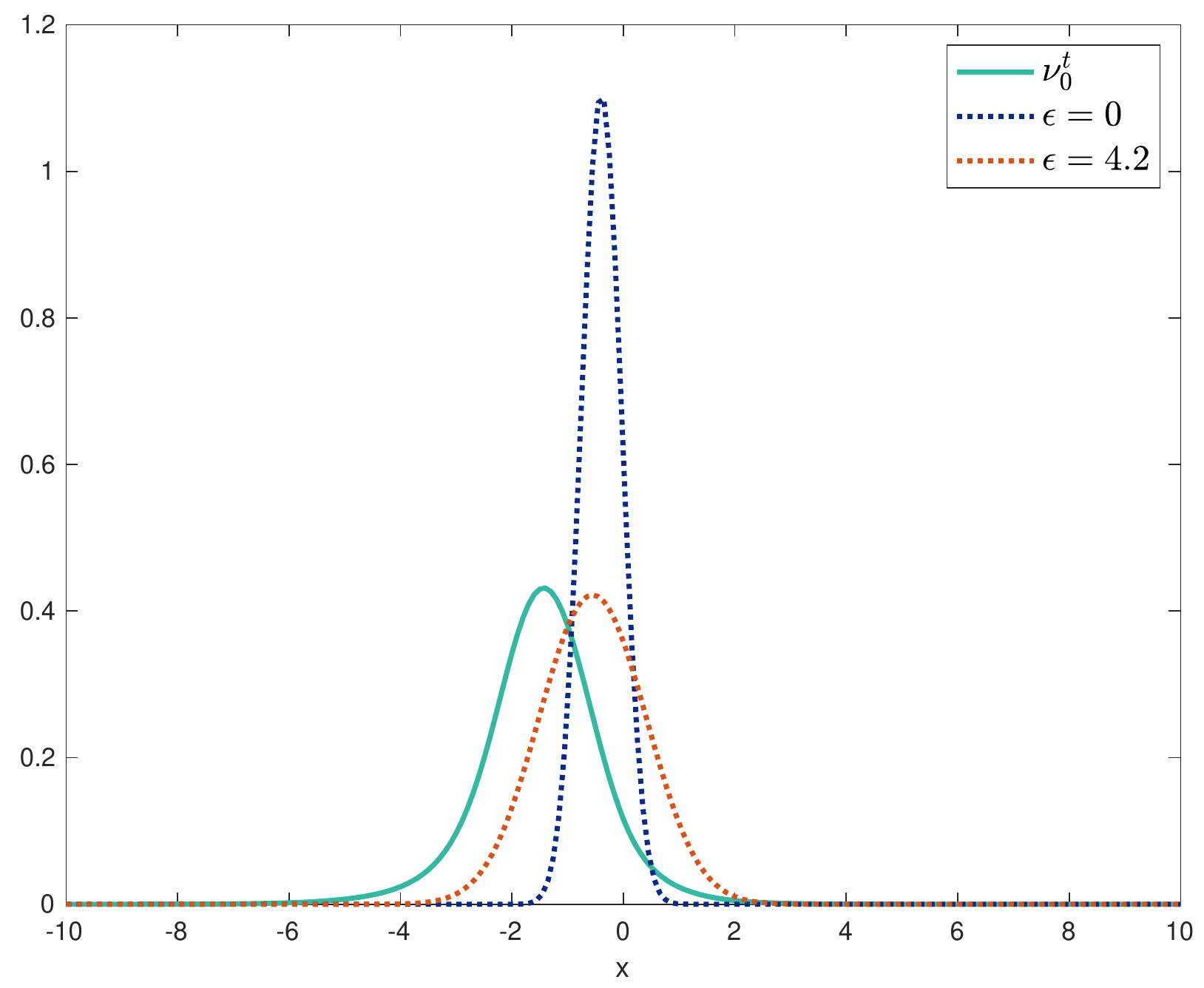}}
	\subfigure[]{\includegraphics[width=0.30\textwidth]{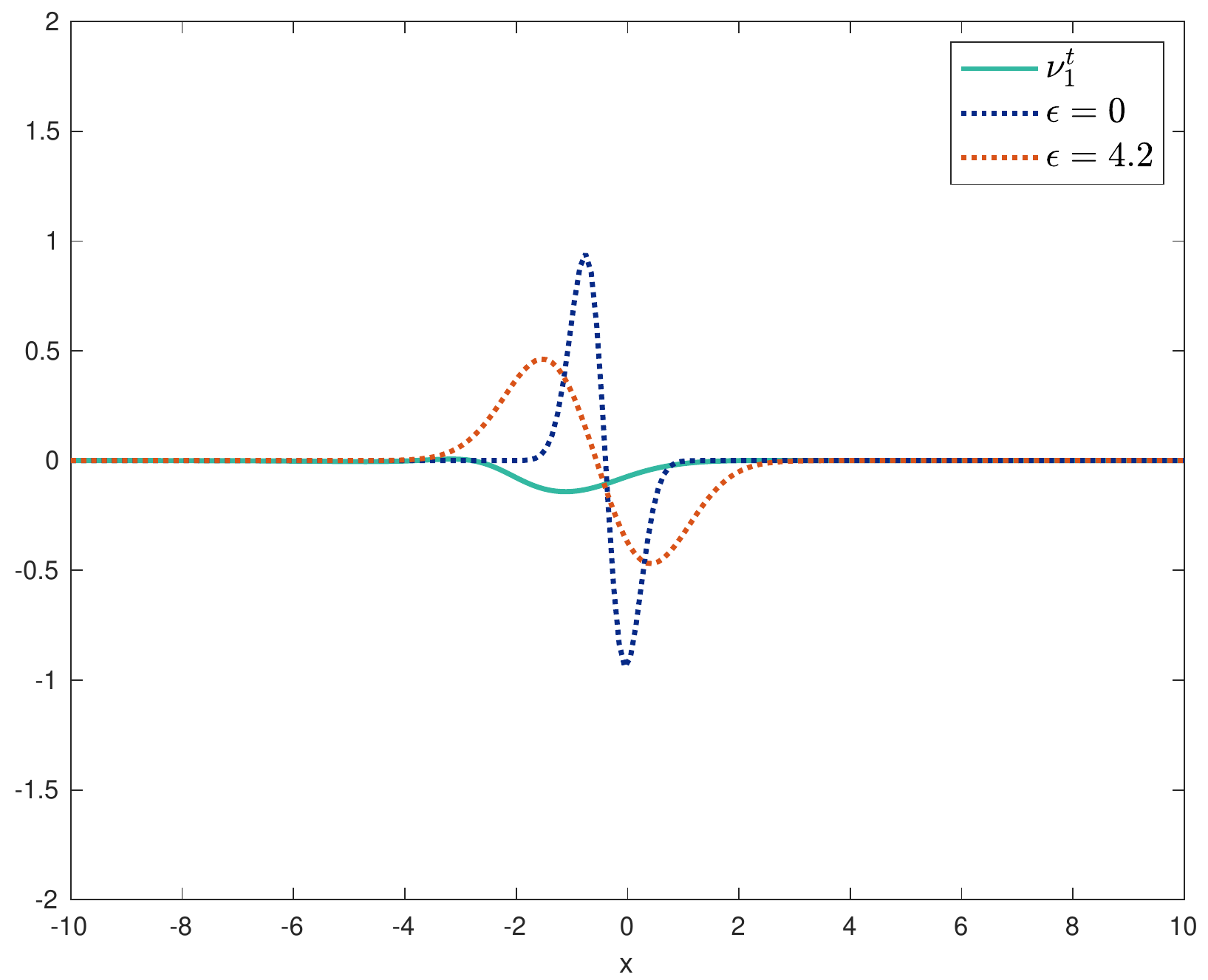}}
	\caption{Test 12: Screenshots of numerical simulation of \eqref{Euler1D_damping_press} for different values of $\varepsilon$. Comparison between zero-th order moment of the solution of (V) and $\mu^t$,  first order moment of the solution of (V) and $Q^t$  at a),b) $t=1$, c),d) $t=2$, e),f) $t=4$. Here $\eta=3$, $\alpha=0$.  }
	\label{fig:2bumps_press}
\end{figure}

\begin{figure}[h!]
	\centering
	\subfigure[]{\includegraphics[width=0.30\textwidth]{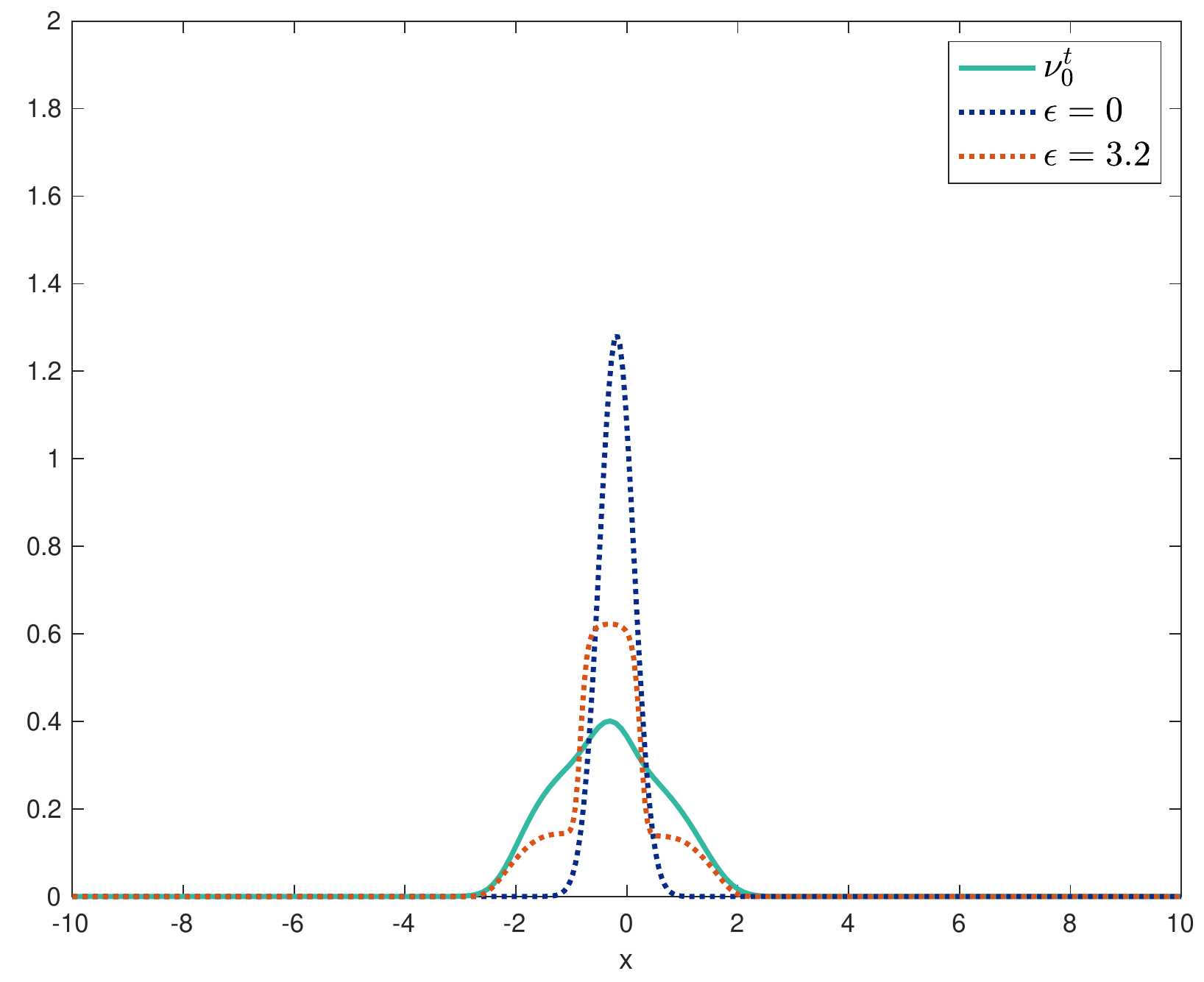}}
	\subfigure[]{\includegraphics[width=0.30\textwidth]{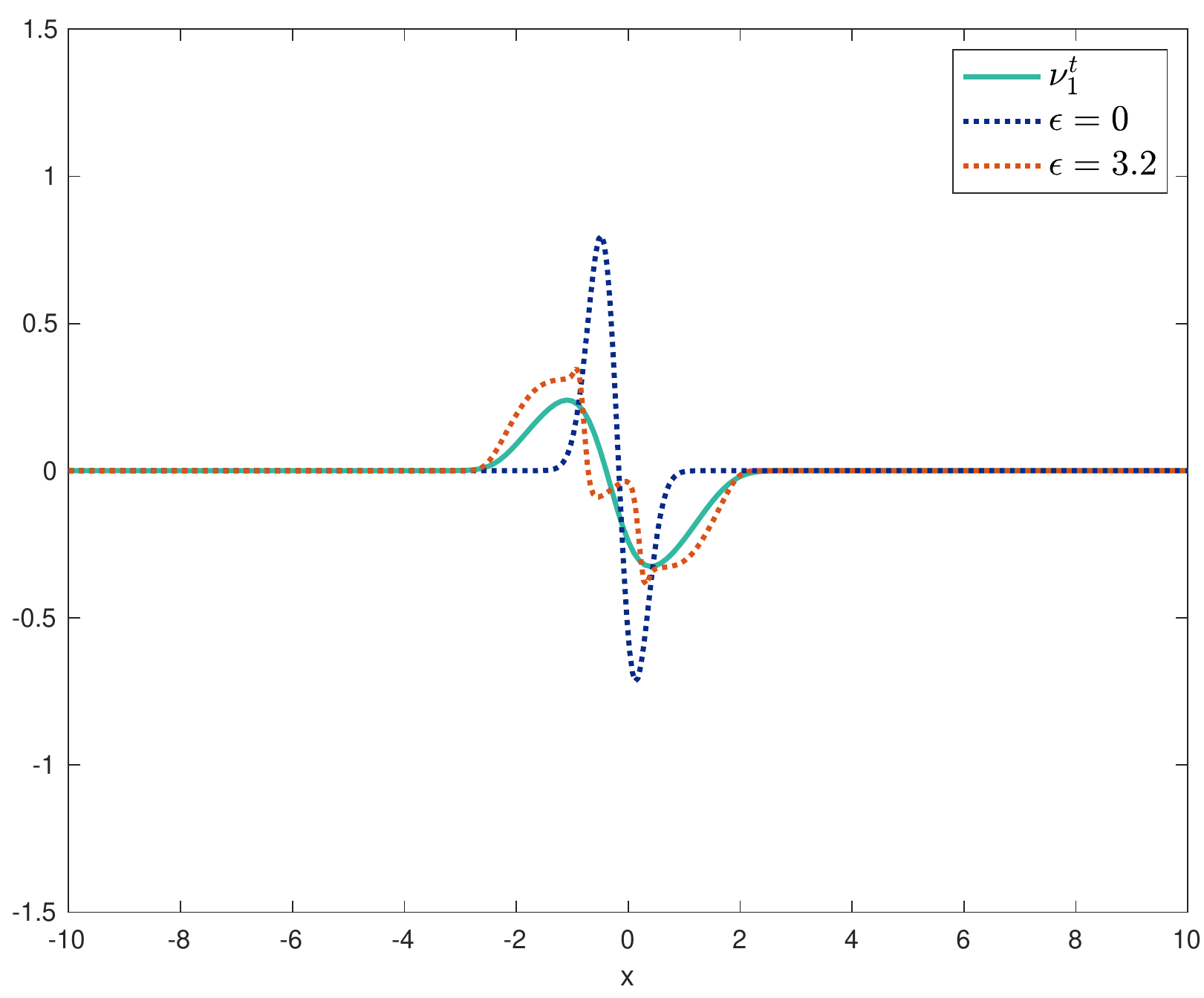}}\\
	\subfigure[]{\includegraphics[width=0.30\textwidth]{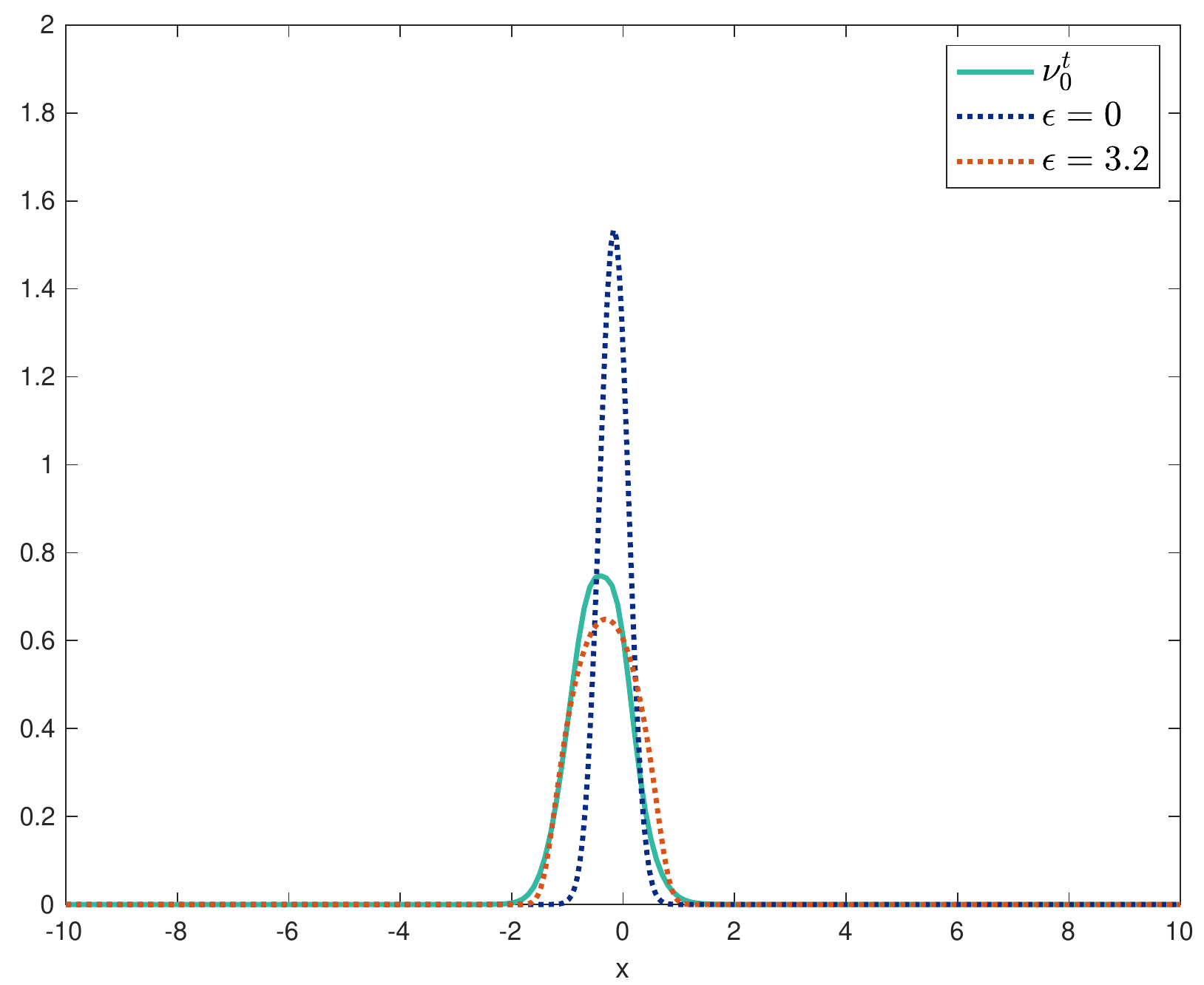}}
	\subfigure[]{\includegraphics[width=0.30\textwidth]{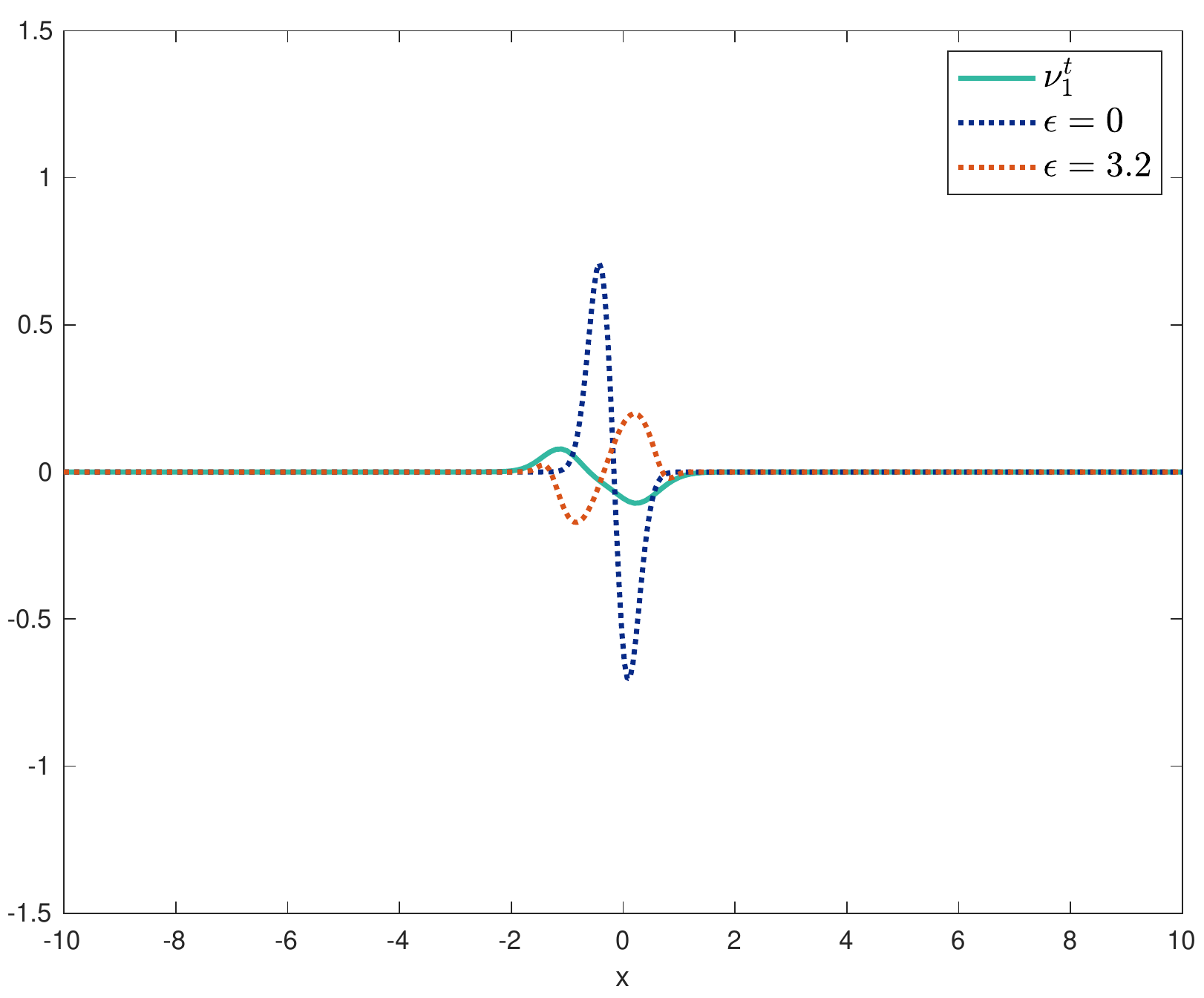}}\\
	\subfigure[]{\includegraphics[width=0.30\textwidth]{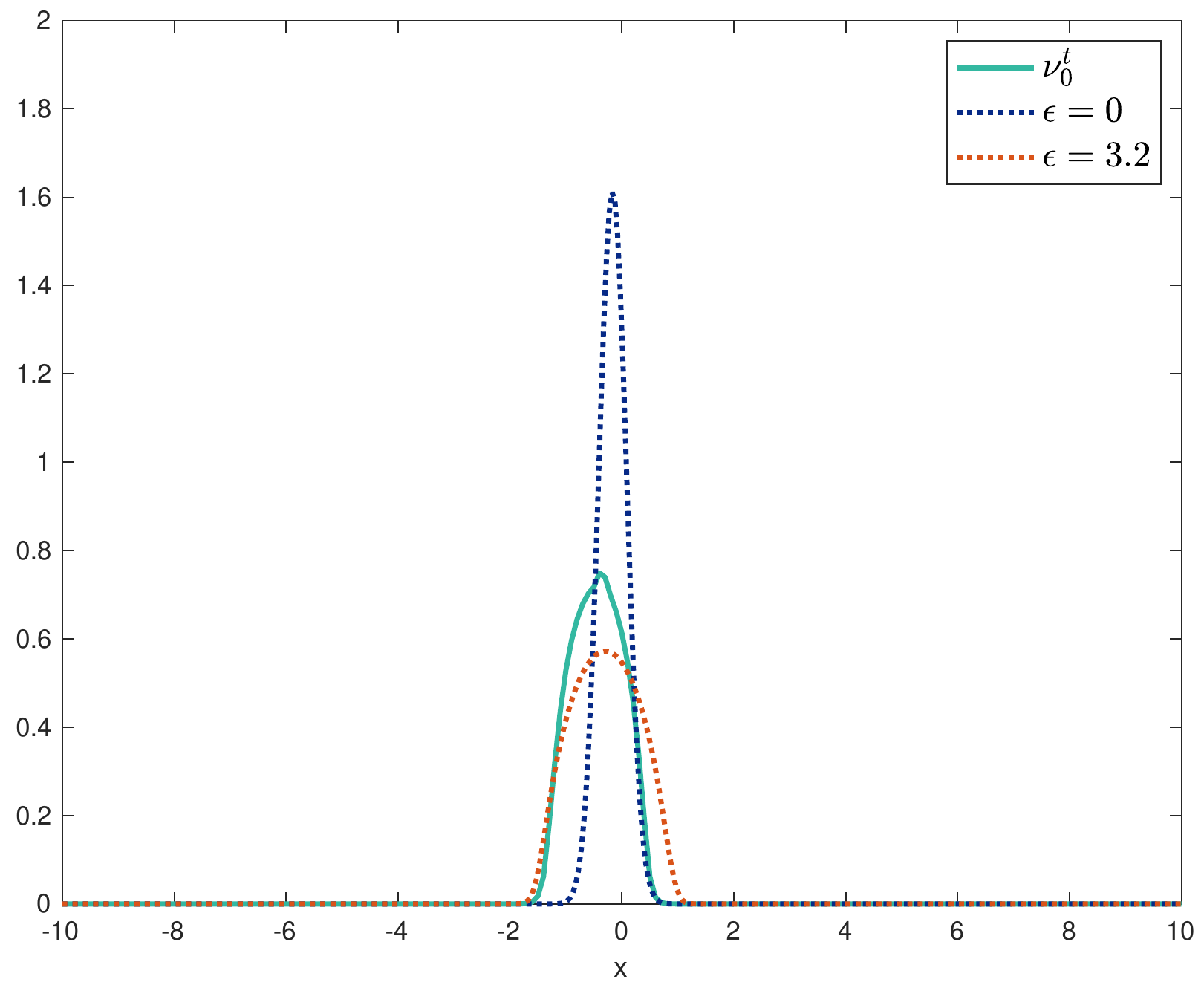}}
	\subfigure[]{\includegraphics[width=0.30\textwidth]{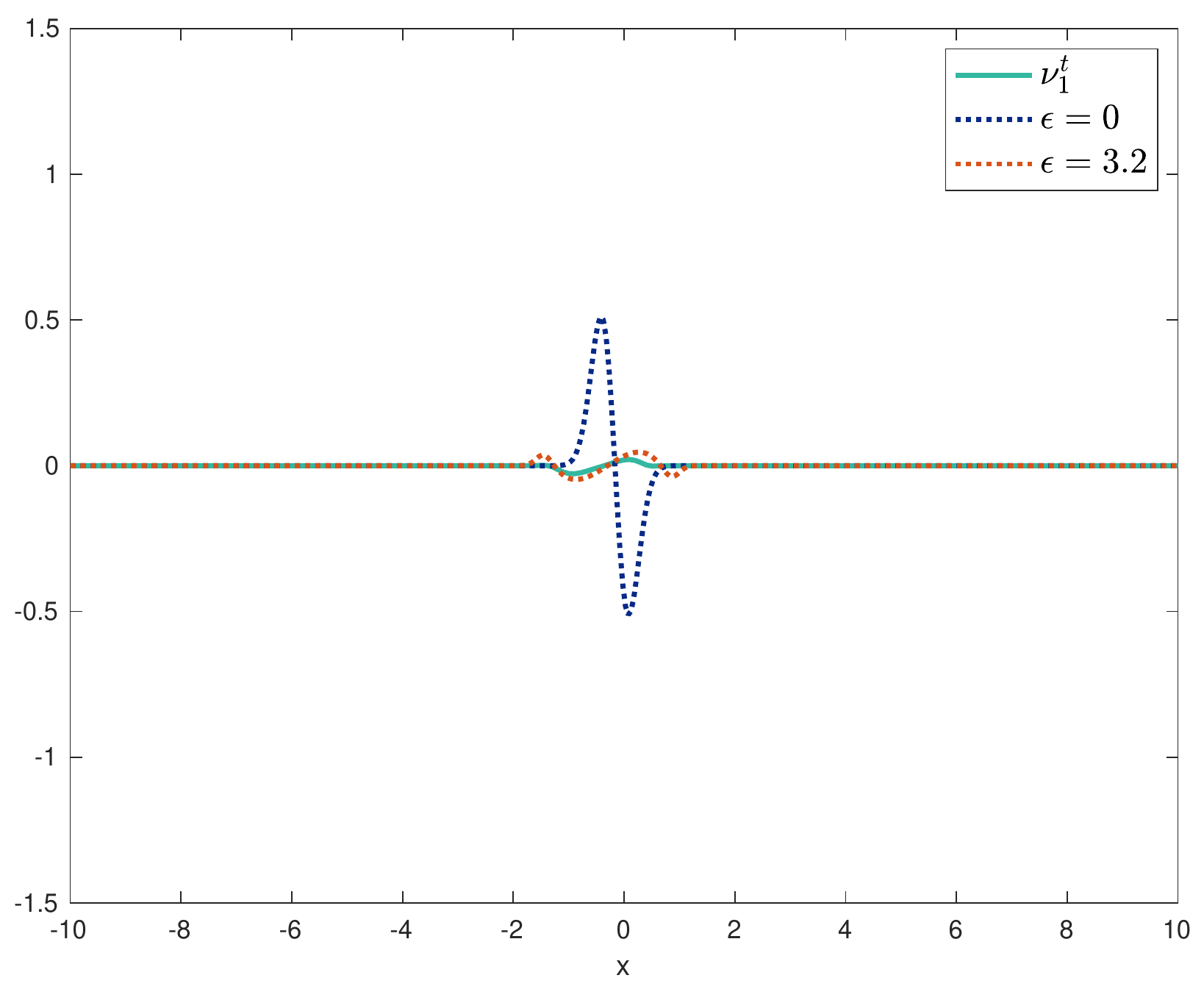}}
	\caption{Test 13: Screenshots of numerical simulation of \eqref{Euler1D_damping_press} for different values of $\varepsilon$. Comparison between zero-th order moment of the solution of (V) and $\mu^t$,  first order moment of the solution of (V) and $Q^t$  at a)-b) $t=1$, c)-d) $t=2$, e)-f) $t=4$. Here $\eta=3$, $\alpha=1$.  }
	\label{fig:2bumps_press_damping}
\end{figure}

Figure \ref{fig:2bumps_press} and Figure \ref{fig:2bumps_press_damping} show three snapshots of the performed numerical simulations without and with damping term, respectively.
In particular we plot $\nu_0$ and $\nu_1$ against the solution of \eqref{Euler1D_damping_press} for two values of $\varepsilon$, corresponding to the case of pressureless Euler system ($\varepsilon=0$) and choosing the value of the parameter computed as the mean value of the optimal values $\varepsilon^*_k$ collected at each iteration.
For the sake of completeness, in Table \ref{tab:2bumps_nodamping} and \REVISION{Table} \ref{tab:2bumps_damping} we present the results obtained for the three sample time instants chosen in the corresponding Figure \ref{fig:2bumps_press} and Figure \ref{fig:2bumps_press_damping}.

\begin{table}[h!]
		\begin{center}
	\begin{tabular}{|l|c|c|c|c|c|c|c|c|c|}
		\hline
		 &
		$\varepsilon^*_k$&
		$E_0^k(\varepsilon^*_k)$ 
		\\ \noalign{\hrule height2pt}
		$t_k=0$  & 0  & 0     \\ \hline 
  		$t_k=1$  & 4.146  & 0.034       \\ \hline 
                 $t_k=2$  & 4.235  & 0.533        \\ \hline 
        		$t_k=4$  & 4.236  & 0.588       \\  
		\hline
		
	\end{tabular}

		\end{center}
\caption{Computation of $\epsilon^*_k$ and corresponding value of $E^k_0$. Time instants correspond to the ones selected for snapshots in Figure \ref{fig:2bumps_press}.   }
\label{tab:2bumps_nodamping}
\end{table}

\begin{table}[h!]
		\begin{center}
	\begin{tabular}{|l|c|c|c|c|c|c|c|c|c|}
		\hline
		 &
		$\varepsilon^*_k$&
		$E_0^k(\varepsilon^*_k)$ 
		\\ \noalign{\hrule height2pt}
		$t_k=0$  & 0  & 0     \\ \hline 
  		$t_k=1$  & 3.891  & 0.154      \\ \hline 
                 $t_k=2$  & 3.082  & 0.094        \\ \hline 
        		$t_k=4$  & 2.639  & 0.178       \\ 
		\hline
	\end{tabular}
		\end{center}
\caption{Computation of $\epsilon^*_k$ and corresponding value of $E^k_0$. Time instants correspond to the ones selected for snapshots in Figure \ref{fig:2bumps_press_damping}.   }
\label{tab:2bumps_damping}
\end{table}

Starting with non monokinetic initial data, the agreement between Vlasov moments and Euler solutions, even adding pressure, is not so evident as in the approximated monokinetic case.
Focusing on the density profile, as shown in Figure \ref{fig:2bumps_press}, first column, choosing an optimal value for $\varepsilon$ helps in suppressing the blow-up occurring with $\varepsilon=0$.
At the beginning, the agreement obtained by adding the pressure term seems promising. As time grows, the blow-up is still reduced but the density and zero-moment profiles are no longer overlapped. In fact, due to nonzero initial momentum, the \REVISION{zero-th order} moment of Vlasov solution travels, moving leftwards.
The comparison for longer time improves introducing also a damping term, considering a nonnull value of $\alpha$ both in \eqref{VlasovChemo} and in \eqref{Euler1D_damping_press}. 
Actually, the pressure term affects only the Euler system, whereas the damping term involves a coherent dissipation both in the Vlasov and in the Euler dynamics. 
The obtained results are shown in Figure \ref{fig:2bumps_press_damping}.
The combined effect of the pressure term, reducing the density blow up, and of the damping term both at kinetic and macroscopic level, results in an improved agreement. Moreover, for longer time, also a good agreement between first moment and momentum is \REVISION{reached}, as shown in Figure \ref{fig:2bumps_press_damping} f).

\newpage

\newpage

\section{Conclusions}\label{compers}

In this paper we investigated numerically three different levels of descriptions of hybrid models for dynamics of interacting cells subject to chemotaxis: the microscopic position-velocity description $(P)$, the phase space density Vlasov type description $(V)$ and the physical space density-velocity field Euler type  description $(E)$.
At particle level, in the last years several analytical and numerical results have been presented \cite{dicostanzo2, mp, mp2, mppsJDE}, whereas kinetic and corresponding hydrodynamic limit have been analytically studied in \cite{rt, rt2}.
In particular, a complete correspondence between Vlasov moments and Euler solutions have been proved in \cite{rt2} under monokinetic assumption, which is a strong requirement in applications.
In our work we start investigating the behavior of the solutions far from the monokinetic case, obtaining numerical insights that can be seen as a first step for a rigorous study in general contexts.
Although each passage between the three levels consists in loosing some information - averaging on every but one particle and letting the number of particles diverge for $(P)\to(V)$ and taking only \REVISION{zero-th} and first moments in velocity for $(V)\to(E)$ -  we show that in many situations the most striking features of the particle level $(P)$ are preserved by passing from the paradigm of $N$ particles to the one of Vlasov at $N\to\infty$. 
Concerning the passage $(V)\to(E)$, numerical evidence shows that the fidelity of Euler versus Vlasov is increased by the presence of the chemical interaction and damping term, leading to positive results even for situations which are not covered by rigorous analytical results.
With respect to previous Euler systems coupled with chemotactic effect, starting with the seminal work in \cite{GambaPercolation}, the Euler system $(E)$ keeps memory of the interactions at the microscopic level through the nonlocal integral term. This is an important feature in biological applications, especially when dealing with experimental data.
However,  two main drawbacks need to be considered. First, it is well-known \REVISION{that} pressureless Euler systems exhibit finite time blow-up of the solution, depending on the initial data (\cite{CarrilloEuler}, Theorem 3.1). Moreover, from a numerical point of view, the approximation of the integral term is computationally costly and need to be improved, aiming at 2D and 3D simulations.
In order to control blow up phenomena, and to reduce the effect of neglected momenta, we \REVISION{investigated the effect of} an additional pressure term, still keeping the nonlocal integral term. 
In particular, we let the effect of the pressure term variable through the parameter $\varepsilon$. Tuning its value we investigate the existence of an optimal value $\varepsilon^*$ realizing an improved correspondence between Vlasov moments and Euler solutions.
From a numerical point of view, the transition $(V)\to(E)$ realizes an economical numerical gain, passing from a $2d$ PDE to a system of two $1d$ PDEs.
Establishing correlations among the different scales of descriptions allows to improve the computational effort, and it is one of the main perspective for applications to biological systems, aiming at developing the most appropriate macroscopic model as their digital twins.

\begin{appendix}

\section{Numerical schemes}\label{schemes}

In the following we detail the numerical schemes used for the numerical simulations of the different equations.

\textbf{Particle level: numerical scheme for (P)} 
For the ODE part of the coupled system,  we implement the well-known implicit Euler method, whereas the chemotaxis equation is discretized computing at each time the value of the solution on grid point discretizing the spatial domain.
We consider the simulation time interval $[0,T]$, for some $T>0$, divided into subintervals of same size $\Delta t$, assuming that $T$ is a multiple of $\Delta t$ to avoid rounding. We define $n_T:=\frac{T}{\Delta t}\in\mathbb N$, and denote the time grid points as $\{t^0,\ldots,t^{n_T}\}$.
Finally, $(x_i^k, v_i^k)$ denote approximations of position and velocity of agent $i$ at time $t^k$, for any $i=1,...,N$, $k=1,...,n_T$ starting with initial condition $(x_i^0, v_i^0)$.

Focusing on the chemotaxis equation, we eliminate the stiff degradation term performing the classical exponential transformation $\varphi^t(x)=e^{-\kappa s} u^t(x)$, where $u$ solves
\begin{equation}
\partial_t u=D\partial_{x}^{2}u+ e^{\kappa t}  \sum\limits_{j=1}^N\chi (x-x_j),
\end{equation}
We approximate the solution by means of a second order centered finite differences method in space and a classical explicit-implicit Crank-Nicholson integration in time, assuming periodic boundary conditions.
The computational domain $\Omega \times [0,T]$, for some $\Omega\subset \mathbb R$, is divided in cells of side $\Delta x \times \Delta t$, where $\Delta x$ denotes the space step, and the time step is the same used for the ODE part.  We denote with $n_{\Omega}$ the total number of the spatial grid points $\left\{\omega_0,\ldots,\omega_{n_{\Omega}}\right\}$.
In greater details, denoting with $u_l^k$ the approximation of $u$ at the grid point $(y_l,t^k)$, the numerical scheme reads
\begin{equation}\label{approximated_parabolic}
\begin{aligned} 
\frac{u^{k+1}_l-u^k_l}{\Delta t} =  D\frac{u^k_{l-1}-2u^k_{l} +u^k_{l+1}}{(\Delta x)^2}
&+ \frac{1}{2} e^{\kappa (k+1) \Delta t}  \sum\limits_{j=1}^N\chi (y_l-x_j^{k+1})\\
&+\frac{1}{2} e^{\kappa k \Delta t}  \sum\limits_{j=1}^N\chi (y_l-x_j^{k}) 
\end{aligned} 
\end{equation}
The coupling between the particle dynamics and the chemotaxis equation is realized by means of the gradient (the spatial derivative) of $\varphi$, evaluated, at each time instant $t_k$, in the position of agent $i$ at that time. 
The values of the gradient at grid points is approximated with centered difference, i.e. $\partial_l \varphi^k \approx \frac{\varphi^k_{l+1}-\varphi^k_{l-1}}{2 \Delta x}$. 
Since $x_i^k$ does not necessarely correspond to one of the grid point $y_l$, we approximate $\partial_x \varphi$ at time $t_k$ with a linear interpolation of the values of the gradient at the nearest grid points, in the following denoted as $\Phi_{x_i^k}(\partial_l \varphi^k)$.

Recalling the choice of $\gamma$ performed in \eqref{gamma_particle}, choosing $R=1$, the numerical scheme implemented for $(P)$ reads

\begin{align}\label{P_discretized}
\left\{
	\begin{array}{lll}
	v_i^{k+1}=v_i^{k}+\frac{\Delta t}{N} \sum\limits_{j=1}^N \frac{v_j^{k+1}-v_i^{k+1}}{\left(1+\left\|x^k_i-x^k_j\right\|^2\right)^\beta}+\eta \Delta t \Phi_{x_i^k}(\partial_l \varphi^k)\\\\
         x_i^{k+1}=x_i^{k}+\Delta t v_i^{k+1} \\\\
         \varphi_l^{k} = e^{-\kappa k \Delta t} u_l^{k} \ \ \mbox{with}\ \ \ \ u_l^{k} \ \ \mbox{approximated as in \eqref{approximated_parabolic}}
		\end{array}
	\right.
\end{align}

\textbf{Kinetic level: numerical scheme for (V)} %
Let now focus on system (V). 
The computational domain $\Omega \times [0,T ]$ is divided in cells of side $\Delta x \times \Delta v \times \Delta t$ where
$\Delta x$, $\Delta v$, are the space and velocity steps and $\Delta t$ is the time step. Let us assume that $T$ is a multiple of $\Delta t$ to avoid rounding. 

We choose $\Omega$ as a rectangular domain of size $n_x\Delta x \times n_v\Delta v$, with $n_x, n_v \in\mathbb N$ the number of inner point discretizing the space and velocity domain. The generic point of the 2-dimensional grid is hence denoted as $(x_i,v_j)$
with $i=0,...,n_x$, $j=0,...,n_v$.
We consider periodic boundary condition with respect to the $x$ variable, and null condition with respect to the $v$ variable.
For any time instant $t_k$ we denote the approximation of $\rho$ and $\nu$ at $(x_i,v_j)$ with $\rho_{i,j}^{k}$ and $\mathcal{V}^k_{i,j}$, respectively. 
The first-order scheme implemented reads:

\begin{align}
\begin{array}{lll}
\displaystyle\frac{\rho_{i,j}^{k+1}-\rho_{i,j}^{k}}{\Delta t}&=
 \displaystyle-\frac{v_j}{2 \Delta x} \left( \rho_{i+1,j}^{k}-\rho_{i-1,j}^{k}\right)
  \displaystyle+\frac{\left|v_j\right|}{2 \Delta x} \left( \rho_{i+1,j}^{k}-2\rho_{i,j}^{k}+\rho_{i-1,j}^{k}\right)\\\\
  &  \displaystyle-\frac{1}{2 \Delta v} \left( \rho_{i+1,j}^{k}\mathcal{V}^k_{i+1,j}-\rho_{i-1,j}^{k}\mathcal{V}^k_{i-1,j}\right)\\\\
  &+ \displaystyle \frac{1}{2 \Delta v} 
  \left( \rho_{i+1,j}^{k} \left|\mathcal{V}^k_{i+1,j} \right| -2\rho_{i,j}^{k} \left|\mathcal{V}^k_{i,j} \right|+\rho_{i-1,j}^{k}\left|\mathcal{V}^k_{i-1,j} \right|\right)
 \end{array}
\end{align}

At each time step, a CFL condition of the form

\begin{equation}
1- \frac{\Delta t}{\Delta x} |v_{max}| - \frac{\Delta t}{\Delta v} \max_{i,j} | \mathcal{V}^k_{i,j} | \ge 0
\end{equation}
 must be satisfied. 

As at particle level, we approximate the solution of the chemotaxis equation by means of a second order centered finite differences method in space and a classical explicit-implicit Crank-Nicholson integration in time. 
The only difference is given by the source term. We approximate the two-dimensional integral in \eqref{VlasovChemo} using a  quadrature formula. In greater details, it holds:
\begin{equation}
\displaystyle \int_{\mathbb{R}} \int_{x-R}^{x+R}\rho^s(y,\xi)dy d\xi \approx \sum_{\substack{m,n:\\  x_m \in [x_i-R \Delta x, x_i+R \Delta x ],\\  n=0,...,n_v}} \rho^k_{m,n} \Delta x \Delta v
\end{equation}

Concerning the coupling through the gradient, the approximation of $\partial_x \psi$ is required at spatial grid point. We approximate its value at each time instant using centered finite difference
$\partial_i \psi^k \approx \frac{\psi_{i+1}^k-\psi_{i-1}^k}{2\Delta x}$.

\textbf{Macroscopic level: numerical scheme for (E)}
We describe the numerical scheme for the most general Euler system simulated in this paper \eqref{Euler1D_damping_press}. The scheme used for system \eqref{Euler1D_nochemo} and \eqref{Euler1D_damping} can be obtained simply considering $\varepsilon=0, \eta=0, \alpha=0$ and $\varepsilon=0$, respectively.

Denoting
$W=W(x,t)=(\mu^t(x),\mu^t(x) u^t(x))^T$ the vector of unknowns, density and momentum, we rewrite the hyperbolic part as
\begin{equation}\label{sis_relaxation}
W_t + A(W)_x=S(W),
\end{equation}
where A and S are defined as follow:
 \begin{equation} 
 A(W) =
      \begin{pmatrix}
       \mu u   \\
       \mu u^2+\varepsilon P(\mu)      
      \end{pmatrix}
    \end{equation}
    \begin{equation} 
 S(W) =
      \begin{pmatrix}
       0   \\
       \mu\int\gamma(\cdot-y,u(\cdot)-u(y))\mu(y)dy
       +\eta\mu\partial_{x}\psi-\alpha \mu u    
      \end{pmatrix}
    \end{equation}

We approximate the solution at each time step starting from the relaxation tecniques originally proposed in \cite{AregbaNatalini}.
In few words, the idea is to approximate hyperbolic system of equations through a BGK relaxation, which leads to a linear advection system with a relaxation source term. 
The advantage lies in the fact that at each time step $t^k$, in order to find the solution at $t^{k+1}$, it is required to solve a system of linear, transport problems, which reduces the
computational complexity.

We denote with $I \subset \mathbb{R}$ the computational domain. We approximate the solution of Euler system $W_i^k=(\mu_i^k,\mu_i^k u_i^k)^T$ and of the parabolic equation $\psi_i^k$ at each point of the grid $x_i=i \Delta x$, $i=1,...,n_x-1$. At boundary, we assume $W_{0}=W_{n_x}=(0,0)^T$.
The parabolic equation is approximated using the classical explicit-implicit Crank- Nicholson method, as previously explained.
The algorithm requires to first solve a homogeneous
system of equations and then a time integration of the source part is performed.
Let now focus on the homogeneous part of \eqref{sis_relaxation}. 
Using BGK approximation with two velocities of equal speed and opposite sign, i.e. $\lambda_1= \lambda=-\lambda_2$ we denote with $\textbf{f}_r=(f_{r}^{\mu}, f_{r}^{\mu u})$, $r=1,2$,  the components of BGK model approximating $W$. For a complete description of the relaxation technique see \cite{shixin}.
 At each time iteration we solve

 \begin{equation} \label{sis_f_omogeneo}
 \partial_t
      \begin{pmatrix}
      \textbf{f}_1  \\
       \textbf{f}_2     
      \end{pmatrix}
    +
      \begin{pmatrix}
      \lambda & 0 & 0 & 0  \\
       0 & \lambda & 0 & 0  \\
       0 & 0 & -\lambda& 0  \\
         0 & 0 & 0 & -\lambda  
      \end{pmatrix}
      \partial_x
      \begin{pmatrix}
      \textbf{f}_1  \\
       \textbf{f}_2     
      \end{pmatrix}
      = \textbf{0}
    \end{equation}    
with initial states given by the Maxwellian functions

 \begin{equation} 
      \begin{pmatrix}
      \textbf{f}_1^0  \\
       \textbf{f}_2^0     
      \end{pmatrix}
      =
       \begin{pmatrix}
     \frac{1}{2}\left( W(x,0) + \frac{A(W(x,0))}{\lambda}\right)  \\
      \frac{1}{2}\left( W(x,0) - \frac{A(W(x,0))}{\lambda}\right) 
      \end{pmatrix}
 \end{equation}   
Denoting with $(\textbf{f}_r)_i^{k}$, the approximated value of $\textbf{f}_r$, at time $t_k$ at grid point $x_i$, we compute $(\textbf{f}_r)_i^{k+1/2}$
running the following numerical upwind scheme:
 $$
( \textbf{f}_r)_{i}^{k+1/2}= (\textbf{f}_r)_{i}^{k}
-\frac{\Delta t}{2\Delta x} \lambda_r^k\left( (\textbf{f}_r)_{i+1}^{k}-(\textbf{f}_r)_{i-1}^{k} \right)
+ \frac{\Delta t}{2\Delta x} |\lambda_r^k| \left( (\textbf{f}_r)_{i+1}^{k}-2(\textbf{f}_r)_{i}^{k}+(\textbf{f}_r)_{i-1}^{k} \right)
 $$
for any $i=1,...,n_x$, $r=1,2$, with 
$\lambda_1^k=\max_{i} \left( |u^k_i| + \sqrt{P^\prime (\mu^k_i)}\right)$, 
$\lambda_2^k=-\lambda_1^k$.
The stability of the scheme is guaranteed choosing, at each time step, $\Delta t ^k=0.9 \frac{\Delta x}{\lambda_1^k}$.

The solution $W$ at time $t_{k+1}$ is obtained adding the discretization of the source term $S(W)$.
The integral term and the chemotaxis gradient are treated explicitly, whereas the damping term in implicit, due to stiffness problem.
The approximated value of the integral term, in the following denoted as $\mathcal{I}_i^k$ , is obtained by a quadrature formula, which has to be solved at each iteration.
Finally we get
\begin{equation}
W_i^{k+1}=( \textbf{f}_1)_{i}^{k+1/2}+( \textbf{f}_2)_{i}^{k+1/2}+ \Delta t \left( 0,  \mu_i^k \mathcal{I}_i^k+\eta \mu_i^k \frac{\psi_{i+1}^k- \psi_{i-1}^k}{2\Delta x} \right).
\end{equation}

\end{appendix}


\begin{thebibliography}{99}

\bibitem{AlbiPareschi}
G. Albi, and L. Pareschi. 
Binary interaction algorithms for the simulation of flocking and swarming dynamics. Multiscale Modeling \& Simulation, 11(1), 1-29 (2013).

\bibitem{AlbiControllo}
G. Albi, C. Young-Pil, A.-S. Haeck. 
Pressureless Euler alignment system with control.
Mathematical Models and Methods in Applied Sciences, 28(09), 1635-1664 (2018).

\bibitem{AregbaNatalini}
 D. Aregba-Driollet, and R. Natalini. 
 Discrete kinetic schemes for multidimensional systems of conservation laws.
 SIAM Journal on Numerical Analysis, 37(6), 1973-2004 (2000).

\bibitem{arboleda}
 Y.~Arboleda-Estudillo, M.~Krieg, J.~St\"{u}hmer, N.~A. Licata, D.~J.
  Muller and C.-P. Heisenberg.
{M}ovement {D}irectionality in {C}ollective {M}igration of {G}erm
  {L}ayer {P}rogenitors,
Current Biology, 20(2), 161-169 (2010).

\bibitem{belmonte}
J.~M. Belmonte, G.~L. Thomas, L.~G. Brunnet, R.~M. de~Almeida and
  H.~Chat\'{e}.
 Self-propelled particle model for cell-sorting phenomena.
 Physical Review Letters 100, 248702 (2008).

\bibitem{sticky2}
Y. Brenier, W. Gangbo, G. Savaré, M. Westdickenberg.
Sticky particle dynamics with interactions.
Journal de Mathématiques Pures et Appliquées,
99(5), 577-617 (2013).

\bibitem{sticky1}
Y. Brenier, and E. Grenier.
Sticky Particles and Scalar Conservation Laws.
SIAM Journal on Numerical Analysis 35(6), 2317-2328 (1998).

\bibitem{axioms}
G. Bretti, A. De Ninno , R. Natalini, D. Peri, N. Roselli.
 Estimation algorithm for a hybrid pde–ode model inspired by
 immunocompetent cancer on-chip experiment,  Axioms, 10(4), 243 (2021).

\bibitem{bouchut}
F. Bouchut.
On zero pressure gas dynamics.
Advances in kinetic theory and computing: selected papers, 171-190 (1994).

\bibitem{carr4} J. A. Carrillo,  Y.-P. Choi.
J. A. Carrillo, and Y.-P. Choi,
Mean-field limits: from particle descriptions to macroscopic equations.
 Archive for Rational Mechanics and Analysis 241, 1529-1573 (2021).

\bibitem{CarrilloVlasov}
J. A. Carrillo, M. Fonasier, G. Toscani, F. Vecil.
 \newblock Particle, kinetic, and hydrodynamic models of swarming
\newblock Mathematical modeling of collective behavior in socio-economic and life sciences (pp. 297-336). Birkhäuser Boston. 2020.

\bibitem{CarrilloVecil}
J.A. Carrillo, and F. Vecil. 
 Nonoscillatory interpolation methods applied to Vlasov-based models.
  SIAM Journal on Scientific Computing 29(3), 1179-1206 (2007).

\bibitem{CarrilloEuler}
J. A. Carrillo, Y.P. Choi, and S. Perez.
A review on attractive–repulsive hydrodynamics for consensus in collective behavior.
Active Particles, Volume 1: Advances in Theory, Models, and Applications, 259-298 (2017).

\bibitem{CarrilloEuler2}
J. A. Carrillo, Y.-P. Choi, E. Tadmor, and C. Tan.
Critical thresholds in 1D Euler equations with nonlocal forces.
 Mathematical Models and Methods in Applied Sciences, 26(1), 185-206 (2016).
 

\bibitem{CS1}
F. Cucker, and S. Smale.
 On the mathematics of emergence. 
 Journal of Mathematics, 2, 197-227 (2007). 
 
\bibitem{CS2} F. Cucker and S. Smale.
Emergent behavior in flocks. 
IEEE Transaction on Automatic Control, 52(5), 852-862 (2007).





\bibitem{DiMarco}
G. Di Marco and L. Pareschi.
Numerical methods for kinetic equations,
Acta Numerica 23, 369-520 (2014).

\bibitem{dorsogna}
 M.~R. D'Orsogna, Y.~L. Chuang, A.~L. Bertozzi and L.~S. Chayes.
 {S}elf-{P}ropelled {P}articles with {S}oft-{C}ore {I}nteractions:
  {P}atterns, {S}tability, and {C}ollapse,
Physical review letters, 96(10), 104302 (2006).

\bibitem{dicostanzo} E. Di Costanzo, R. Natalini, Roberto, L. Preziosi.
A hybrid mathematical model for self-organizing cell migration in the zebrafish lateral line. 
Journal of mathematical biology 71(1), 171-214 (2015).


\bibitem{dicostanzo2}
 E. Di Costanzo, M. Menci, E. Messina, R. Natalini, A. Vecchio.
 A hybrid model of collective motion of discrete particles under alignment and continuum chemotaxis, Discrete \& Continuous Dynamical System-Series B, 25(1)  443-472 (2020).

\bibitem{dicostanzo_cardio}
E. Di Costanzo,  A. Giacomello,  E. Messina,  R. Natalini, G. Pontrelli,  F. Rossi,  R. Smits,  M. Twarogowska.
 A discrete in continuous mathematical model of cardiac progenitor cells formation and growth as spheroid clusters (Cardiospheres), 
 Mathematical Medicine and Biology: A Journal of the IMA, 35(1), 121-144 (2018).
 
 \bibitem{Elchina}
 N. V. Elkina and J. Buchner.
 A new conservative unsplit method for the solution of the Vlasov equation, Journal of Computational Physics, 213(2), 862-875 (2006).

%

\bibitem{Filbet1}
F. Filbet, E. Sonnendrücker,
Comparison of eulerian vlasov solvers,
Computer Physics Communications, 150(3), 247-266 (2003).


\bibitem{GambaPercolation} 
A. Gamba, D. Ambrosi, A. Coniglio, A. de Candia, S. Di Talia, E. Giraudo, G. Serini, L. Preziosi, and F. Bussolino. 
Percolation, morphogenesis, and burgers dynamics in blood vessels formation. 
Physical Review Letters, 90(11):118101 (2003).

\bibitem{gilmour06} P. Haas P, D. Gilmour.
  Chemokine signaling mediates self-organizing tissue migration in the zebrafish
lateral line. 
Developmental cell, 10(5), 673–680 (2006).

  
\bibitem{haliu} S.-Y. Ha and J.-G. Liu.
 A simple proof of Cucker-Smale flocking dynamics and mean-field limit.
Communications in Mathematical Sciences, 7(2), 297-325 (2009).
 

\bibitem{hatzikirou}
 H.~Hatzikirou and A.~Deutsch.
Cellular automata as microscopic models of cell migration in heterogeneous environments,
Current topics in developmental biology, 81, 401-434 (2008).


%
  

\bibitem{shixin}
S. Jin and Z. Xin.
The relaxation schemes for systems of conservation laws in arbitrary space dimensions.
Communications on pure and applied mathematics 48(3), 235-276 (1995).

\bibitem{gilmour08} V. Lecaudey, G.C. Akdogan, WHJ Norton and D.Gilmour.   Dynamic Fgf signaling couples morphogenesis
and migration in the zebrafish lateral line primordium. 
Development 13, 2695–2705 (2008).

    
 
  
  \bibitem{vicsek-collective-cell}
 E.~M\'{e}hes and T.~Vicsek.
 Collective motion of cells: from experiments to models.
Integrative biology 6(9), 831--854 (2014).

\bibitem{mp} M. Menci and M. Papi.
 Global solutions for a path-dependent hybrid system of differential equations under parabolic signal.
 Nonlinear Analysis 184(7), 172–192 (2019).

\bibitem{mp2}
M. Menci and M. Papi. 
Existence of solutions for hybrid systems of differential equations under exogenous information with discontinuous source term.
 Nonlinear Analysis 221(8), 112885 (2022).

\bibitem{mppsJDE}
M. Menci, M. Papi, M. M. Porzio and F. Smarrazzo.
 On a coupled hybrid system of nonlinear differential equations with a nonlocal concentration. Journal of Differential Equations, 361, 288-338 (2023).



\bibitem{rt}  R. Natalini and T. Paul.
 On The Mean Field limit for Cucker-Smale models.
  Discrete \& Continuous Dynamical Systems - Series B, 27(5), 2873-2889 (2022).
  

\bibitem{rt2} R. Natalini and T. Paul.
 The Mean-Field limit for  hybrid  models of collective motions with chemotaxis.
SIAM Journal on Mathematical Analysis, 55(2), 900-928 (2023).


\bibitem{prt} B. Piccoli, F. Rossi, E. Tr\'elat, 
Control to flocking of the kinetic Cucker-Smale model.
SIAM Journal on Mathematical Analysis, 47(6), 4685–4719 (2014).



\bibitem{strombom}
 D.~Str\"{o}mbom,
Collective motion from local attraction,
Journal of theoretical biology, 283(1), 145-151 (2011).



\bibitem{tadmor2022} E. Tadmor.
Swarming: hydrodynamic alignment with pressure.
Bulletin of the American Mathematical Society, 60(3), 285-325 (2023).

\bibitem{vicsek}
T.~Vicsek, A.~Czir\`ok, E.~Ben-Jacob, I.~Cohen and O.~Shochet.
{N}ovel {T}ype of {P}hase {T}ransition in a {S}ystem of
  {S}elf-{D}riven {P}articles,
Physical Review Letters, 75(6),  1226--1229 (1995).

\bibitem{vicsek-collective}
 T.~Vicsek and A.~Zafeiris.
Collective motion,
Physics Reports, 517(3-4), 71--140 (2012).

\bibitem{VillaniAMS}
C. Villani.
Topics in Optimal Transportation, 
American Mathematical Society, 58 (2021).

\bibitem{VillaniTOT}
C. Villani
Optimal Transport: old and new.
Berlin: springer, 338 (2009).


\end{thebibliography}
\end{document}